\documentclass[10pt]{article}
\pdfoutput=1

\usepackage[unicode]{hyperref}
\usepackage{amsmath}
\usepackage{graphicx}								
\usepackage{subfig}
\usepackage{amsfonts} 							
\usepackage{dsfont}
\usepackage{bigints} 							
\usepackage{algorithm}              
\usepackage{algorithmic}
\usepackage{booktabs} 							
\usepackage{multirow}								
\usepackage{setspace}
\usepackage{mathrsfs}
\usepackage{amssymb}
\usepackage{tikz}
\usepackage{color}
\usepackage{xcolor}
\usepackage{mathtools}
\usepackage{xfrac}
\usepackage{upgreek}

\usepackage{framed}
\usepackage{afterpage}
\usepackage{capt-of}
\usepackage{stackengine}

\usepackage{geometry}
\usepackage{tikz}	
\usetikzlibrary{arrows,chains,matrix,positioning,scopes}

\makeatletter
\tikzset{join/.code=\tikzset{after node path={%
\ifx\tikzchainprevious\pgfutil@empty\else(\tikzchainprevious)%
edge[every join]#1(\tikzchaincurrent)\fi}}}

\makeatother
\tikzset{>=stealth',every on chain/.append style={join}, every join/.style={->}}
\tikzstyle{labeled}=[execute at begin node=$\scriptstyle,   execute at end node=$]
	
\newcommand {\Int}   {\int\limits}
\newcommand {\Sum}   {\sum\limits}	
\newcommand {\eps}   {\varepsilon}

\newcommand {\IntQT} {\Int_{Q_T}}

\newcommand {\CFriedrichs}    {C_{\mathrm{F}}}
\newcommand {\CF}    {C_{\mathrm{F\Omega}}}
\newcommand {\CPoincare}    {C_{\mathrm{P}}}

\newcommand {\Ctr}   {C_{\mathrm{tr}}}
\newcommand {\Rd}    {{\mathds{R}}^d}

\newcommand {\Ieff}  {{I_{\rm eff}}}

\newcommand {\bftau} {\boldsymbol{\tau}}

\def \NormA#1  {{\mid\!\mid\!\mid #1 \mid\!\mid\!\mid}^2_A }  
\def \NormAinverse#1  { {\mid\!\mid\!\mid #1 \mid\!\mid\!\mid}^2_{A^{-1}} }  
\def \Normenergy#1  {\mid\!\mid\!\mid #1 \mid\!\mid\!\mid}  

\def \NormQT#1 {{\mid\!\mid\!\mid #1 \mid\!\mid\!\mid}^2_{Q_T}}   
\def \NormQO#1 {{\mid\!\mid\!\mid #1 \mid\!\mid\!\mid}^2_{Q^0}}   
\def \NormQk#1 {{\mid\!\mid\!\mid #1 \mid\!\mid\!\mid}^2_{Q^k}}   
\def \Normf#1  {\Big \lceil #1 \Big\rceil_\Omega}
\def \dvrg     {\mathrm{div}}	

\def \dt       {\mathrm{\:d}t}
\def \dx       {\mathrm{\:d}x}

\def \dxt      {\mathrm{\:d}x\mathrm{d}t}
\def \dst      {\mathrm{\:d}s\mathrm{d}t}
\def \ds       {\mathrm{\:d}s}

\def \Normt#1  {\mid\!\mid\!\mid #1 \mid\!\mid\!\mid}

\def\L#1{L^{#1}}
\def\H#1{H^{#1}}
\def\HD#1#2{H^{#1}_{#2}}

\def\H#1{H^{#1}}

\def\maj#1{{\overline{\mathrm M}^2}^{#1}}

\def\Min{{\underline{\mathrm M}^2}}
\def\error{{[e]\,}^2}

\newcommand {\meas} {\mathrm{meas}}

\newcommand*\rfrac[2]{{}^{#1}\!/_{#2}}
\newcommand{\minus}{\scalebox{0.5}[1.0]{\( - \)}}

\def \Mean#1#2 {{ \Big \{ #1 \Big\} }_{#2}}
\def \smallMean#1#2 {{ \big \{  #1 \big\} }_{#2}}

\def\Pone{{\rm P}_{1}}
\def\Ptwo{{\rm P}_{2}}
\def\Pthree{{\rm P}_{3}}

\def\RTone{{\rm RT}_{1}}

\def\ProofBegin{\noindent{\bf Proof:} \:}
\def\ProofEnd{{\hfill $\square$}}

\newtheorem{theorem}{Theorem}{\bf}{\it}
{\bf}{\it}
{\bf}{\it}
\newtheorem{example}{Example}{\bf}{}
\newtheorem{remark}{Remark}{\bf}{\it}

\usepackage{color}
\definecolor{deepblue}{rgb}{0,0,0.5}
\definecolor{deepred}{rgb}{0.6,0,0}
\definecolor{deepgreen}{rgb}{0,0.5,0}
\usepackage{listings} 

\newcommand\pythonstyle{
\lstset{
		language=Python,
		basicstyle=\footnotesize\ttfamily,
		otherkeywords={self},          
		keywordstyle=\bfseries\color{deepblue}\bfseries,
		commentstyle=\itshape\color{purple!40!black},
		stringstyle=\color{deepgreen},
		frame=single,                         
		showstringspaces=false,            
		belowcaptionskip=.75\baselineskip
}}
\lstnewenvironment{python}[1][]
{
\pythonstyle
\lstset{#1}
}
{}

\newcommand\pythoninline[1]{{\pythonstyle\lstinline!#1!}}

\usepackage{colortbl}
\definecolor{gainsboro}{rgb}{0.86, 0.86, 0.86}
\definecolor{lightgray}{rgb}{0.83, 0.83, 0.83}
\definecolor{silver}{rgb}{0.75, 0.75, 0.75}
\definecolor{ashgrey}{rgb}{0.7, 0.75, 0.71}
\definecolor{battleshipgrey}{rgb}{0.52, 0.52, 0.51}

\def \EI           {{\rm E \!\!\! Id}}
\newcommand {\flux}     {\boldsymbol {y}}

\title{
On the a posteriori error analysis for linear Fokker-Planck
models in convection-dominated diffusion problems
}
\author{Svetlana Matculevich \thanks{RICAM Linz, Johann Radon Institute, Linz, Austria}
\and Monika Wolfmayr \thanks{Faculty of Information Technology, University of Jyv\"askyl\"a, Jyv\"askyl\"a, Finland}}
\date{}

\begin{document}
	
\maketitle 

\begin{abstract}	
This work is aimed at the derivation of reliable and efficient a posteriori error estimates for 
convection-dominated diffusion problems motivated by a linear Fokker-Planck problem appearing in 
computational neuroscience. We obtain computable error bounds of functional type for the static and 
time-dependent case and for different boundary conditions (mixed and pure Neumann boundary conditions).
Finally, we present a set of various numerical examples including discussions on mesh adaptivity and 
space-time discretisation. The numerical results confirm the reliability and efficiency of the error estimates 
derived.
\end{abstract}


\section{Introduction}

Convection-diffusion systems are widely used in natural sciences and financial mathematics for 
describing physical flows and stochastically changing systems, respectively. The Fokker-Planck 
(or Kolmogorov) equation is one of the examples of this class illustrating the velocity 
of a particle in the Brownian motion, along with the Black--Scholes equation governing the price 
evolution, and the Navier--Stokes equation that describes the motion of viscous fluid substances. 
There are many more problems appearing in physics or biology that are concerned with modelling 
of semiconductors or of decision-making processes in neuroscience that can be targeted by this type 
of systems. The considered equation is sometimes also called the drift-diffusion equation. 

This paper derives and discusses a posteriori error estimates of the functional type 
for convection-diffusion problems, which are motivated by the decision-making Fokker-Planck model 
problem discussed by Carrillo et al. \cite{CarrilloCordierMancini2011} appearing in computational 
neuroscience. 
The mathematical model discussed in \cite{CarrilloCordierMancini2011} is 
associated with stochastic dynamical systems, modelling the evolution of the decision-making 
(average firing rates) of two interacting neurone populations, see also \cite{BrownHolmes2001} and
\cite{DecoMarti2007}. In \cite{CarrilloCordierMancini2011}, the existence of a unique solution for stationary 
as well as evolutionary linear Fokker-Planck equations is discussed and proved under the 
assumption that the (regular enough) flux or drift is incoming in the bounded domain.
The proofs are mainly based on the work by Wloka \cite{Wloka1987} as well as  Evans
\cite{Evans2010}. In \cite{Wloka1987}, a more general result is proved, which can be then applied 
in order to show existence and uniqueness in the time-dependent case too. 
Carrillo et al. \cite{CarrilloCordierMancini2011} discusses this issue together with 
the large time behaviour of the problem and proves for that the convergence of the solution 
to the stationary state.

The main contribution of this work is the derivation of two-sided functional 
a posteriori error estimates (the so-called majorant and minorant) for stationary as well as 
time-dependent linear convection-diffusion problems. The ideas of functional 
type error bounds derivation were introduced by S. Repin in the 90's, see, e.g., 
\cite{Repin1997, RepinPowerFunc1997, Repin1999, Repin2000, NeittaanmakiRepin2004}.
Naturally, there are alternative approaches in a posteriori error control such as the residual-based or 
the goal-oriented estimation techniques, see, e.g., the books \cite{Verfurth1996, BangerthRannacher2003}.
Moreover, we refer to the works \cite{LazarovTomov2002, Verfurth2005} dedicated to convection-diffusion 
problems. Compared to other methods, the functional type a posteriori estimation techniques provide a 
general (only functional based) framework to derive 
independent and guaranteed 
upper bounds.
The initial work addressing elliptic convection-diffusion problems was done 
in \cite{RepinDeGruyter2008}, where the convection-dominated aspect was considered in detail
(see, e.g., \cite[Section 4.3.]{RepinDeGruyter2008}). The derivation of majorants for the time-dependent 
reaction-convection-diffusion problem has been discussed in \cite{RepinTomar2010}. However, 
the numerical aspects as well as properties of the majorant applied to this class of problems 
has not been studied so far. With the current work, we 
fill this gap and consider various 
aspects of two-sided error bounds in applications.
In particular, we present a set of numerical examples for the linear Fokker-Planck model problem
for different parameters and boundary conditions. A specially interesting case is the 
convection-dominated (CD) setting of the problem, which
may follow from various applications, such as the linearised Navier-Stokes equation with high 
Reynolds number or drift-diffusion equations of a semiconductor device modelling.
Due to the small diffusion coefficient, the solution of the discussed problem has singularities in
form of boundary or interior layers. The classical numerical approximations relying on equidistant 
meshes are not able to capture these layers unless the mesh contains an unacceptably large amount 
of nodes. If the mesh is not fine enough, the obtained approximations will include the non-physical 
oscillations polluting the data. 
The techniques developed to trigger the above described issue include the upwind scheme 
\cite{Heinrichetall1977},  streamline diffusion finite element method (SDFEM), also known
as streamline-upwind/Petrov-Galerkin formulation (SUPG) \cite{BrooksHughes1982}, 
residual free bubble (RFB) functions \cite{Brezzietall1997, BrezziFrancaRusso1998,Brezzietall1997, 
Brezzietall1999,BrezziMariniSuli2000}, 
Galerkin Least Squares (GLS) \cite{HughesFrancaHulbert1989, ArayaBehrensRodriguez2005, 
RoosStynesTobiska1996},
Continuous Interior Penalty (CIP) \cite{DouglasDupont1976}, 
Edge Stabilisation (ES) \cite{BurmanHansbo2004},
exponential fitting \cite{Banketall1990, Brezzietall1989, Nooyen1995, XuZikatanov1999}, 
discontinuous Galerkin methods \cite{BaumannOden1999, HoustonSchwabSuli2002}, 
and spurious oscillations at layers diminishing (SOLD) methods \cite{JohnKnobloch2007, 
JohnKnobloch2008}.
This list is by no means exhaustive, and there exist other techniques tackling this issue. 
In this work, we use the SUPG method introduced by Brooks and Hughes \cite{BrooksHughes1982}.

An alternative approach to handle the state Fokker-Planck boundary value problem,
see later 
\eqref{eq:fokker-planck-state-equation}--\eqref{eq:dual-state-equation},
relies on an adaptation (grading) of the underlying meshes in order to capture the boundary layers. 
The mesh adaptation for convection-dominated problems (CDPs) affects both the stability and 
the accuracy of the scheme. In practice, it was observed that once the accuracy is improved through 
the mesh adaption, the scheme gets stabilised as well (see, e.g., \cite{BanschMorinNochetto2002, 
RoosStynesTobiska1996,ZhangTang2002}). According to the numerical results presented in 
\cite{LiTangZhang2002,RoosStynesTobiska1996, ZhangTang2002} and the theoretical 
justification for a model problem discussed in \cite{RoosStynesTobiska1996, Zhang2003}, 
meshes can be adapted to boundary levels in such a way that even a standard finite element method 
can provide the optimal or quasi-optimal approximation property. 
However, the question of constructing such a grid, correctly adapted to the particular problem, still 
remains open. Moreover, in \cite{SunChenXu2010} the authors conclude that the stabilisation of the CDP scheme
is necessary.
Since a priori error analysis does not provide directions on mesh adaptation, 
an appropriate a posteriori error estimate, extracting corresponding information from the 
obtained numerical approximations, would be a crucial tool to construct desirable meshes. 
In this work's numerical examples, we perform an adaptive strategy 
based on a bulk marking criterion for selecting the elements which have to be refined, 
see \cite{Doerfler1996}. 

According to Sun et al. \cite{SunChenXu2010}, most a posteriori estimators constructed for 
elliptic type equations do not work well for the convection-dominated problems. John \cite{John2000} 
presents a thorough numerical comparison of some standard a posteriori error estimators 
(e.g., residual-based error estimator, Zienkiewicz-Zhu estimator and error estimator based on 
the solution of local Neumann problems) for CDPs. Moreover, it is shown that no error estimator 
performs equally good for all numerical examples due to the presence of the term 
$\tfrac{1}{\varepsilon}$ (see, e.g., Kang and Yu \cite{KangYu2001}). For $\varepsilon$ close to zero, 
the error indicator might indicate wrong regions for the refinement.

This work is structured as follows:  In Section \ref{Sec:FuncSpaceModelProb}, the relevant 
function spaces as well as the stationary and time-dependent linear Fokker-Planck model problems
are presented. The derivation of two-sided estimates for the stationary problems
is discussed in Section \ref{Sec:MajMin}, followed by the corresponding results for the 
time-dependent case in Section \ref{sec:time-dependent}. Finally, first numerical results for various 
given data sets are presented and discussed in Section \ref{Sec:NumRes}.

\section{Function spaces and Fokker-Planck model problems}
\label{Sec:FuncSpaceModelProb}

This section presents a short summary of the functional spaces' definitions needed throughout the 
paper. We also state the model problem for the static and time-dependent cases, which will be 
addressed from the a posteriori error control point of view in the upcoming sections.

\subsection{Function Spaces}

Let $\Omega \subset \Rd$ , $d = \{1, 2, 3\}$, be a bounded domain with Lipschitz 
boundary $\Gamma = \partial \Omega$, where $\overline{\Omega}$ is the closure of $\Omega$.
We denote by $\Gamma_N$ and $\Gamma_D$ with $\Gamma = \Gamma_N \cup \Gamma_D$, 
$\meas_{d - 1} \Gamma_N > 0$ and $\meas_{d - 1} \Gamma_D > 0$, the parts of the boundary 
$\Gamma$, where Robin or Dirichlet boundary conditions (BCs) are imposed, respectively.
In this work, we consider problems with mixed and Robin boundary conditions.
We denote by $L_2(\Omega)$ the Banach space of all Lebesgue-measurable functions $v$ defined 
on $\Omega$ for which
\begin{align*}
 \| v \|_{\Omega} := 
 \left( \int_\Omega |v(x)|^2 \dx \right)^{1/2} < \infty.
\end{align*}
In addition, we define the following weighted $L_2$ norm:
\begin{alignat}{2}
\| v \|_{\varepsilon,\Omega} := \varepsilon \left( \int_\Omega |v(x)|^2 \dx \right)^{1/2}.
\label{eq:weigthedL2nrom}
\end{alignat}
We define the $H^1$ spaces
\begin{alignat}{2}
V_0 := \{v \in  H^{1}(\Omega) \;\big|\; v = 0 \text{ on } \Gamma_D \},
\label{eq:V0}
\end{alignat}
and 
\begin{alignat}{2}
\tilde V := \{ v \in H^{1}(\Omega) \;\big|\; \{  v \}_{\Omega} = 0 \}, \, \text{ where } \,
\{  v \}_{\Omega} := \frac{1}{|\Omega|} \int_{\Omega} v \dx
\label{eq:definition:zeromean}
\end{alignat}
with underlying standard $L_2$ spaces. Here, the first space includes the homogeneous Dirichlet 
boundary conditions and the second one contains functions with zero mean in $\Omega$.
Moreover, we introduce the space 
\begin{alignat}{2}
	Y_{\dvrg}(\Omega, \Gamma_N) := 
	\Big \{  \flux \in L_{2}\big(\Omega, \Rd\big) \;\big|\; 
	& \dvrg \, \flux \in L_{2}(\Omega),
	\, \, \flux \cdot {\bf n} \in L_{2}(\Gamma_N) \Big \}
	\label{eq:y-set-div}
\end{alignat} 	
of vector-valued functions satisfying the identity
\begin{equation}
	\int_{\Omega} \dvrg \, (\flux \, w) \dx
	= \int_{\Omega} (\dvrg \, \flux \, w + \flux \cdot \nabla w) \dx
	= \int_{\Gamma_N} (\flux \cdot {\bf n}) \, w  \ds, \quad
	\forall w \in V_0.
\label{eq:mix:y-identity}
\end{equation}
Let us also define the space
\begin{alignat}{2}
L^{\infty}_{[0, 1]} (\Omega)
:= \big\{ \chi \in \L{\infty}(\Omega) \; | \;  0 \leq \chi \leq 1 \;\; \mbox{a.a.} \;\; \Omega \big \}.
\label{eq:L-inf-space}
\end{alignat} 	
Let $(0, T)$, $0 < T < +\infty$, be a given time interval. 
Then, we denote by $Q_T := \Omega \times (0, T)$ and $S_T = \partial \Omega \times (0, T)$ 
the space-time cylinder and its lateral surface, respectively. We define the space 
$$H^{1, 1}(Q_T) = \big\{u \in L^2(Q_T) \, | \, 
\nabla_x u \in [L^2(Q_T)]^d, \partial_t u \in L^2(Q_T) \big\}$$
in space and time with zero mean in $\Omega$ for a.a $t \in \, (0, T)$ as follows
\begin{alignat}{2}
\tilde V_0 =  
\big\{ v \in H^{1, 1}(Q_T) \, | \, \{  v \}_{\Omega} = 0 \, \text{ for a.a. } t \in \, (0, T) \big\}.
\label{eq:H11tilde}
\end{alignat}
Analogously to \eqref{eq:y-set-div}, we define the space
\begin{alignat*}{2}
	\flux \in Y_{\dvrg}(Q_T) := 
	\Big \{  \flux \in L_{2} \big(0, T; L_{2}\big(\Omega, \Rd\big)\big)\;\big|\; 
	& \dvrg \, \flux \in L_{2} \big(0, T; L_{2}(\Omega)\big), \,\, 
	\flux \cdot {\bf n} \in L_{2} \big(0, T; L_{2}(\Gamma)\big) \Big \}
	\label{eq:y-set-div-time}
\end{alignat*} 	
by means of equation
\begin{equation}
	\int_{Q_T} \dvrg \, (\flux \, w) \dxt
	= \int_{Q_T} (\dvrg \, \flux \, w + \flux \cdot \nabla w) \dxt 
	= \int_{S_T} \flux \cdot {\bf n}  \dst, \quad
	\forall w \in \tilde V_0.
	\label{eq:y-identity-time}
\end{equation}
Finally, we define several inequalities, and corresponding to them constants,
that are used throughout the work. By $\CFriedrichs$  \cite{Friedrichs1937} 
we denote the constant in Friedrichs' inequality given as follows
\begin{equation}
\| w \|_{Q_T} \leq \CFriedrichs \, \| \nabla w \|_{Q_T}, \quad 
\forall w \in \HD{1, 0}{0}(Q_T) := \Big\{ \, w \in L_{2}(Q_T) \; : \nabla_x w \in [L_{2}(Q_T)]^d, 
\; w = 0 \; {\rm on} \; \Sigma \,\Big\}.
\label{eq:friedrichs-inequality}
\end{equation}
The Poincar\'{e} inequality \cite{Poincare1890,Poincare1894} reads as 
\begin{equation}
\| w \|_{Q_T} \leq C_{\mathrm{P}} \, \| \nabla w \|_{Q_T}, \quad
\forall w \in \tilde{H}^{1, 0}_{0}(Q_T) := \Big\{ \, w \in \HD{1, 0}{0}(Q_T) \; : \; \{ w \}_{\Omega} = 0 \Big\}.
\label{eq:classical-poincare-constant}
\end{equation}
Finally, the classic trace inequality is given by
\begin{equation}
	\| u \|_{L_{2}(S_T)} \leq C_{{\rm tr}} \| \, u \,\|_{\H{1}(Q_T)} \,,
	\qquad \forall u \in C^{1}(\overline{Q_T}).
	\label{eq:trace-inequality}
\end{equation}

\subsection{Fokker-Planck model problems}

The state Fokker-Planck boundary value problem (BVP) reads as follows
\begin{alignat}{2}
	- \, \dvrg \, {\boldsymbol p} + \dvrg ({\bf F} u) & =\, f,	      
	& \quad x \in \Omega,\label{eq:fokker-planck-state-equation}\\
	{\boldsymbol p} & = \varepsilon \, \nabla u, & \quad x \in \Omega, \label{eq:dual-state-equation}
\end{alignat}
where $f \in L_{2}(\Omega)$.
In the following, we consider the problem 
\eqref{eq:fokker-planck-state-equation}--\eqref{eq:dual-state-equation}
for two cases, which differ with respect to (w.r.t.) the boundary conditions imposed on it.
In the first problem, we impose mixed BCs of the form
\begin{alignat}{2}
	({\bf F} u - \, {\boldsymbol p}) \cdot {\bf n} & =\, 0,			
	&& \qquad x \in \Gamma_N, \label{eq:mix:neumann-bc}\\	
							u & = \,0,
	&& \qquad x \in \Gamma_D, \label{eq:mix:dirichlet-bc}
\end{alignat}
on \eqref{eq:fokker-planck-state-equation}--\eqref{eq:dual-state-equation}, and in the second 
one, Robin BCs as follows
\begin{alignat}{3}
	({\bf F} u - \, {\boldsymbol p}
	) \cdot {\bf n} & =\, 0,					
	& \quad x \in \Gamma. \label{eq:robin-bc}
\end{alignat}
Here, ${\bf n}$ denotes the vector of unit outward normal to $\partial\Omega$,
and the function ${\bf F}$ satisfies the following regularity requirements
\begin{alignat}{1}
	{\bf F} \in L^{\infty} (\Omega, \Rd), \quad \dvrg \,{\bf F} \in L^{\infty} (\Omega), 
         \label{eq:static-coefficients-condition}
\end{alignat}
and the incoming flux conditions 
\begin{alignat}{2}
	\quad {\bf F} \cdot {\bf n} < 0  \quad \mbox{on} \quad  \Gamma_N, \quad \text{and} \quad
	\quad {\bf F} \cdot {\bf n} < 0  \quad \mbox{on} \quad  \Gamma,
	\label{eq:static-coefficients-condition:incoming-flux-condition}
\end{alignat}
in the case of mixed and Robin BCs, 
respectively.

The time-dependent version of the problem 
\eqref{eq:fokker-planck-state-equation}--\eqref{eq:dual-state-equation}
with Robin boundary conditions is discussed in section \ref{sec:time-dependent} as well.
It reads as follows: find $u \in H^{1, 1}(Q_T)$ satisfying
\begin{alignat}{3}
	u_t  -  \dvrg \, {\boldsymbol p}  + 
	{\dvrg \, }({\bf F} u) & =\, f,	      
	& \quad (x, t) \in Q_T,\label{eq:fokker-planck-equation}\\
	{\boldsymbol p}  & = \varepsilon \,  \nabla u, 
	& \quad (x, t) \in Q_T, \label{eq:dual-equation}\\	
	u(x, 0) & =\, u_0,		
	& \quad x \in \Omega,
	\label{eq:parabolic-initial-condition}
	\\
  ({\bf F} u - {\boldsymbol p}
  ) \cdot {\bf n} & =\, 0,					
	& \quad (x, t) \in S_T,\label{eq:parabolic-robin-bc}
\end{alignat}
where
\begin{equation}
f \in L_{2}(Q_T), \quad \mbox{and} \quad
u_0 \in \tilde V. 
\label{eq:problem-condition}
\end{equation}
Function ${\bf F}$ satisfies \eqref{eq:static-coefficients-condition} for a.a 
$t \in \, (0, T)$ and the incoming flux condition ${\bf F} \cdot {\bf n} < 0$ on $S_T$.

\section{Two-sided estimates for the static case}
\label{Sec:MajMin}

The current section is dedicated to the derivation of the majorant and minorant of the distance between 
any approximation function $v$ and the exact solution $u$ of static or time-dependent Fokker-Planck problem.
Presented upper bounds can be viewed as a generalisation of the majorant studied in 
\cite[Section 4.3.]{RepinDeGruyter2008} and \cite{RepinTomar2010}. The derived lower bound of the error 
is a completely original result. 

\subsection{Error majorant in case of mixed boundary conditions}

By testing (\ref{eq:fokker-planck-state-equation}) with a function $\eta \in V_0$, we arrive at the 
generalised formulation of \eqref{eq:fokker-planck-state-equation}--\eqref{eq:dual-state-equation}
with boundary conditions \eqref{eq:mix:neumann-bc}--(\ref{eq:mix:dirichlet-bc}): 
find $u \in V_0$ satisfying the integral identity
\begin{equation}
	a(u,\eta) = \, \langle f, \eta \rangle_{\Omega}, 
	\quad \forall \eta \in V_0,
	\label{eq:mix:state-generalized-statement}
\end{equation}
where the non-symmetric bilinear form $a(\cdot,\cdot)$ and the right-hand side (RHS) are given by
\begin{equation}
	a(u,\eta) := \varepsilon \int_{\Omega} \nabla u 
	\cdot \nabla \eta \dx \, 
	- \, \int_{\Omega} ({\bf F} u) \cdot \nabla \eta \dx \qquad \text{and} \qquad
	\langle f, \eta \rangle_{\Omega} := \, \int_{\Omega} f \eta \dx,
	\label{eq:mix:nonsymBFAndRHS}
\end{equation}
respectively.
Assume that function $v \in V_0$ is the approximation of solution $u$ and the distance 
$e = u - v \in V_0$ between them is measured in terms of the norm
\begin{equation}
|\!|\!| \, e \, |\!|\!|^2_{\overline{\rm M}} 
\equiv
	|\!|\!| \, e \, |\!|\!|^2_{(\nu, \delta, \chi)} 
	:= \!\nu \left \| \, \nabla e \, \right\|^2_{\Omega} \,
	+ \,  \left\| \, \delta \, e \right\|^2_\Omega
	+ \,  \left\| \, \chi \, e \right\|^2_{\Gamma},	
	\label{eq:static-energy-norm}
\end{equation}
where $\nu := \eps$ is a 
positive weight parameter, and 
$\delta^2 := \tfrac{1}{2} \, \dvrg {\bf F}$ and 
$\chi^2 := \big(-\tfrac{1}{2} \, {\bf F} \cdot {\bf n}\big)^{\rfrac{1}{2}}$ are positive 
weight functions.
By transforming \eqref{eq:mix:state-generalized-statement} and assuming $\eta = e$,
we arrive at the following equality
\begin{equation}
	\varepsilon \int_{\Omega} \nabla e \cdot \nabla e \dx \, 
	- \, \int_{\Omega} ({\bf F} e) \cdot \nabla e \dx \, \\
	= \, \int_{\Omega} (f \, e - \varepsilon \, \nabla v \cdot \nabla e + 
	({\bf F} v) \cdot \nabla e)\dx.
	\label{eq:mix:e-equality}
\end{equation}
Substituting the identity
\begin{equation}
- \int_{\Omega} ({\bf F} e) \cdot \nabla e \dx \, 
= \tfrac{1}{2} 
  \Big( \int_{\Omega} \dvrg {\bf F} \, e^2 \dx
         - \int_{\Gamma_N} ({\bf F} \cdot {\bf n}) \, e^2 \ds \Big)
\label{eq:mix:F-e-nabla-e}
\end{equation}
into \eqref{eq:mix:e-equality}, we obtain the relation known as `error-identity' for 
$e \in V_0$
\begin{equation}
	\varepsilon \int_{\Omega} |\nabla e|^2 \dx \, + \, 
	\tfrac{1}{2} \Big( \int_{\Omega} \dvrg {\bf F} \, e^2 \dx
         - \int_{\Gamma_N} ({\bf F} \cdot {\bf n}) \, e^2 \ds \Big) \, \\
	= \, \int_{\Omega} (f \, e - \varepsilon \, \nabla v \cdot \nabla e + 
	({\bf F} v) \cdot \nabla e)\dx. 
	\label{eq:mix:e-equality-2}
\end{equation}
Following the ideas exposed in monographs 
\cite{NeittaanmakiRepin2004, RepinDeGruyter2008, Malietall2014}, we introduce an auxiliary 
vector-valued function \linebreak $\flux \in Y_{\dvrg}(\Omega, \Gamma_N)$ satisfying the identity 
\eqref{eq:mix:y-identity}. By means of substituting \eqref{eq:mix:y-identity} into 
\eqref{eq:mix:e-equality-2}, we arrive at
\begin{multline}
	|\!|\!| \, e \, |\!|\!|^2_{(\nu, \delta, \chi)} 
	:= \| \nabla \, e \|^2_{\eps, \Omega} \, 
	+ \, \| \delta \, e \|^2_{\Omega} 
	+ \, \| \chi \, e \|^2_{\Gamma_N}\, \\
	= \, \int_{\Omega} (f -	\dvrg ({\bf F} v)  
	+ \, \dvrg \flux) \, e \dx + 	
	\int_{\Omega} (\flux - \varepsilon \nabla v) \cdot \nabla e \dx +
	 \int_{\Gamma_N} (({\bf F} v -\, \flux) \cdot {\bf n}) \, e \ds.
	\label{eq:mix:e-equality-3}
\end{multline}
From now on,
let ${\bf r}_{\rm eq}$, ${\bf r}_{\rm d}$, and ${\bf r}_{\rm \Gamma_N}$ denote the residuals in 
the RHS of \eqref{eq:mix:e-equality-3}
\begin{equation}
{\bf r}_{\rm eq} :=  f - \dvrg ({\bf F} v) + \, \dvrg \flux, \quad 
{\bf r}_{\rm d}   :=  \flux - \varepsilon \nabla{v}, \quad \mbox{and} \quad
{\bf r}_{\rm \Gamma_N} := ({\bf F} v - \, \flux) \cdot {\bf n},
\end{equation}
which also
denote the residuals of the equations \eqref{eq:fokker-planck-state-equation}, 
\eqref{eq:dual-state-equation}, and \eqref{eq:mix:neumann-bc}, respectively.
In order to estimate the RHS of \eqref{eq:mix:e-equality-3}, we apply the H\"{o}lder inequality, i.e., 
\begin{equation}
\int_{\Omega} {\bf r}_{\rm d} \cdot \nabla e \dx 
\leq \| {\bf r}_{\rm d} \|_{\varepsilon^{-1}, \Omega}\, \| \nabla e\|_{\varepsilon, \Omega}.
\label{eq:holder-1}
\end{equation}
For treatment of the terms ${\bf r}_{\rm eq}$ and ${\bf r}_{\rm \Gamma_N}$, we introduce 
weight functions $\mu, \theta$ from the space $L^{\infty}_{[0, 1]} (\Omega)$
in order to obtain robust majorants w.r.t. drastic changes in terms of
$\dvrg {\bf F}$ and ${\bf F} \cdot {\bf n}$.
Then, the terms with residuals ${\bf r}_{\rm eq}$ and ${\bf r}_{\rm \Gamma_N}$ can be 
estimated as follows:
\begin{alignat}{2}
\int_{\Omega} {\bf r}_{\rm eq} \, e \dx
= \int_{\Omega} \mu \, {\bf r}_{\rm eq} \, e \dx 
  + \int_{\Omega} (1 - \mu) \,{\bf r}_{\rm eq}\, e  \dx
\leq \| \tfrac{\mu}{\delta} \, {\bf r}_{\rm eq} \|_{\Omega} \| \delta e\|_{\Omega} 
+ \tfrac{\CFriedrichs}{\sqrt{\varepsilon}} \| (1 - \mu) \,{\bf r}_{\rm eq} \|_{\Omega} \,  \| \nabla e \|_{\varepsilon, \Omega}
\label{eq:holder-2}
\end{alignat}
and 
\begin{alignat}{2}
\int_{\Omega} {\bf r}_{\rm \Gamma_N} \, e \dx
= \int_{\Omega} \theta \, {\bf r}_{\rm \Gamma_N} \, e \dx + \int_{\Omega} (1 - \theta) \,{\bf r}_{\rm \Gamma_N} \, e \dx
\leq \| \tfrac{\theta}{\chi} \, {\bf r}_{\rm \Gamma_N} \|_{\Omega} \| \chi e\|_{\rm \Gamma_N}
+ \tfrac{\Ctr }{\sqrt{\varepsilon}}\, \| (1 - \theta) \,{\bf r}_{\rm \Gamma_N} \|_{\rm \Gamma_N} \, \| \nabla e \|_{\varepsilon, \Omega}.
\label{eq:holder-3}
\end{alignat}
The resulting majorant follows {from the collection of} 
\eqref{eq:holder-1}--\eqref{eq:holder-3} 
obtained by means of Young's inequality:
\begin{alignat*}{2}
|\!|\!| \, e \, |\!|\!|^2_{(\nu, \delta, \chi)} 
& \leq \Big( \| {\bf r}_{\rm d}\|_{\varepsilon^{-1}, \Omega}
		 + \tfrac{\CFriedrichs}{\sqrt{\varepsilon}} \, \| (1 - \mu) \,{\bf r}_{\rm eq} \|^2_{\Omega}
		 + \tfrac{\Ctr }{\sqrt{\varepsilon}}\, \| (1 - \theta) \,{\bf r}_{\rm \Gamma_N} \|_{\Gamma_N} \Big)^2
		+ \| \tfrac{\mu}{\delta} \, {\bf r}_{\rm eq} \|^2_{\Omega} 
		+ \| \tfrac{\theta}{\chi} \, {\bf r}_{\rm \Gamma_N} \|^2_{\Omega} \\
 & 
 \leq \Big(
 \tfrac{\beta}{2} \,\| {\bf r}_{\rm d} \|^2_{\varepsilon^{-1}, \Omega} 
 + \tfrac{1}{2\, \beta} \, \big( \tfrac{\zeta}{2} \, \tfrac{\CFriedrichs^2}{\varepsilon} \, \| (1 - \mu) \,{\bf r}_{\rm eq} \|^2_{\Omega} 
 + \tfrac{1}{2\, \zeta}\tfrac{\Ctr^2}{{\varepsilon}}\,\| (1 - \theta) \,{\bf r}_{\rm \Gamma_N} \|^2_{\Gamma_N} \big)\Big)
 + \| \tfrac{\mu}{\delta} \, {\bf r}_{\rm eq} \|^2_{\Omega} 
 + \| \tfrac{\theta}{\chi} \, {\bf r}_{\rm \Gamma_N} \|^2_{\Gamma_N},
\end{alignat*}
where $\beta$ and $\zeta$ are positive parameters. Therefore, 
\begin{equation}
|\!|\!| \, e \, |\!|\!|^2_{\overline{\rm M}} 
\equiv 
|\!|\!| \, e \, |\!|\!|^2_{(\nu, \delta, \chi)} 
\leq \inf_{\flux \in Y_{\dvrg}(\Omega)} 
       \overline{\mathrm M}^2_{\Gamma_N, \Omega}(v, \flux; \beta, \zeta),
       \label{eq:majorant-static}
\end{equation}
where
\begin{equation}
\overline{\mathrm M}^2_{\Gamma_N, \Omega} (v, \flux; \beta, \zeta) 
:= \tfrac{1}{2} \bigg( \beta \,\| {\bf r}_{\rm d}\|^2_{\varepsilon^{-1}, \Omega} + 
    \tfrac{1}{2\, \beta}
    \Big( \zeta \, \tfrac{\CFriedrichs^2}{\varepsilon} \, \| {\bf r}_{\rm eq} \|^2_{\Omega} + 
	   \tfrac{1}{\zeta} \, \tfrac{\Ctr^2}{{\varepsilon}} \, \| {\bf r}_{\rm \Gamma_N} \|^2_{\Gamma_N}\Big)
    \bigg) 	
    + \| \tfrac{\mu}{\delta} \, {\bf r}_{\rm eq} \|^2_{\Omega} 
    + \| \tfrac{\theta}{\chi} \, {\bf r}_{\rm \Gamma_N} \|^2_{\Gamma_N}.
\end{equation}

\subsection{Error majorant in case of Robin boundary conditions}

Next, we consider the state Fokker-Planck BVP
\eqref{eq:fokker-planck-state-equation}--\eqref{eq:dual-state-equation}
with Robin BCs \eqref{eq:robin-bc}. 
By testing (\ref{eq:fokker-planck-state-equation}) with a function $\eta \in \tilde V$, 
we arrive at the generalized formulation of 
\eqref{eq:fokker-planck-state-equation}-\eqref{eq:dual-state-equation} with (\ref{eq:robin-bc}): 
find $u \in \tilde V$ satisfying the integral identity 
$a(u,\eta) = \, \langle f, \eta \rangle_{\Omega}$ for all $\eta \in \tilde V$, where
$a(u,\eta)$ and $\langle f, \eta \rangle_{\Omega}$ are defined in 
\eqref{eq:mix:nonsymBFAndRHS}.
Assume that function $v \in \tilde V$ is the approximation of solution $u$ and the error 
$e = u - v \in \tilde V$ is measured in terms of the norm \eqref{eq:static-energy-norm}.
By transforming $a(u,\eta) = \, \langle f, \eta \rangle_{\Omega}$ and assuming $\eta = e$, 
we arrive at the error-identity
\begin{equation}
	\varepsilon \int_{\Omega} \nabla e \cdot \nabla e \dx \, + \, 
	\tfrac{1}{2} \Big( \int_{\Omega} \dvrg {\bf F} \, e^2 \dx
         - \int_{\Gamma} ({\bf F} \cdot {\bf n}) \, e^2 \ds \Big) \, \\
	= \, \int_{\Omega} (f \, e - \varepsilon \, \nabla v \cdot \nabla e + 
	({\bf F} v) \cdot \nabla e)\dx,
	\label{eq:e-equality-2}
\end{equation}
{which is similar to \eqref{eq:mix:e-equality-2}.}
Analogously to the previous subsection, we introduce an auxiliary vector-valued function 
$\flux \in Y_{\dvrg}(\Omega, \Gamma)$, which yields the relation
\begin{equation}
	\| \nabla e \|^2_{\varepsilon, \Omega} \, + \, 
	\tfrac{1}{2} \| (\dvrg {\bf F})^{\rfrac{1}{2}} e \|^2_{\Omega} + 
  \tfrac{1}{2} \| (- {\bf F} \cdot {\bf n})^{\rfrac{1}{2}} \, e \|^2_{\Gamma}\,
	= \, \int_{\Omega} {\bf r}_{\rm eq} \, e \dx + 	
	 \int_{\Omega} {\bf r}_{\rm d} \cdot \nabla e \dx +
	 \int_{\Gamma} ({\bf r}_{\rm \Gamma} \cdot {\bf n}) \, e \ds.
	\label{eq:e-equality-3}
\end{equation}
In order to estimate the RHS of \eqref{eq:e-equality-3}, we apply the H\"{o}lder, Young, 
Poincar\'e, and trace inequalities leading to
\begin{equation}
{|\!|\!|e|\!|\!|^2_{\overline{\rm M}} \equiv }
|\!|\!| \, e \, |\!|\!|^2_{(\nu, \delta, \chi)} 
\leq \inf_{\flux \in Y_{\dvrg}(\Omega)} 
\overline{\mathrm M}^2_{\Gamma, \Omega} (v, \flux; \beta, \zeta) ,
\end{equation}
where
\begin{multline}
\overline{\mathrm M}^2_{\Gamma, \Omega} (v, \flux) 
:= \tfrac{1}{2}\, \Big( \beta \, \| {\bf r}_{\rm d} \|^2_{\varepsilon, \Omega} + 
    \tfrac{1}{2\, \beta} \, 
    \Big( \zeta\, \tfrac{\CPoincare^2}{\varepsilon} \, \| (1 - \mu) {\bf r}_{\rm eq} \|^2_{\Omega} + 
	       \tfrac{1}{\zeta} \, \tfrac{\Ctr^2}{{\varepsilon}} \, \| (1 - \chi){\bf r}_{\rm \Gamma_N} \|^2_{\Gamma_N} \big)\Big)
	+ \| \tfrac{\mu}{\delta} \, {\bf r}_{\rm eq} \|^2_{\Omega} 
		+ \| \tfrac{\theta}{\chi} \, {\bf r}_{\rm \Gamma_N} \|^2_{\Gamma_N}.
\end{multline}

\subsection{Error minorant in case of mixed boundary conditions}

In the following, we derive a lower bound of the approximation error in terms of a different quantity.
\begin{theorem}
For $v \in V_0$, the following  inequality holds:
\begin{equation}
\sup\limits_{w \in V_0} \underline{\rm M}^2_{\Omega}(v, w) 
{:= |\!|\!|e|\!|\!|^2_{\underline{\rm M}} }
\equiv 
|\!|\!| \, e \, |\!|\!|^2_{(\tilde \nu, \tilde \delta, \chi)},
\label{estimate:minorant}
\end{equation}
where 
$\tilde \nu = \tfrac{\varepsilon}{2}$, 
$\tilde \delta^2 = \tfrac{1}{2} \Big(\tfrac{|{\bf F}|^2}{\varepsilon} + \dvrg {\bf F} \Big)$, 
$\chi^2 = - \tfrac{1}{2} {\bf F} \cdot {\bf n}$,
and 
\begin{equation}
\underline{\rm M}^2_{\Omega}(v, w) := 
\int_{\Omega} (f w - \varepsilon \nabla v \cdot \nabla w + {\bf F} v \cdot \nabla w) \dx 
	- \frac{1
	}{2} \|\nabla w\|^2_{\varepsilon, \Omega}.
\label{eq:minorant}
\end{equation}
\end{theorem}
\ProofBegin
The supremum can be estimated as follows
\begin{align}
\begin{aligned}
       \sup_{w \in V_0} & \left\{ \int_{\Omega} (f w - \varepsilon \nabla v \cdot \nabla w + {\bf F} v \cdot \nabla w) \dx 
	- \frac{\varepsilon}{2} \|\nabla w\|^2_{\Omega} \right\} \\
	& \qquad \qquad
	= \sup_{w \in V_0} \left\{ \int_{\Omega} \Big(\varepsilon \nabla (u-v) - {\bf F} (u-v)\Big) \cdot \nabla w \dx 
	- \frac{\varepsilon}{2} \|\nabla w\|^2_{\Omega} \right\} \\
	& \qquad \qquad 
	\leq \sup_{\tau \in L^2(\Omega)} \left\{ \int_{\Omega} (\varepsilon \nabla e - {\bf F} e) \cdot \tau \dx 
	- \frac{\varepsilon}{2} \|\tau\|^2_{\Omega} \right\} \\
	& \qquad \qquad 
	\leq \frac{1}{2 \varepsilon} \| \varepsilon \nabla e - {\bf F} e \|^2_{\Omega}
	=  \frac{\varepsilon}{2} \|\nabla e\|^2_{\Omega} + \frac{1}{2 \varepsilon} \| {\bf F} e \|^2_{\Omega}
	- \int_{\Omega} \nabla e \cdot {\bf F} e \dx.
\end{aligned}
\end{align}
Substituting identity \eqref{eq:mix:F-e-nabla-e} then yields the following estimate:
\begin{align}
\begin{aligned}
       \sup_{w \in V_0} & \left\{ \int_{\Omega} (f w - \varepsilon \nabla v \cdot \nabla w + {\bf F} v \cdot \nabla w) \dx 
	- \frac{1}{2} \|\nabla w\|^2_{\eps,\Omega} \right\} \\
	& \qquad \qquad 
	\leq 
         \frac{1}{2} \|\nabla e\|^2_{\eps,\Omega} + \frac{1}{2 \varepsilon} \| {\bf F} e \|^2_{\Omega}
	+ \frac{1}{2} \| (\dvrg {\bf F})^{\rfrac{1}{2}} e \|^2_{\Omega}
         + \frac{1}{2} \| (- {\bf F} \cdot {\bf n})^{\rfrac{1}{2}} \, e \|^2_{\Gamma} \\
	& \qquad \qquad
	=  \frac{1}{2} \|\nabla e\|^2_{\eps,\Omega} + \frac{1}{2} \Big\| \Big(\tfrac{|{\bf F}|^2}{\varepsilon} + \dvrg {\bf F}\Big)^{\rfrac{1}{2}} e \Big\|^2_{\Omega}
         + \frac{1}{2} \| (- {\bf F} \cdot {\bf n})^{\rfrac{1}{2}} \, e \|^2_{\Gamma},
\end{aligned}
\end{align}
which finally leads to equation \eqref{estimate:minorant}.

\ProofEnd

\begin{remark}
Note that the new norm differs from the norm $|\!|\!| \, e \, |\!|\!|^2_{(\nu, \delta, \chi)}$ only by the 
weighting parameters that are given by
$\tilde \nu = \frac{\nu}{2}$ and 
$\tilde \delta^2 = \delta^2 + \tfrac{|{\bf F}|^2}{2 \varepsilon}$.
The relation can also be written as follows
\begin{align}
\begin{aligned}
         |\!|\!| \, e \, |\!|\!|^2_{(\tilde \nu, \tilde \delta, \chi)} =
	|\!|\!| \, e \, |\!|\!|^2_{(\nu, \delta, \chi)} + \frac{1}{2 \, \varepsilon} \| {\bf F} e \|^2_{\Omega}
	-  \frac{\varepsilon}{2} \|\nabla e\|^2_{\Omega}.
\end{aligned}
\end{align}
\end{remark}

When maximising the minorant functional, the auxiliary function $w \in V_0$ must be chosen from the 
richer approximation space in comparison to the approximation $v \in V_0$. Otherwise, the 
minorant vanishes. Analogously to the majorant functional, the minorant includes only known data and, 
therefore, is fully computable. By an appropriate selection of $w$, we can find a lower bound arbitrarily
close to the exact error. In particular, if $w = u - v$, the minorant coincides with the error norm.

\section{
Two-sided estimates for the time-dependent case}
\label{sec:time-dependent}

\subsection{Upper estimates for the time-dependent case}

We consider now the parabolic initial-boundary value problem.
After multiplying (\ref{eq:fokker-planck-equation}) by a test function 
$\eta \in \tilde V_0$, we arrive at the generalized formulation of 
(\ref{eq:fokker-planck-equation})--(\ref{eq:parabolic-robin-bc}): 
find $u \in \tilde V_0$ satisfying the integral identity
\begin{multline}
	\int_{\Omega} \big( (u\eta)(x, T) - (u\eta)(x, 0) \big) \dx \, 
	+ \, \int_{Q_T} u \eta_t \dxt \,
	+ \, \int_{Q_T} \varepsilon \, \nabla u \cdot \nabla \eta \dxt \, 
	- \, \int_{Q_T} ({\bf F} u) \cdot \nabla \eta \dxt \, \\
	= \, \int_{Q_T} f \eta \dxt, 
	\quad \forall \eta \in \tilde V_0.
	\label{eq:generalized-statement-ibvp}
\end{multline}
We present the functional error estimate, which provides a guaranteed upper bound of the 
deviation $e = u - v$ between the generalised solution of (\ref{eq:generalized-statement-ibvp}) and 
function $v \in \tilde V_0$ (generated by any numerical method) measured in terms of the norm
\begin{equation}
	\error_{(\nu, \theta, \chi, \zeta)}
	{\equiv \error_{\overline{\rm M}, Q_T}}
	:= \nu \left \| \, \nabla e \, \right\|^2_{Q_T} \,
	+ \, \left\| \, \theta \, e \right\|^2_{Q_T}  
	+ \, \left\| \, \chi \, e \right\|^2_{S_T}
	+ \, \zeta \! \, \left \| \, e (\cdot, T) \, \right \|^2_{\Omega},	
	\label{eq:energy-norm}
\end{equation}
where $\nu$ and $\zeta$ are positive weight parameters, and $\theta, \chi$ are weight functions.

The initial step in the derivation of both upper estimates is the transformation of 
\eqref{eq:generalized-statement-ibvp} into the integral identity
\begin{multline}
	\int_{\Omega} \big( (e\eta)(x, T) - (e\eta)(x, 0) \big) \dx \, 
	+ \, \int_{Q_T} e \eta_t \dxt \,
	+ \, \int_{Q_T} \varepsilon \, \nabla e \cdot \nabla \eta \dxt \, 
	- \, \int_{Q_T} ({\bf F} e) \cdot \nabla \eta \dxt \, \\	
	= \int_{Q_T} \!\! \left(f  \eta + v_t  \eta + ({\bf F} v) \cdot \nabla \eta
	                         - \varepsilon \, \nabla{v} \cdot \nabla \eta \right) \dxt.
	\label{eq:energy-balance-equation}
\end{multline}
By substituting $\eta = e$ and using the relation
\begin{equation}
\int_{\Omega} \big( (e\eta)(x, T) - (e\eta)(x, 0) \big) \dx \, + \, \int_{Q_T} e \eta_t \dxt \,
= \, \tfrac{1}{2} \! \, \Big(\left \| \, e (\cdot, T) \, \right \|^2_{\Omega} 
- \left \| \, e (\cdot, 0) \, \right \|^2_{\Omega}\Big), 
\end{equation}
we obtain 
\begin{multline}
	\tfrac{1}{2} \, \big(\left \| \, e (\cdot, T) \, \right \|^2_{\Omega} 
	- \left \| \, e (\cdot, 0) \, \right \|^2_{\Omega}\big)
	+ \, \varepsilon \int_{Q_T} |\nabla e |^2 \dxt \, 
	- \, \int_{Q_T} ({\bf F} e) \cdot \nabla e \dxt \, \\	
	= \int_{Q_T} \!\! \left( f e + v_t e - \dvrg ({\bf F} v) e
	                         - \varepsilon \, \nabla{v} \cdot \nabla e \right) \dxt +
		\int_{S_T} ({\bf F} v) \cdot {\bf n} \, e \dst.											
	\label{eq:energy-balance-equation-1}
\end{multline}
Taking into account the equality
\begin{equation}
- \int_{Q_T} ({\bf F} e) \cdot \nabla e \dxt \, 
= \tfrac{1}{2} 
  \Big( \int_{Q_T} \dvrg {\bf F} \, e^2 \dxt 
         - \int_{S_T} ({\bf F} \cdot {\bf n}) \, e^2 \dst \Big), 
\end{equation}
we rewrite the LHS of \eqref{eq:energy-balance-equation-1} as follows
\begin{multline}
	\tfrac{1}{2} \, \big(\left \| \, e (\cdot, T) \, \right \|^2_{\Omega} 
	- \left \| \, e (\cdot, 0) \, \right \|^2_{\Omega}\big)
	+ \, \varepsilon \int_{Q_T} |\nabla e |^2 \dxt \, + 
	\tfrac{1}{2} \Big( \int_{Q_T} \dvrg {\bf F} \, e^2 \dxt 
         + \int_{S_T} ( - {\bf F} \cdot {\bf n}) \, e^2 \dst \Big) \, \\	
	= \int_{Q_T} \!\! \Big( \big(f \, e - v_t \, e - \nabla ({\bf F} v) \, e \big) \dxt
	                         - \varepsilon \, \nabla{v} \cdot \nabla e \Big) \dxt +  
		\int_{S_T} ({\bf F} v) \cdot {\bf n} \, e \dst.											
	\label{eq:energy-balance-equation-1-2}
\end{multline}
We follow the ideas from \cite{Repin2002} and introduce {an arbitrary} vector-valued function 
$\flux \in Y_{\dvrg}(Q_T)$ leading to the equation
\begin{multline}
	\tfrac{1}{2} \! \, (\left \| \, e (\cdot, T) \, \right \|^2_{\Omega} - 
											\left \| \, e (\cdot, 0) \, \right \|^2_{\Omega})
	+ \, \varepsilon \int_{Q_T} |\nabla e |^2 \dxt \, 
	+ \tfrac{1}{2} \Big( \int_{Q_T} \dvrg {\bf F} \, e^2 \dxt 
         + \int_{S_T} ( - {\bf F} \cdot {\bf n}) \, e^2 \dst \Big) \, \\	
	= \int_{Q_T}
	  \big(f - v_t - \nabla ({\bf F} v)
		     + \, \dvrg \flux \big) \, e \dxt
	  + \int_{Q_T} 
		  \big( \flux - \varepsilon \nabla{v}\big) \cdot \nabla e \dxt +
		\int_{S_T} ({\bf F} v - \flux) \cdot {\bf n} \, e \dst.											
	\label{eq:energy-balance-equation-2}
\end{multline}
The residuals
\begin{equation}
{\bf r}_{\rm eq} =  f - v_t - \nabla ({\bf F} v) + \dvrg \flux, \qquad 
{\bf r}_{\rm d} = \flux - \varepsilon \, \nabla{v}, \qquad
{\bf r}_{\rm s_T} = ({\bf F} v - \flux) \cdot {\bf n}
\end{equation}
in the RHS of \eqref{eq:energy-balance-equation-2} {correspond in a natural manner}
to the equations 
\eqref{eq:fokker-planck-equation}, \eqref{eq:dual-equation}, and 
\eqref{eq:parabolic-robin-bc}, respectively.
We proceed by estimating the terms in \eqref{eq:energy-balance-equation-2} using the H\"{o}lder, 
Young, Friedrichs as well as trace inequality, i.e., 
\begin{alignat}{2}
\int_{Q_T} {\bf r}_{\rm eq} \, e \dxt & 
\leq \| {\bf r}_{\rm eq} \|_{Q_T} \, \| e \|_{\varepsilon, Q_T} 
\leq \tfrac{\alpha_1}{2} \, \tfrac{\CF^2}{\varepsilon} \, \| {\bf r}_{\rm eq} \|^2_{Q_T} +
 \tfrac{1}{2 \, \alpha_1}\|\nabla e \|^2_{\varepsilon, Q_T}, \\
\int_{Q_T} {\bf r}_{\rm d} \cdot \nabla e \dxt & 
\leq \| {\bf r}_{\rm d} \|_{Q_T} \, \| \nabla e \|_{\varepsilon, Q_T} 
\leq
\tfrac{\alpha_2}{2} \, \| {\bf r}_{\rm d} \|^2_{\varepsilon^{\minus1}, Q_T}  
+ \tfrac{1}{2\, \alpha_2}\|\nabla e \|^2_{\varepsilon, Q_T}, \\
\int_{Q_T} {\bf r}_{\rm s_T} \, e \dst & 
\leq \| {\bf r}_{\rm s_T} \|_{\rm S_T} \, \| e \|_{S_T} 
\leq 
\tfrac{\alpha_3}{2} \, \tfrac{\Ctr^2}{\varepsilon} \, \| {\bf r}_{\rm s_T} \|^2_{S_T} 
+ \tfrac{1}{2 \, \alpha_3}\|\nabla e \|^2_{\varepsilon, Q_T},
\end{alignat}
where functions $\alpha_i >0$, $i = 1, 2, 3$.
Finally, we arrive at the estimate
\begin{multline}
	\left \| \, e (\cdot, T) \, \right \|^2_{\Omega} + 
	(2 - \Sum_{i = 1}^{3}\, \tfrac{1}{\alpha_i}) \|\nabla e \|^2_{\varepsilon, Q_T} + 
	\| (\dvrg {\bf F})^{\rfrac{1}{2}} \, e \|^2_{Q_T} + 
  \| (- {\bf F} \cdot {\bf n})^{\rfrac{1}{2}} \, e \|^2_{S_T} 
          =: \error_{\overline{\rm M}, Q_T} \\	
	\leq \overline{\rm M}^2_{Q_T} 
	:= \left \| \, e (\cdot, 0) \, \right \|^2_{\Omega} + \Big(
	     \alpha_1 \, \tfrac{\CF^2}{\varepsilon} \, \| {\bf r}_{\rm eq} \|^2_{Q_T}
	     + \alpha_2 \| {\bf r}_{\rm d} \|^2_{\varepsilon^{\minus1}, Q_T}
	     + \alpha_3 \, \tfrac{\Ctr^2}{\varepsilon} \, \| {\bf r}_{\rm S_T} \|^2_{S_T}\Big),
	\label{eq:estimate-space-time}
\end{multline}
where parameter $\delta \in (0, 2]$.

\subsection{Lower estimates for the time-dependent case}
\label{sec:time-dependent-minorant}

Minorants of the error between the approximate and exact solutions for the evolutionary 
reaction-convection-diffusion problem provide useful information while testing upper 
error-bounds (when $u$ is not available). The quality of the majorant is evaluated by 
considering its ratio w.r.t. the minorant.
The first results related to lower error bounds for evolutionary problems were derived 
and numerically tested in \cite{MatculevichRepin2014}. The result presented 
in the current section can be viewed as {generalisation of 
Section} 3 in 
\cite{MatculevichRepin2014}.

\begin{theorem}
\label{th:theorem-minimum-of-minorant}
For any $v, \: \eta \in \tilde V_0$ the following estimate holds:
\begin{multline}
	\sup_{\eta \in \tilde V_0} \Min (v, \eta) := 
	\sup_{\eta \in \tilde V_0} \bigg \{  
	\sum_{i = 1}^{3} G_{i}(v, \eta) + F_{fu_0}(\eta) 
	\bigg \} \\
	\leq 
	\error_{\underline{\rm M}, Q_T} : = 
    \tfrac{1}{2} \| \nabla e \|^2_{\eps,Q_T} 
	+ \tfrac{1}{2} \, \Big\| \big(\tfrac{|{\bf F}|}{\varepsilon} + \dvrg {\bf F} + 1\big)^{\rfrac{1}{2}} \, e \Big\|^2_{Q_T} 
	+ \tfrac{1}{2} \, \Big\| ( - {\bf F} \cdot {\bf n})^{\rfrac{1}{2}} \, e \Big\|^2_{S_T}
	+ \tfrac{1}{2} \| e(x, T) \|^2_{\Omega}, 
	\label{eq:lower-estimate}
\end{multline}
where
\begin{alignat}{2}
G_{1}(v, \eta)
& := \int_{Q_T} \Big( - \nabla \eta \cdot (\varepsilon \nabla v + {\bf F} \, v) - \tfrac{\varepsilon}{2} | \nabla \eta |^2   \Big) \dxt, \nonumber \\
G_{2}(v, \eta)
& := \int_{Q_T} \Big( \eta_t v - \tfrac{1}{2} \, |\eta_t|^2 \Big) \dxt, \nonumber \\
G_{3}(v, \eta)
& := \int_{\Omega} \Big( - v(x, T) \eta(x, T) - \tfrac{1}{2} \, |\eta(x, T)|^2 \Big) \dx, \nonumber \\
F_{fu_0}(\eta)
& := \int_{Q_T} f \eta \dxt + \int_{\Omega} u_0\, \eta(x, 0) \dx.
\end{alignat}
\end{theorem}

\ProofBegin
It is not difficult to see that    
\begin{alignat}{2}
	\mathcal{M} (e) & := 
	\sup\limits_{\eta \in \tilde V_0}
	\bigg \{
	\int_{Q_T} \Big(
	\nabla \eta \cdot (\varepsilon \nabla e - {\bf F} \, e) - \tfrac{\varepsilon}{2} | \nabla \eta |^2 
     - \eta_t e - \tfrac{1}{2} |\eta_t|^2
	\Big) \dxt 
	+ \int_\Omega \Big( e(x, T) \eta (x, T) - \tfrac{1}{2} |\eta (x, T)|^2 \Big)\dx 
	\bigg \} 
    \nonumber \\ 	
	& \; \leq \sup\limits_{ {{\bftau} \in [L_{2}(Q_T)]^d}} \!
	\bigg \{ \int_{Q_T} \! \! 
	\Big( {\bf \tau} \cdot (\varepsilon \nabla e  - {\bf F} \, e ) 
	- \tfrac{\varepsilon}{2} | \tau |^2 \Big) \dxt \bigg \}
	+ \sup\limits_{\xi \in \tilde V_0} \!
	\bigg \{ \int_{Q_T} \! \! \Big( - \xi e - \tfrac{1}{2} |\xi|^2 \Big) \dxt \bigg \} 
	\nonumber\\
	& \qquad \qquad \qquad \qquad \qquad \qquad
	+ \sup\limits_{{\eta(\cdot,T) \in V_0}} \bigg \{ \int_\Omega \! \! 
	\Big( e(x, T) \eta (x, T) - \tfrac{1}{2} |\eta (x, T)|^2 \Big)\dx  \bigg \}.
 	\label{eq:quadratic-func-inequality}
\end{alignat}
Since
\begin{alignat*}{2}
\sup\limits_{{{\bftau} \in [L_{2}(Q_T)]^d}} 
\Bigg \{ \int\limits_{Q_T} \bigg( 
{\bf \tau} \cdot (\varepsilon \nabla e  - {\bf F} \, e ) - \tfrac{\varepsilon}{2} | \tau |^2 
\bigg) \dxt \Bigg \} & \leq \tfrac{1}{2\,\varepsilon} \| \varepsilon \nabla e  - {\bf F} \, e \|^2_{Q_T} \,, \nonumber \\
\sup\limits_{\xi \in \tilde V_0} \Bigg \{ \int\limits_{Q_T} \Big( - \xi e - \tfrac{1}{2} |\xi|^2 \Big) \dxt \Bigg \} & \leq \tfrac{1}{2} \| e \|^2_{Q_T}, \nonumber \\
\sup\limits_{{\eta(\cdot,T) \in V_0}} \Bigg \{ \int\limits_{\Omega} \Big( e(x, T) \eta (x, T) - \tfrac{1}{2} |\eta (x, T)|^2 \Big)\dx \Bigg \} & \leq \tfrac{1}{2} \| e(x, T) \|^2_{\Omega}, \nonumber\\
 \end{alignat*}
we find that 
\begin{alignat}{2}
	\mathcal{M} (e) 
	& \leq 
	\tfrac{1}{2 \, \varepsilon} \| \varepsilon \nabla e  - {\bf F} \, e \|^2_{Q_T}
	+ \tfrac{1}{2} \| e \|^2_{Q_T} 
	+ \tfrac{1}{2} \| e(x, T) \|^2_{\Omega} \nonumber\\
	& = \tfrac{1}{2} \| \nabla e \|^2_{\eps,Q_T} 
	+ \tfrac{1}{2 \, \varepsilon} \| {\bf F} \, e \|^2_{Q_T} 
	+ \tfrac{1}{2}\, \| (\dvrg {\bf F})^{\rfrac{1}{2}} \, e \|^2_{Q_T}
	+ \tfrac{1}{2} \, \| ( - {\bf F} \cdot {\bf n})^{\rfrac{1}{2}} \, e \|^2_{S_T}
	+ \tfrac{1}{2} \, \| e \|^2_{Q_T} 
	+ \tfrac{1}{2} \, \| e(x, T) \|^2_{\Omega} \nonumber\\
	& = \tfrac{1}{2} \| \nabla e \|^2_{\eps,Q_T} 
	+ \tfrac{1}{2} \, \Big\| \big(\tfrac{|{\bf F}|}{\varepsilon} + \dvrg {\bf F} + 1\big)^{\rfrac{1}{2}} \, e \Big\|^2_{Q_T} 
	+ \tfrac{1}{2} \, \| ( - {\bf F} \cdot {\bf n})^{\rfrac{1}{2}} \, e \|^2_{S_T}
	+ \tfrac{1}{2} \| e(x, T) \|^2_{\Omega}{,} 
	\label{eq:quadratic-func-inequality-1} 
\end{alignat}
{which is exactly defined as $\error_{\underline{\rm M}, Q_T}$ in (\ref{eq:lower-estimate}).}
On the other hand, by using {\eqref{eq:generalized-statement-ibvp}, we see that for any 
$\eta$ {the} functional 
\begin{equation}
	\mathcal{M} (e)
	= \\
	\sup\limits_{{\eta \in \tilde V_0}} \bigg \{  
	\sum_{i = 1}^{3} G_{i}(v, \eta) + F_{fu_0}(\eta) \bigg \}
\end{equation}
generates {a} lower bound for the norm $\error_{\underline{\rm M}, Q_T}$. 
Thus, we arrive at (\ref{eq:lower-estimate}).
\ProofEnd

\begin{remark}
In \eqref{eq:quadratic-func-inequality-1}, the term $\tfrac{1}{2} \, \| e \|^2_{Q_T}$ can 
be combined with $\tfrac{1}{2} \| \nabla e \|^2_{\eps,Q_T}$ by means of Friedrichs'
inequality, see \eqref{eq:friedrichs-inequality}.
\end{remark}

\begin{remark}
Since the majorants and minorants are equal to norms, see
\eqref{eq:majorant-static}, 
\eqref{estimate:minorant}, 
\eqref{eq:energy-norm} together with 
\eqref{eq:estimate-space-time}, and 
\eqref{eq:lower-estimate},
they are nonnegative.
In \cite{CarrilloCordierMancini2011}, positivity for the solution in the static as well as time-dependent
case is proved under the condition that 
$\boldsymbol{F} \in C^1 (\Omega, \mathbb{R}^2)$
and
$\boldsymbol{F} \cdot \boldsymbol{n} < 0  \mbox{ on }  \Gamma_N \text{ or } 
\boldsymbol{F} \cdot \boldsymbol{n} < 0  \mbox{ on }  \Gamma$
in the case of mixed or Robin BCs, respectively.
For the time-dependent case, the nonnegativity is provided for given nonnegative initial values $u_0$. 
Existence of a minimizing sequence for the majorant as
well as a maximizing sequence for the minorant, see for instance \cite{Malietall2014}, and corresponding 
infimum of the majorant and the supremum of the minorant, leads to the estimates of the accuracy 
of the approximation w.r.t. exact solution. However, it does not automatically imply positivity (or non-negativity).
%
Together with the results on the existence of a minimizing sequence for the majorant as
well as a maximizing sequence for the minorant, see for instance \cite[Section 3.6]{RepinDeGruyter2008},
the infimum of the majorant and the supremum of the minorant both
lead to the exact (nonnegative) solution.
\end{remark}

\begin{remark}
Note that the model problems 
\eqref{eq:fokker-planck-state-equation}--\eqref{eq:dual-state-equation}
and \eqref{eq:fokker-planck-equation}--\eqref{eq:parabolic-robin-bc} 
are written in 
form of a first-order system 
of conservation laws supplied by a proper constitutive relation for the flux.
This could yield to the extension of this work's
functional-based error analysis to the use of a dual mixed weak formulation of the PDE system,
which would be of great interest in the implementation of computations for more 
realistic applications, where the approximation of the total flux $({\bf F} u - \, {\boldsymbol p})$ 
is often more important than the approximation of $u$ itself. Studies concerning the
functional-based error analysis devoted to a posteriori error estimates of mixed approximations
for the Poisson equation with Dirichlet and mixed Dirichlet/Neumann boundary conditions can be found in 
\cite{Repinetall2007} and \cite{RepinSmolianski2005}, respectively. Moreover, in the latter work it is 
shown that after an appropriate scaling of the coordinates and the equation, the ratio of the upper and
lower bounds for the error in the product norm never exceeds the constant $3$. Such an a posteriori error 
approach is computationally very cheap and can also be used for the indication of the local error 
distribution. As an application, 
\cite{RepinSmolianski2005} 
considers the linear elasticity problem, where the reconstruction of the
stresses (dual variable) is equally important to the reconstruction of the displacement (primal variable).
\end{remark}

\section{Numerical examples}
\label{Sec:NumRes}

We have presented a generalized error control method, 
namely,
the functional error estimates, 
for convection-dominated diffusion problems in
Sections \ref{Sec:MajMin} and \ref{sec:time-dependent}.
They are derived independently from the numerical method used 
to reconstruct the approximation. {As a result},
one does not need 
to adapt the obtained error estimates to the stabilization technique chosen 
in the solver but can just use the reconstructed approximation. 

The current section {presents numerical results} for the problems
\eqref{eq:fokker-planck-state-equation}--\eqref{eq:dual-state-equation}
and \eqref{eq:fokker-planck-equation}--\eqref{eq:parabolic-robin-bc}
in the static and time-dependent case, respectively.
The examples were computed for various parameters and boundary conditions. 
A particularly interesting setting (discussed throughout this section) is the convection-dominated 
problem, namely, when $\varepsilon \ll  \| {\bf F}\|_{L^{\infty}}$. In that case, the solution 
contains a boundary layer, which can be either of the width $O(\varepsilon)$ (the so-called 
regular boundary layer) or of the width $O(\sqrt{\varepsilon})$ (parabolic boundary layer). 
The detailed discussion addressing properties for both of these layers can be found in 
\cite{Stynes2005}.

The stabilised form of \eqref{eq:mix:state-generalized-statement} can be written as follows:
$$
a(u, \eta) + \sum_{K \in \mathcal{K}_h} \delta_K \, s_K(u,  \eta) =  \langle \, f, \eta \, \rangle_{\Omega},
$$
where $s_K(u,  \eta)$ depends on the approach used and $\delta_K$ is a non-negative stabilisation 
parameter, which can vary from one method to another. In this work, we {use the Streamline}
Upwind Petrov-Galerkin (SUPG) introduced by Brooks and Hughes \cite{BrooksHughes1982}, where 
$$ s_K(u, v)  = \int_K \Big( - \Delta u + \,\nabla ({\bf F} u), {\bf F} \cdot \nabla u \Big).$$
The stabilisation parameter $\delta_K$ is defined depending on the {cell's} Pecl\'et number, i.e., 
$$
P_{K} := \tfrac{ \| {\bf F}\|_{\infty, K} \, h_K}{2 \, \varepsilon},
$$
where 
$\| \cdot \|_{\infty, K}$ is the norm in $[L^{\infty}(K)]^{{d}}$, 
$h_K$ denotes the diameter of a finite element $K$, and 
${\bf F}$ prescribes the direction of the convective flow. Then, 
$$
\delta_K := \begin{cases} 
\delta_0 \, \tfrac{h_K}{\| {\bf F}\|_{\infty, K}}, \quad \mbox{if} \quad P_{K} > 1 \quad \mbox{(convection-dominated case), }\\ 
\delta_1 \, \tfrac{h^2_K}{\varepsilon},            \quad \mbox{if } \quad P_{K} \leq 1 \quad \mbox{(diffusion-dominated case), }
\end{cases}
$$
where $\delta_0, \delta_1 > 0$ are appropriately chosen constants. 
Set $\delta_K:= \tfrac{h_K}{2 \, | {\bf F} |}$, where
$| {\bf F} |$ stands for the magnitude of the vector. 
All the numerical results obtained in the section are carried out in Python and C++
in the framework of The FEniCS Project
\cite{Fenicsproject, LoggMardalWells2012}.

\begin{example}
\label{ex:example-40}
\rm 
Consider the one-dimensional classical example on the unit interval $\Omega = (0, 1)$:
\begin{equation}
- \varepsilon \, u_{xx} + a \, u_x = 0, \qquad u(0) = u(1) = 0, 
\end{equation}
with the solution $u = \tfrac{1 - e^{a\, x /  \varepsilon }}{ - e^{a\,/ \varepsilon}}$ known to have 
boundary layer in the neighborhood of the right part of the boundary $x = 1$. We select
$a = 1.0$ and $\varepsilon = 1e\minus2$.
Its approximation $v$ is reconstructed for two cases
with Lagrangian finite elements (FEs) of order 1 ($v \in \Pone$) and of order 2
($v \in \Ptwo$).
The correctness of the numerical solver is confirmed by the error order of convergence (e.o.c.)
for the uniform refinement procedure executed for both $v \in \Pone$ and $v \in \Ptwo$
(see Table \ref{tab:example-40-convergence-v-Pone-vPtwo}). The same results are illustrated in 
Figure \ref{fig:example-40-convergence-v-Pone-vPtwo}. In both cases, the initial mesh contains of 8 
elements and 9 nodes, and 
we proceed for 12 global refinement steps. 
As expected, 
the error behaves as $O(h)$ and $O(h^2)$ in the first and the second case, respectively. 
Here, 
the flux $\flux$ is reconstructed with Raviart--Thomas elements of the first order denoted by $\RTone$ 
as follows
\begin{align*}
RT_1(\mathcal{K}_h) := \{&\boldsymbol{y} \in H(\dvrg, K): 
\, \forall \, K \in \mathcal{K}_h \quad \exists \, \boldsymbol{q} \in (\Pone)^d, r \in \Pone
\quad \forall \, \boldsymbol{x} \in K: \, 
\boldsymbol{y}(\boldsymbol{x}) = \boldsymbol{q}(\boldsymbol{x}) 
+ \boldsymbol{x} \, r(\boldsymbol{x})   
\}.
\end{align*}
The auxiliary function $w$ for the minorant reconstruction is approximated by $\Pthree$-Lagrangian
FEs. 
\begin{table}[!t]
\centering
\footnotesize
\begin{tabular}{c|cc|cc|cc|cc|cc|cc}
& \multicolumn{6}{c|}{ $v \in \Pone$} & \multicolumn{6}{c}{$v \in \Ptwo$} \\
\midrule
ref. $\#$ & 
 $[ e ]_{\overline{\rm M}}$ & \parbox[c]{0.7cm}{\centering e.o.c. $([ e ]_{\overline{\rm M}})$} & 
 $\overline{\rm M}$      & \parbox[c]{0.7cm}{\centering e.o.c. $(\overline{\rm M})$} & 
 $\underline{\rm M}$      & \parbox[c]{0.7cm}{\centering e.o.c. $(\underline{\rm M})$} & 
 $[ e ]_{\overline{\rm M}}$ & \parbox[c]{0.7cm}{\centering e.o.c. $([ e ]_{\overline{\rm M}})$} & 
 $\overline{\rm M}$      & \parbox[c]{0.7cm}{\centering e.o.c. $(\overline{\rm M})$} & 
 $\underline{\rm M}$      & \parbox[c]{0.7cm}{\centering e.o.c. $(\underline{\rm M})$} \\
\midrule
 5 & 5.66e-01 &     0.31 & 1.29e+00 &   0.98 & 1.06e+00 &   0.85 & 2.71e-01 &     0.75 & 5.15e-01 &   1.27 & 4.54e-01 &   1.22 \\
 7 & 2.97e-01 &     0.92 & 3.24e-01 &   1.01 & 3.19e-01 &   1.00 & 5.64e-02 &     1.86 & 5.90e-02 &   1.90 & 5.86e-02 &   1.90 \\
9 & 7.94e-02 &     0.99 & 7.99e-02 &   1.00 & 7.97e-02 &   1.00 & 3.99e-03 &     1.99 & 4.01e-03 &   1.99 & 4.00e-03 &   1.99 \\
 11 & 1.99e-02 &     1.00 & 1.99e-02 &   1.00 & 1.99e-02 &   1.00 & 2.51e-04 &     2.00 & 2.51e-04 &   2.00 & 2.51e-04 &   2.00 \\
\end{tabular}
\caption{Example \ref{ex:example-40}. Error order of convergence (e.o.c.) 
for approximations $v \in \Pone$ and $\flux \in \RTone$ as well as 
$v \in \Ptwo$ and $\flux \in \RTone$ w.r.t. refinement steps.}
\label{tab:example-40-convergence-v-Pone-vPtwo}
\end{table}
\begin{figure}[!t]
	\centering
	\subfloat[$v \in \Pone$]{\includegraphics[width=6.1cm, trim={6.5cm 10.5cm 7cm 11cm}, clip]{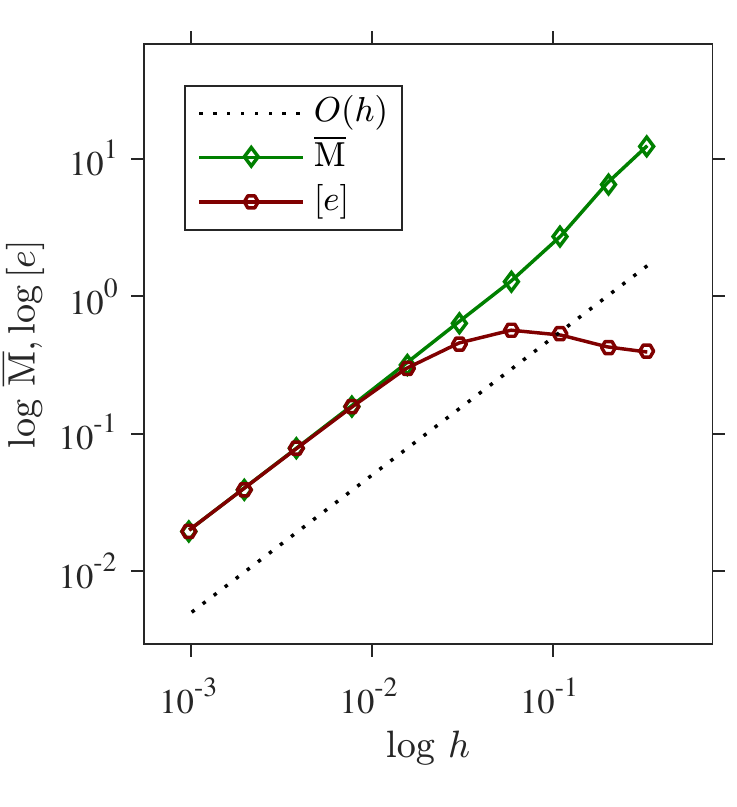}} \quad
	\subfloat[$v \in \Ptwo$]{\includegraphics[width=6.1cm, trim={6.5cm 10.5cm 7cm 11cm}, clip]{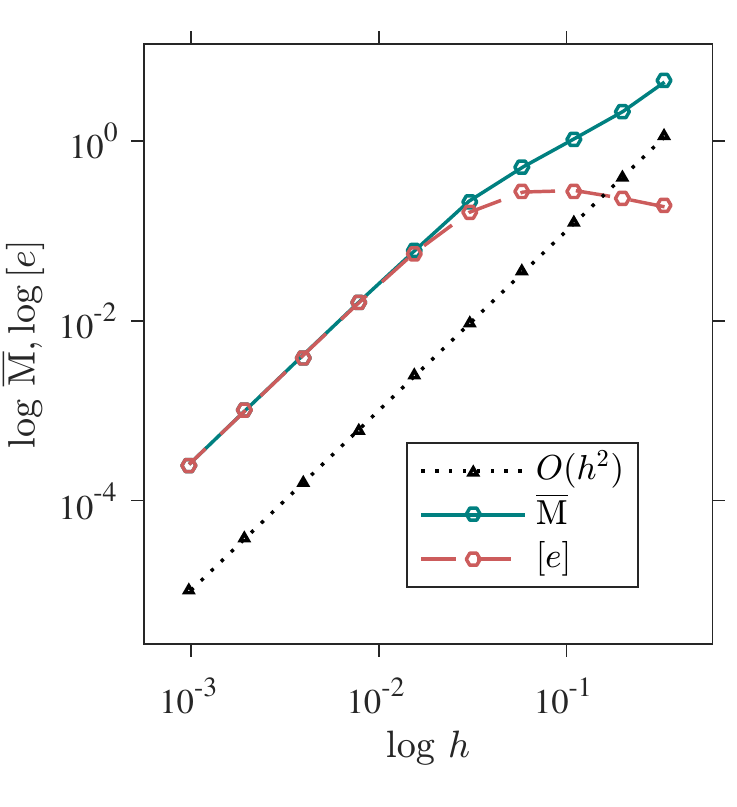}}
	\caption{Example \ref{ex:example-40}. Error order of convergence (e.o.c.) for 
	(a) $v \in \Pone$ and $\flux \in \RTone$ as well as 
	(b) $v \in \Ptwo$ and $\flux \in \RTone$ w.r.t. degrees of freedom (d.o.f.).}	
	\label{fig:example-40-convergence-v-Pone-vPtwo}	
\end{figure}
The values of the majorants and minorants as well as their efficiency indices w.r.t. 
several refinement steps are presented in Table \ref{tab:example-40-convergence-v-Pone-vPtwo}. 
It is easy to see that both upper and lower bounds stay rather sharp w.r.t. the error measures 
$[ e ]_{\overline{\rm M}}$ and $[ e ]_{\underline{\rm M}}$, respectively.
\begin{table}[!t]
\centering
\footnotesize
\begin{tabular}{c|c|ccc|ccc}
ref. $\#$ & d.o.f. $\#$  & 
$[ e ]_{\overline{\rm M}}$ & $\overline{\rm M}$      & $I_{\rm eff}(\overline{\rm M})$ & 
$[ e ]_{\underline{\rm M}}$ & $\underline{\rm M} $     & $I_{\rm eff}(\underline{\rm M})$ \\
\midrule
      0 &        9 &  8.08e-01 & 8.45e+00 &      10.46 &  2.81e+00 & 2.75e+00 &       0.98 \\
      2 &       33 &  5.09e-01 & 1.31e+00 &       2.56 &  6.36e-01 & 6.38e-01 &       1.00 \\
      4 &      129 &  1.57e-01 & 2.50e-01 &       1.59 &  1.59e-01 & 1.59e-01 &       1.00 \\
      6 &      513 &  3.98e-02 & 5.21e-02 &       1.31 &  3.99e-02 & 3.99e-02 &       1.00 \\
      8 &     2049 &  9.97e-03 & 1.15e-02 &       1.16 &  9.97e-03 & 9.97e-03 &       1.00 \\
     10 &     8193 &  2.49e-03 & 2.68e-03 &       1.07 &  2.49e-03 & 2.49e-03 &       1.00 \\
     12 &    32769 &  6.23e-04 & 6.44e-04 &       1.03 &  6.23e-04 & 6.23e-04 &       1.00 \\
\end{tabular}
\caption{Example \ref{ex:example-40}. Majorant, minorant, and corresponding efficiency indices.}
\label{tab:example-40-efficiency-indices}
\end{table}

\end{example}

\begin{example}
\label{ex:example-44}
\rm
Next, we consider the two-dimensional example defined on the unit square domain $\Omega = (0, 1)^2$:
\begin{equation}
- \varepsilon \, \Delta \, u + {\bf F} \cdot \nabla \, u = x, \qquad 
u = 0 \quad \mbox{on}\quad \partial \Omega, 
\end{equation}
where $\varepsilon = 1e\minus3$ and ${\bf F} = (-x, 0)^{\rm T}$.
In this case, the exact solution is not known. However, by following the analysis presented in 
\cite{RoosStynesTobiska1996}, one can predict the location of the boundary layer. 
The approximation $v$ is reconstructed by $\Pone$-elements, whereas 
the majorant is reconstructed with an auxiliary function $\flux \in \RTone$ and 
the minorant
with an auxiliary function $w \in \Pthree$.

\begin{figure}[!t]
	\centering
	\captionsetup[subfigure]{oneside, margin={0.9cm,0cm}}
	\subfloat[${[ e ]}_{\overline{\rm M}}$ and $\overline{\rm M}$]{
	\includegraphics[scale=0.8]{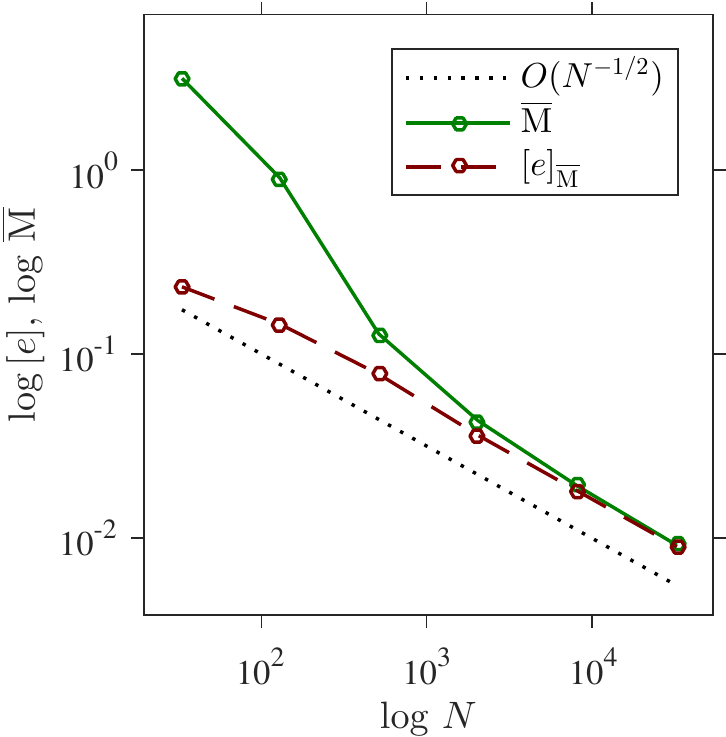}} \quad
	\subfloat[${[ e ]}_{\underline{\rm M}}$ and $\underline{\rm M}$]{
	\includegraphics[scale=0.8]{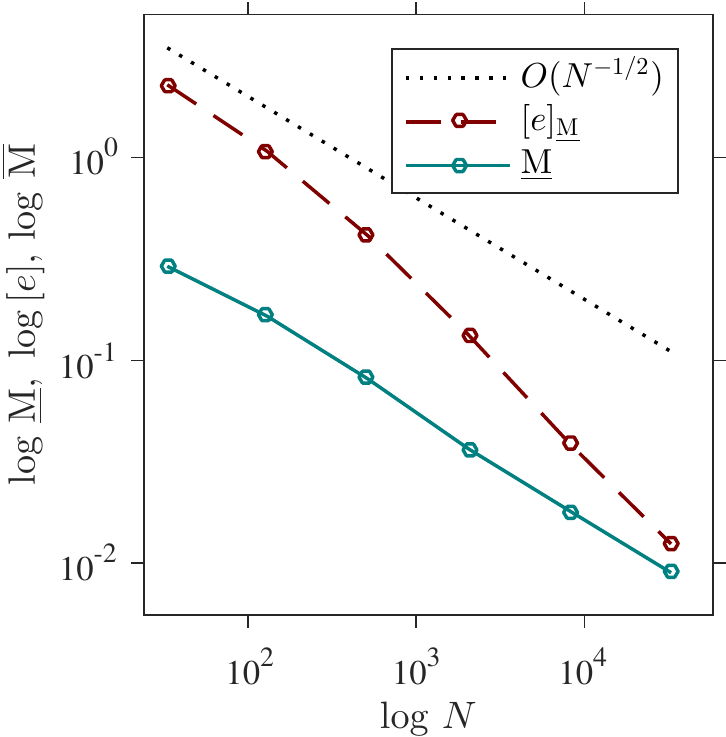}}
	\caption{Example \ref{ex:example-44}. 
	Convergence of 
	(a) ${[ e ]}_{\overline{\rm M}}$ and $\overline{\rm M}$ as well as 
	(b) ${[ e ]}_{\underline{\rm M}}$ and $\underline{\rm M}$.}	
	\label{fig:example-44-convergence-maj-min}	
\end{figure}

\begin{table}[!t]
\centering
\footnotesize
\begin{tabular}{c|c|ccc|ccc}
ref. $\#$ & d.o.f. $\#$    & 
$[ e ]_{\overline{\rm M}}$ & $\overline{\rm M}$      & $I_{\rm eff}(\overline{\rm M})$ & 
$[ e ]_{\underline{\rm M}}$ & $\underline{\rm M} $     & $I_{\rm eff}(\underline{\rm M})$ \\
\midrule
      0 &       25 &  2.32e-01 & 3.15e+00 &      13.57 &  2.29e+00 & 2.91e-01 &       0.13 \\
      2 &      289 &  7.77e-02 & 1.29e-01 &       1.66 &  4.13e-01 & 8.12e-02 &       0.20 \\
      4 &     4225 &  1.80e-02 & 1.92e-02 &       1.06 &  3.88e-02 & 1.81e-02 &       0.47 \\
      6 &    66049 &  4.47e-03 & 4.49e-03 &       1.01 &  4.95e-03 & 4.46e-03 &       0.90 \\
\end{tabular}
\caption{Example \ref{ex:example-44}. Majorant, minorant, and corresponding efficiency indices.}
\label{tab:example-44-efficiency-indices}
\end{table}

Let us first consider the numerical experiments 
obtained by globally refining 
the mesh. 
The order of convergence of $\underline{\rm M}$ and $\overline{\rm M}$ is confirmed in Figure 
\ref{fig:example-44-convergence-maj-min}. Due to their high quantitative efficiency, 
the lines of the error estimates and corresponding errors almost coincide on the plots.  
Table \ref{tab:example-44-efficiency-indices} illustrates values of both minorant and majorant 
as well as their efficiency indices w.r.t. the error measures that they bound.

Next, we test the performance of an adaptive strategy of the refinement. It is based on  
the local distribution of the indicator, following from the majorant functional. 
We choose a widely-used bulk marking criterion for selecting elements to 
be refined denoted by ${\mathds{M}}_{{\rm BULK}}(\theta)$, see \cite{Doerfler1996}.
Here,  
the bulk parameter is set to $\theta = 0.3$. Table 
\ref{tab:example-44-efficiency-indices-adaptive} provides the 
values for the efficiency indices of the majorant and minorant,
which are getting sharper as the refinement steps proceed. The 
evolution of the meshes 
is illustrated 
in Figure \ref{fig:example-44-meshes}. Here, the 
left column presents the meshes constructed on refinement steps 3--4, whereas the right column
illustrates the distribution of the local error values on the elements of the mesh. The elements 
with the highest and lowest errors are indicated by the range of colours between red and blue, 
respectively. From the meshes, one can also see that refinement based on the majorant 
resolves the boundary layers at $x = 0$.

\begin{table}[!t]
\centering
\footnotesize
\begin{tabular}{c|c|ccc|ccc}
ref. $\#$ & d.o.f. $\#$    & 
$[ e ]_{\overline{\rm M}}$ & $\overline{\rm M}$      & $I_{\rm eff}(\overline{\rm M})$ & 
$[ e ]_{\underline{\rm M}}$ & $\underline{\rm M} $     & $I_{\rm eff}(\underline{\rm M})$ \\
\midrule
      0 &       25 &  4.06e-02 & 4.66e-02 &       1.15 &  4.06e-02 & 4.06e-02 &       1.00 \\
      2 &      109 &  1.50e-02 & 1.56e-02 &       1.04 &  1.50e-02 & 1.50e-02 &       1.00 \\
      4 &      576 &  6.53e-03 & 6.66e-03 &       1.02 &  6.53e-03 & 6.53e-03 &       1.00 \\
      6 &     2870 &  2.99e-03 & 3.02e-03 &       1.01 &  2.99e-03 & 2.99e-03 &       1.00 \\
      8 &    14527 &  1.35e-03 & 1.35e-03 &       1.00 &  1.35e-03 & 1.35e-03 &       1.00 \\
\end{tabular}
\caption{Example \ref{ex:example-44}. Majorant, minorant, and corresponding efficiency indices for an 
adaptive refinement strategy.}
\label{tab:example-49-efficiency-indices-adaptive}
\end{table}

\begin{figure}[!t]
	\centering
	\captionsetup[subfigure]{oneside, labelformat=empty}
	{
	\subfloat[ref. \# 3:]{
	\includegraphics[height=4.5cm, trim={2cm 1cm 2cm 1cm}, clip]{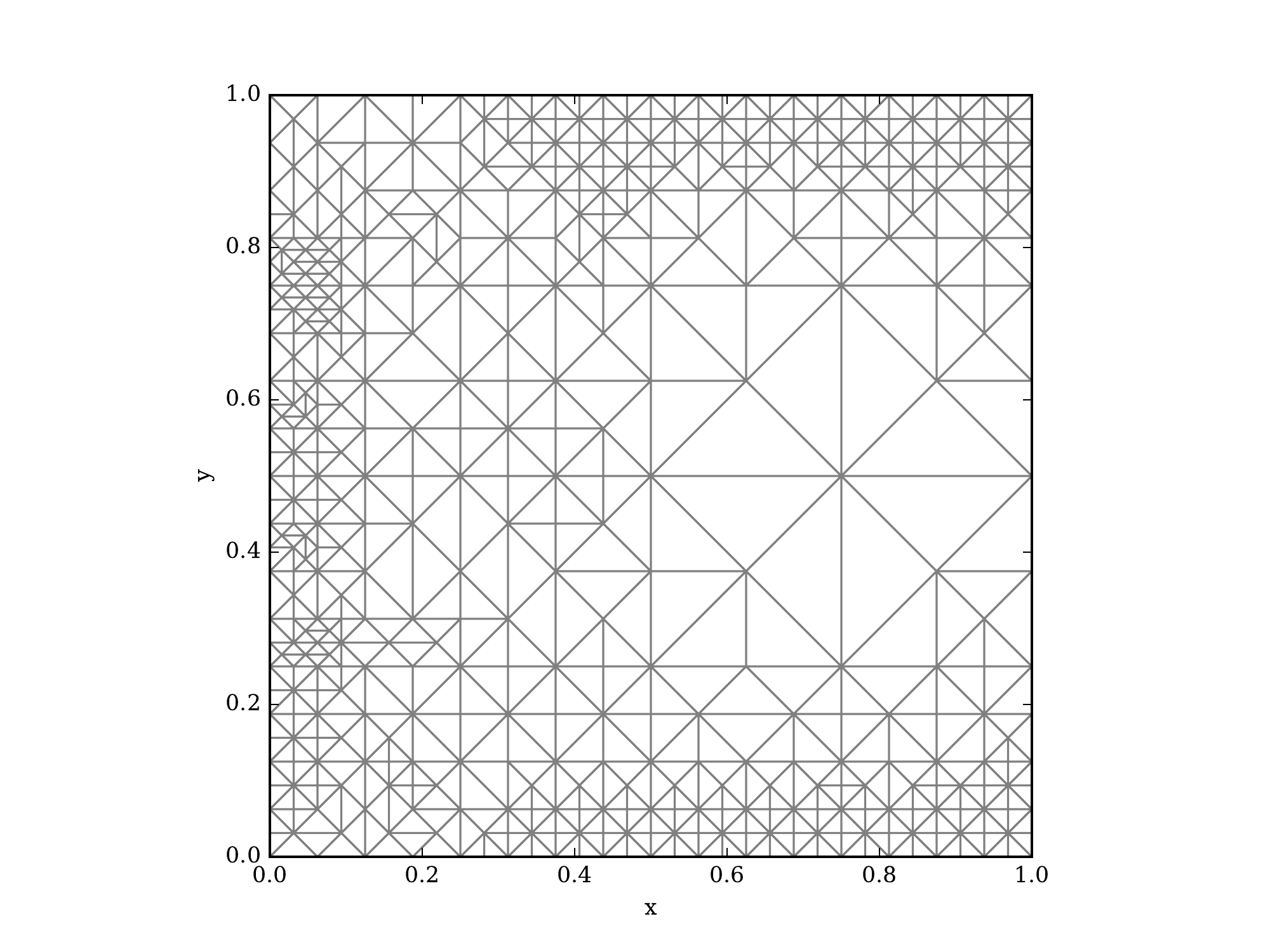}}
	\qquad
	\subfloat[ref. \# 3:]{
	\includegraphics[height=4.5cm, trim={1.5cm 1cm 1.5cm 1cm}, clip]{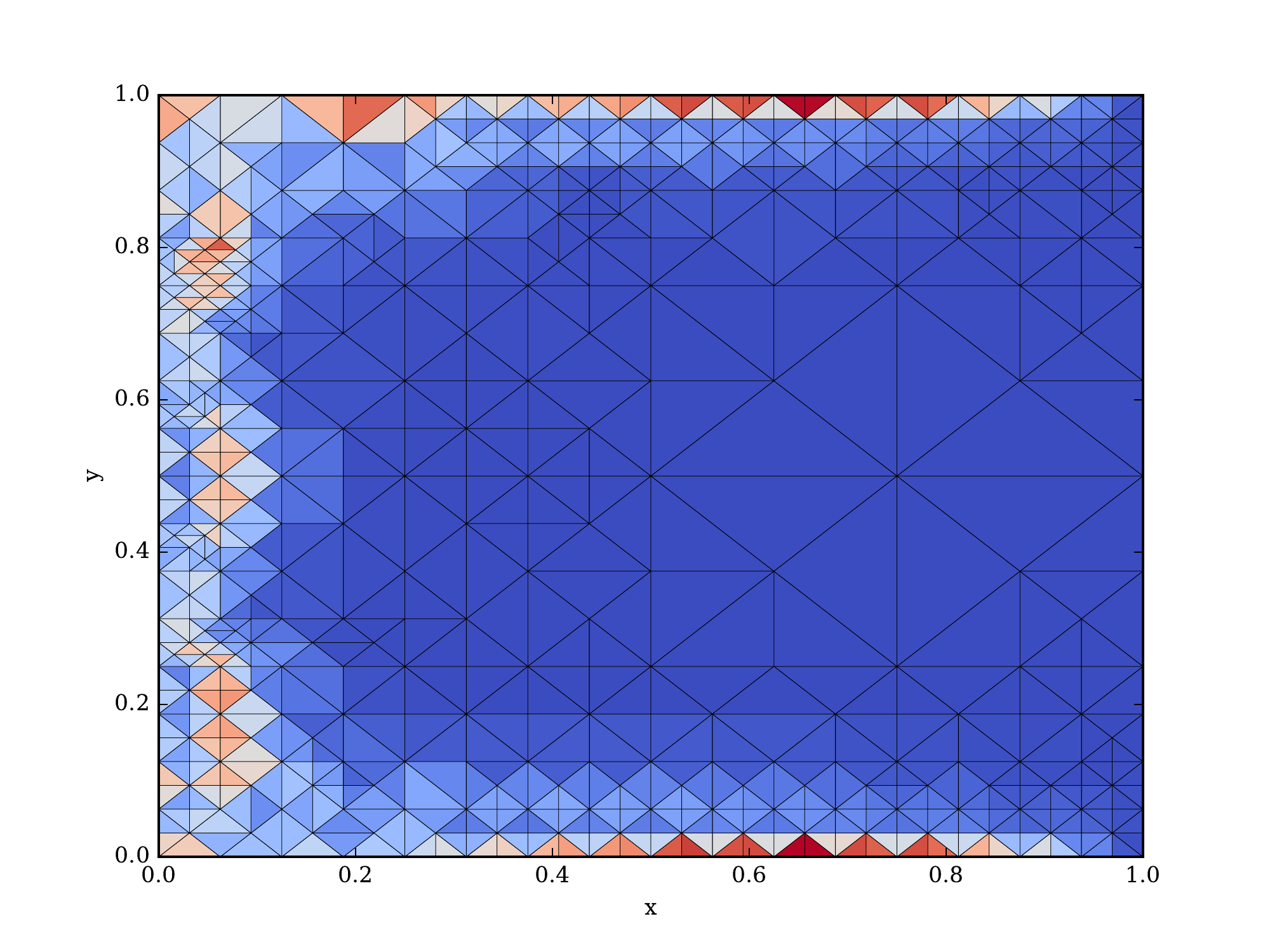}} 
	}
	{
	\subfloat[ref. \# 4:]{
	\includegraphics[height=4.5cm, trim={2cm 1cm 2cm 1cm}, clip]{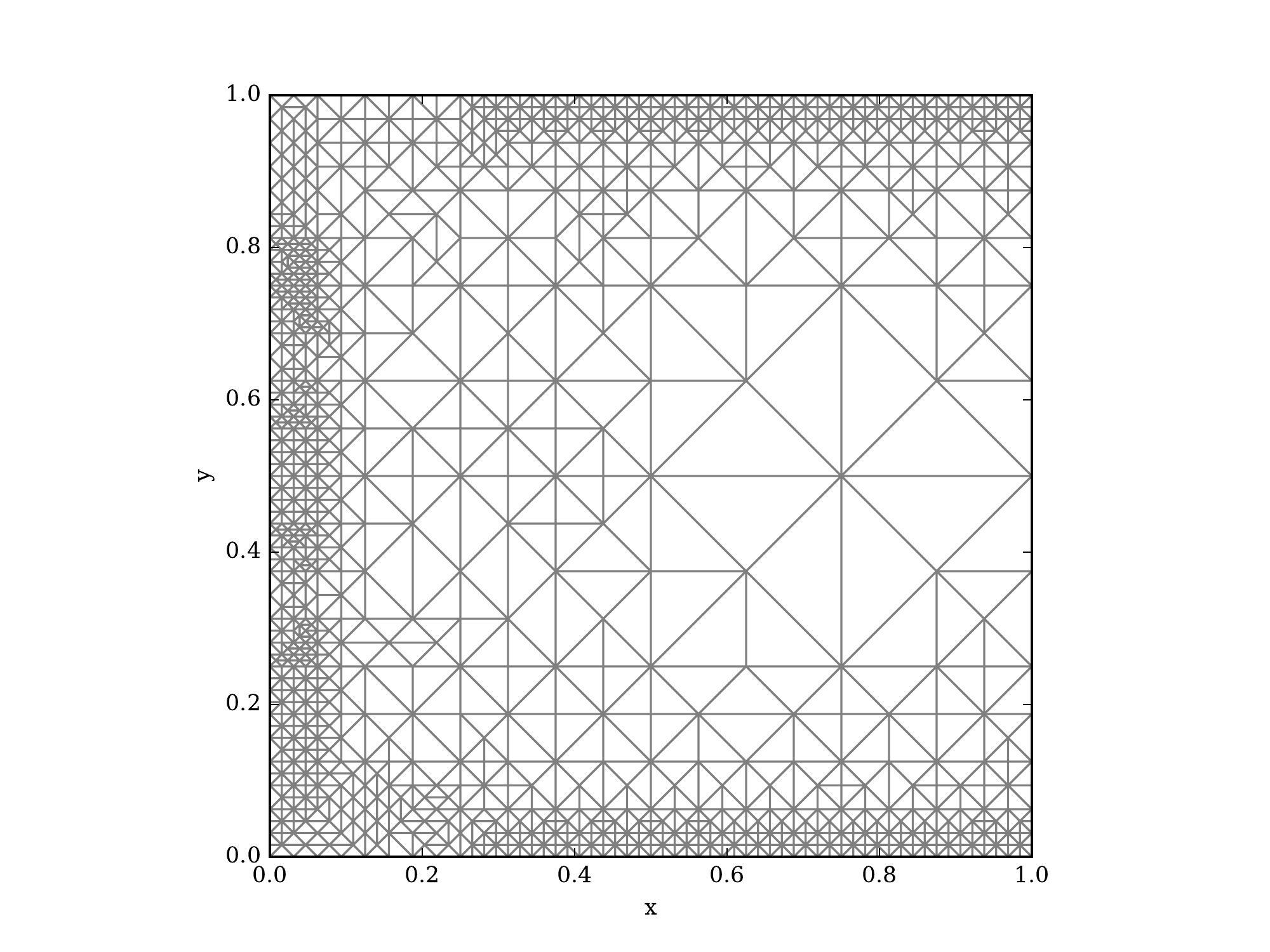}}
	\qquad
	\subfloat[ref. \# 4:]{
	\includegraphics[height=4.5cm, trim={1.5cm 1cm 1.5cm 1cm}, clip]{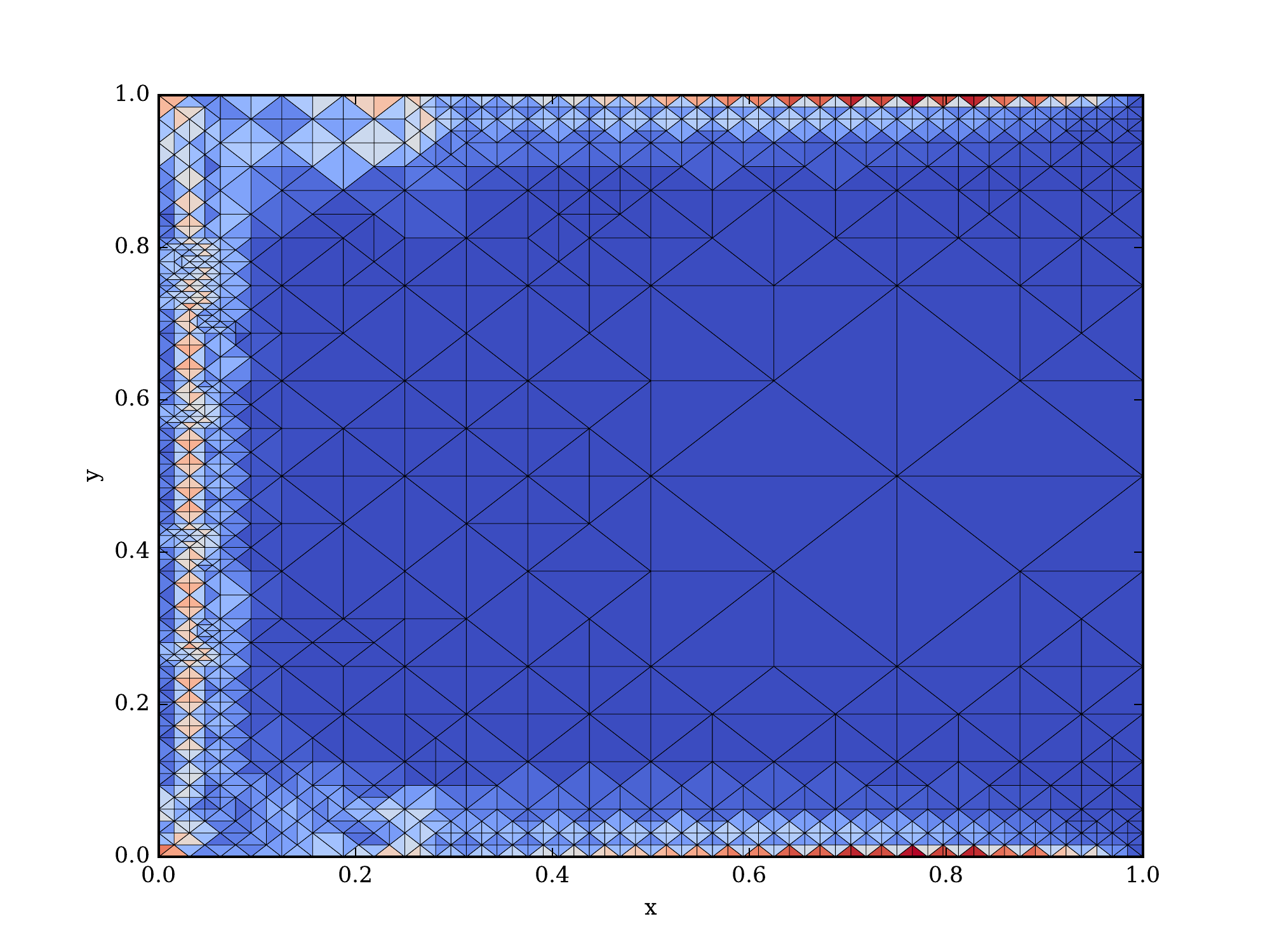}} 
	}
	\caption{Example \ref{ex:example-44}. 
	Meshes generated by the adaptive refinement strategy and the distribution of 
	the local error (right column) with ${\mathds{M}}_{{\rm BULK}}(0.3)$.}
	\label{fig:example-44-meshes}
\end{figure}	

\end{example}

\begin{example}
\label{ex:example-49}
\rm
Next, we consider yet another convection-dominated example on the unit square $\Omega = (0, 1)^2$:
\begin{equation}
- \dvrg ( A  \nabla \, u) + \nabla ({\bf F} \, u) = f, \qquad 
u = 0 \quad \mbox{on} \quad \partial \Omega, 
\end{equation}
where $A = \begin{bmatrix}10 & 0 \\ 0 & 1 \end{bmatrix}$, 
$f = x + y$, and 
${\bf F} = \big(-\cos(\tfrac{\pi}{3})\, x, \, \sin(\tfrac{\pi}{3}) \,y\big)^{\rm T}$. Again, let the 
approximation $v$ be reconstructed with $\Pone$-elements, $\flux$ by $\RTone$-elements, and $w \in \Pthree$ 
(the optimal error convergence is presented 
in Figure \ref{fig:example-49-convergence-v-Pone}). 
The approximate solution 
obtained with an adaptive refinement is illustrated in 
Figure \ref{fig:example-49-approximate-solution-adapt-ref-5-c}, where
one can see that the solution convects towards the boundary $y = 1$.

\begin{figure}[!t]
	\centering
	\subfloat[$v$, adaptive ref. $\#$ 6]{
	\includegraphics[width=7.8cm, trim={2cm 1.5cm 2cm 2cm}, clip]{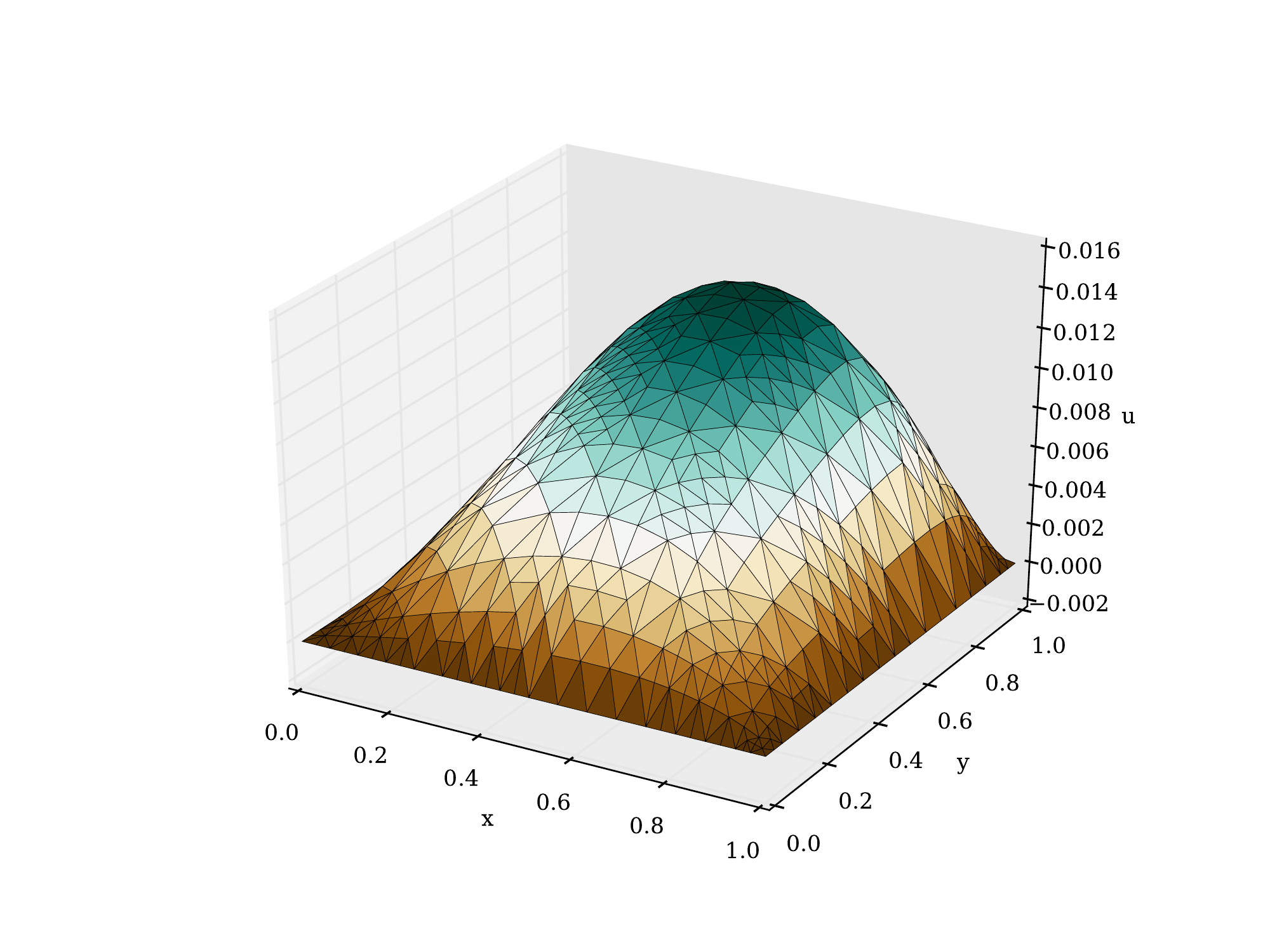}
	\label{fig:example-49-approximate-solution-adapt-ref-5-c}}
	\subfloat[$v$, adaptive ref. $\#$ 6]{
	\includegraphics[width=6.0cm]{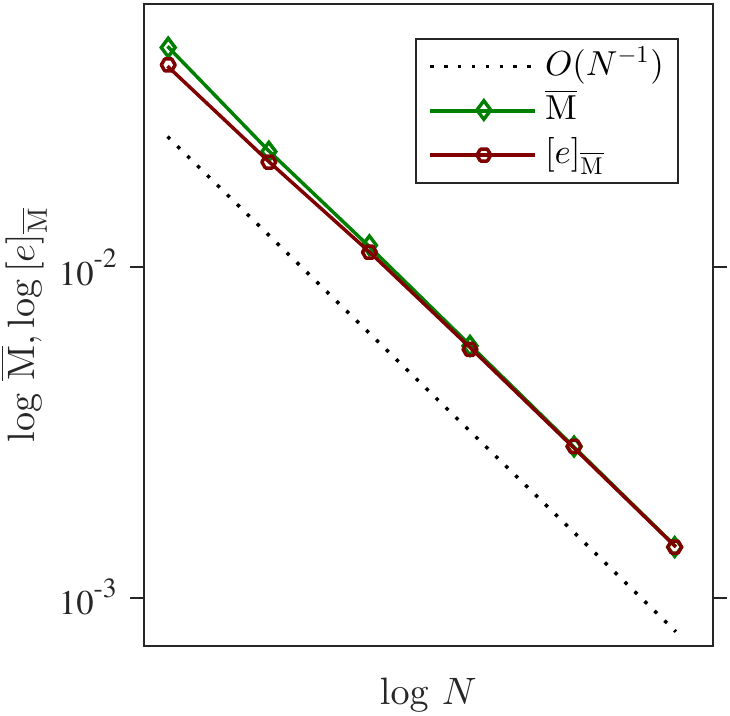}
	\label{fig:example-49-convergence-v-Pone}}
	\caption{Example \ref{ex:example-49}. 
	(a) Approximate solution on the adaptive refinement step $\#$ 4.
	(b) Optimal order of convergence together with error and majorant order of convergence.}	
	\label{fig:ex:example-49}	
\end{figure}
 
By the current example, we aim to show that even for the anisotropic $A$ and heterogeneous 
${\bf F}$ the performance of the error estimates remains rather sharp. In particular, 
Table \ref{tab:example-44-efficiency-indices-adaptive} confirms that the efficiency indices of the 
minorant and the majorant lie close to $1$ as the refinements evolve. 

\begin{table}[!t]
\centering
\footnotesize
\begin{tabular}{c|c|ccc|ccc}
ref. $\#$ & d.o.f. $\#$    & 
$[ e ]_{\overline{\rm M}}$ & $\overline{\rm M}$      & $I_{\rm eff}(\overline{\rm M})$ & 
$[ e ]_{\underline{\rm M}}$ & $\underline{\rm M} $     & $I_{\rm eff}(\underline{\rm M})$ \\
\midrule
      0 &       25 &  2.32e-01 & 3.15e+00 &      13.57 &  2.29e+00 & 2.91e-01 &       0.13 \\
      2 &      101 &  1.18e-01 & 4.47e-01 &       3.80 &  9.98e-01 & 1.31e-01 &       0.13 \\
      4 &      469 &  4.05e-02 & 5.19e-02 &       1.28 &  1.41e-01 & 4.11e-02 &       0.29 \\
      6 &     2319 &  1.27e-02 & 1.36e-02 &       1.07 &  1.95e-02 & 1.27e-02 &       0.65 \\
      8 &    11753 &  5.43e-03 & 5.51e-03 &       1.02 &  6.04e-03 & 5.43e-03 &       0.90 \\
\end{tabular}
\caption{Example \ref{ex:example-49}. Majorant, minorant, and corresponding efficiency indices for an 
adaptive refinement strategy.}
\label{tab:example-44-efficiency-indices-adaptive}
\end{table}

\end{example}

\subsection{Static convection-dominated problems}

Next, we discuss the following \textbf{set of example problems}: consider the 
reaction-convection-diffusion equation 
\begin{equation}
\label{eq:generalized-statement-static}
- \varepsilon \, \Delta \, u + {\bf a} \cdot \nabla \, u + \lambda\, u = 0
\end{equation}
(which basically yields from \eqref{eq:fokker-planck-state-equation}--\eqref{eq:dual-state-equation}) 
with homogeneous Dirichlet BCs. However, in contrast to the 
previous examples, we assume now that reaction and 
convection are independent of each other. 

\begin{example}
\label{ex:example-50-sicom-paper}
\rm
Assume that $d = 2$, ${\bf a} = \big(1, 0 \big)^{\rm T}$, $\varepsilon = 1e\minus{3}$, $\lambda = 1$, 
and the exact solution is defined by the function 
$u = \bigg(x + 
\tfrac{e^{(x - 1)/{\varepsilon}} - e^{- {1}/{\varepsilon}} }{e^{- {1}/{\varepsilon}} - 1}\bigg) 
\, y \, (1 - y)$ 
with a boundary layer in the neighbourhood of $x = 1$. Its approximation $v$ is reconstructed by 
${\rm P}_1$ elements (see Figure \ref{fig:example-50-ref-sol-5}). 
Since the uniform refinement does not provide expected orders of convergence, we apply an adaptive 
refinement with bulk marking criterion ${\mathds{M}}_{{\rm BULK}}(\theta)$. 
We choose two different parameters, i.e., $\theta = 0.3$ and 
$\theta = 0.4$, to confirm optimal convergence order for $\overline{\rm M}$ 
(see Table \ref{tab:example-50-sicom-paper-convergence-rates} 
and Figure \ref{fig:example-50-sicom-paper-convergence-rates}). 
The values of the total majorants and 
minorants as well as corresponding efficiency indices are presented in Table 
\ref{tab:eq:example-50-sicom-paper-efficiency-indices} w.r.t. several refinement steps. It is easy to 
see that even though the efficiency indices are high on the initial steps, both majorant and minorant 
become rather sharp towards the $9$-th refinement step.
The evolution of the meshes corresponding to the different $\theta$ parameters is presented in Figure 
\ref{fig:example-50-sicom-paper-mesh-bulk-30-40}.

\begin{figure}[!t]
	\centering
	\subfloat[$v \in P^1$, ref. $\#$ 5]{
	\includegraphics[width=8.3cm, trim={2cm 1.5cm 2cm 2cm}]{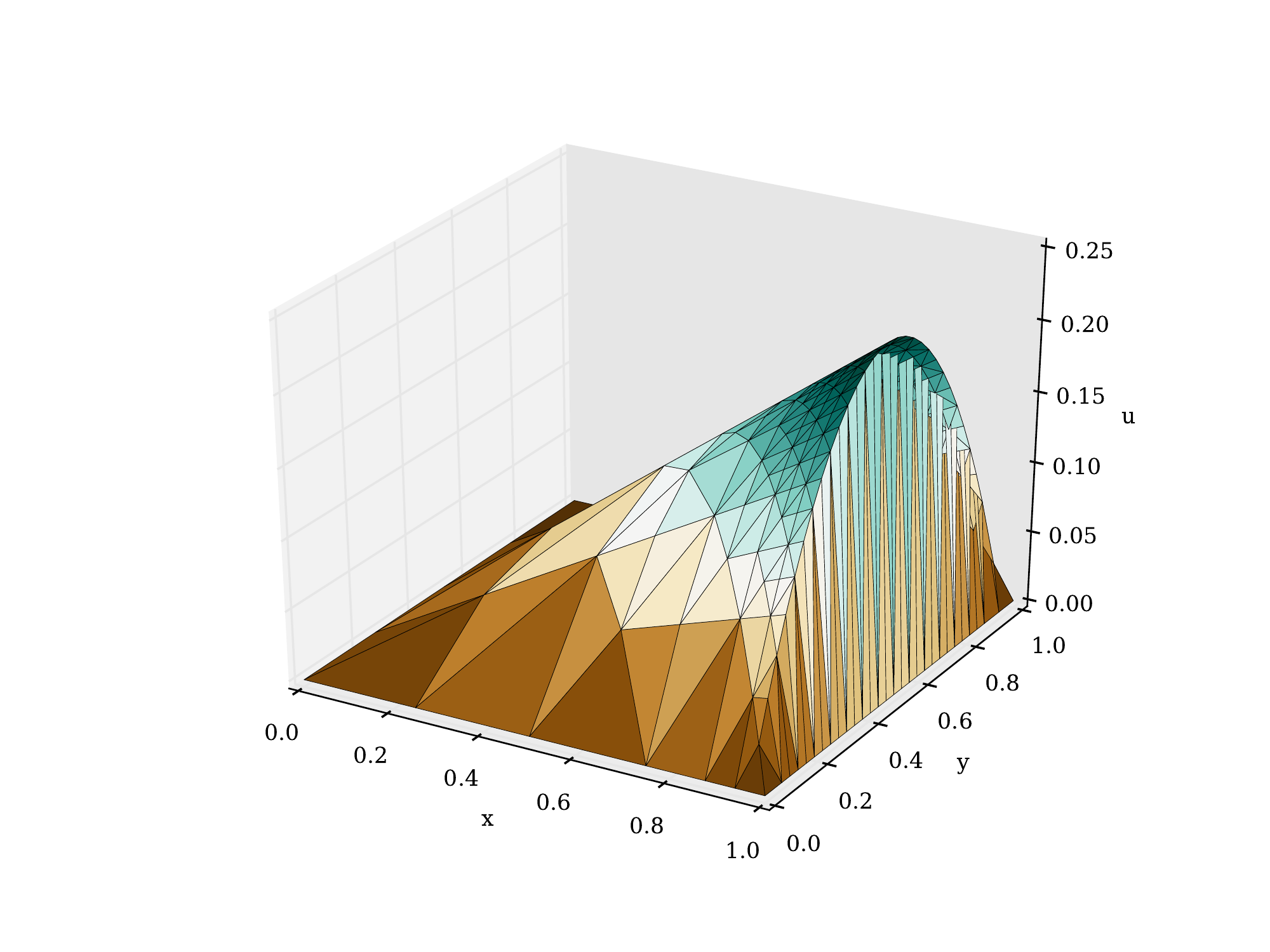}
	\label{fig:example-50-ref-sol-5}}
	\quad
	\subfloat[adaptive refinement]{
	\includegraphics[width=6.0cm]{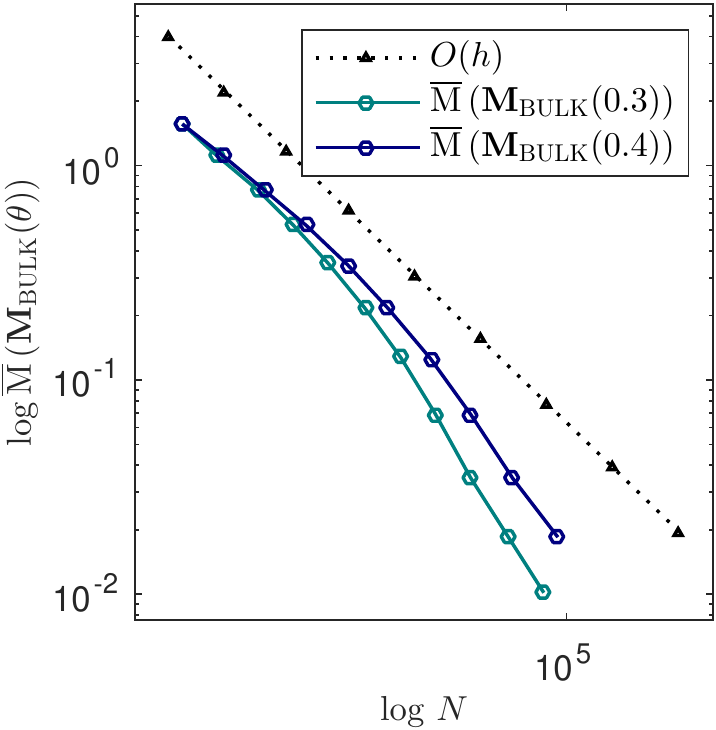}
	\label{fig:example-50-sicom-paper-convergence-rates}}
	\caption{Example \ref{ex:example-50-sicom-paper}. 
	(a) The approximate solution on the mesh (1089 nodes, 2048 elements).
	(b) Majorant's order of convergence for bulk parameters 
	$\theta = 0.3$ and $\theta = 0.4$.}
\end{figure}

\begin{table}[!t]
\centering
\footnotesize
\begin{tabular}{c|cc|cc|cc}
 $N^{-1/d}$ & 
 $[ e ]$ & \parbox[c]{1.7cm}{\centering e.o.c. $([ e ])$} & 
 $\overline{\rm M}$      & \parbox[c]{1.7cm}{\centering e.o.c. $(\overline{\rm M})$} & 
 $\underline{\rm M}$      & \parbox[c]{1.7cm}{\centering e.o.c. $(\underline{\rm M})$} \\
\midrule
\multicolumn{7}{c}{ $\theta = 0.3$} \\
\midrule
    191 & 7.66e-02 &   4.57 & 5.26e-01 &   1.17 & 4.83e-01 &   1.05 \\
    808 & 1.08e-01 &   2.27 & 2.17e-01 &   1.44 & 2.13e-01 &   1.41 \\
   3478 & 6.01e-02 &   1.93 & 6.96e-02 &   1.88 & 6.91e-02 &   1.86 \\
  15724 & 1.84e-02 &   1.66 & 1.86e-02 &   1.63 & 1.85e-02 &   1.65 \\
  32918 & 1.01e-02 &   1.03 & 1.02e-02 &   0.98 & 1.01e-02 &   1.02 \\
  74789 & 6.60e-03 &   0.92 & 6.79e-03 &   0.78 & 6.63e-03 &   0.91 \\
\midrule
\multicolumn{7}{c}{ $\theta = 0.4	$} \\
\midrule
    244 & 7.65e-02 &   3.86 & 5.26e-01 &   1.01 & 4.84e-01 &   0.89 \\
   1290 & 1.08e-01 &   1.95 & 2.15e-01 &   1.24 & 2.12e-01 &   1.22 \\
   7056 & 5.99e-02 &   1.62 & 6.91e-02 &   1.58 & 6.86e-02 &   1.57 \\
  16822 & 3.39e-02 &   1.37 & 3.47e-02 &   1.35 & 3.47e-02 &   1.36 \\
\midrule
\end{tabular}
\caption{Example \ref{ex:example-50-sicom-paper}. 
Error order of convergence for approximations $v \in \Pone$ and $\flux \in \RTone$ for 
bulk parameters $\theta = 0.3$ and $\theta = 0.4$.}
\label{tab:example-50-sicom-paper-convergence-rates}
\end{table}

\begin{table}[!t]
\centering
\footnotesize
\begin{tabular}{c|c|ccc|ccc}
ref. $\#$ & d.o.f. $\#$    & 
${[ e ]}_{\overline{\rm M}}$   & $\overline{\rm M}$    & $I_{\rm eff}(\overline{\rm M})$ & 
${[ e ]}_{\underline{\rm M}}$ & $\underline{\rm M}$  & $I_{\rm eff}(\underline{\rm M})$ \\
\midrule
\multicolumn{8}{c}{ $\theta = 0.3$ } \\
\midrule
      1 &       45 &  5.19e-02 & 1.11e+00 &      21.41 &  1.09e+00 & 6.53e-01 &       0.60 \\
      3 &      191 &  7.66e-02 & 5.26e-01 &       6.86 &  5.22e-01 & 4.83e-01 &       0.93 \\
      5 &      808 &  1.08e-01 & 2.17e-01 &       2.02 &  2.20e-01 & 2.13e-01 &       0.97 \\
      7 &     3478 &  6.01e-02 & 6.96e-02 &       1.16 &  7.01e-02 & 6.91e-02 &       0.99 \\
      9 &    15724 &  1.84e-02 & 1.86e-02 &       1.01 &  1.86e-02 & 1.85e-02 &       1.00 \\
     11 &    74789 &  6.60e-03 & 6.79e-03 &       1.03 &  6.71e-03 & 6.63e-03 &       0.99 \\
\midrule
\multicolumn{8}{c}{ $\theta = 0.4$ } \\
\midrule
      1 &       51 &  5.25e-02 & 1.11e+00 &      21.09 &  1.10e+00 & 6.53e-01 &       0.59 \\
      3 &      244 &  7.65e-02 & 5.26e-01 &       6.87 &  5.20e-01 & 4.84e-01 &       0.93 \\
      5 &     1290 &  1.08e-01 & 2.15e-01 &       2.00 &  2.15e-01 & 2.12e-01 &       0.98 \\
      7 &     7056 &  5.99e-02 & 6.91e-02 &       1.15 &  6.95e-02 & 6.86e-02 &       0.99 \\
      9 &    42827 &  1.83e-02 & 1.85e-02 &       1.01 &  1.84e-02 & 1.84e-02 &       1.00 \\
\end{tabular}
\caption{Example \ref{ex:example-50-sicom-paper}. 
Majorant, minorant, and corresponding efficiency indices.}
\label{tab:eq:example-50-sicom-paper-efficiency-indices}
\end{table}

\begin{figure}[!t]
\centering
	\subfloat[ref. $\#$ 2, $\theta$ = 0.3]{\includegraphics[width=5cm, trim={3cm 1cm 3cm 1cm}, clip]{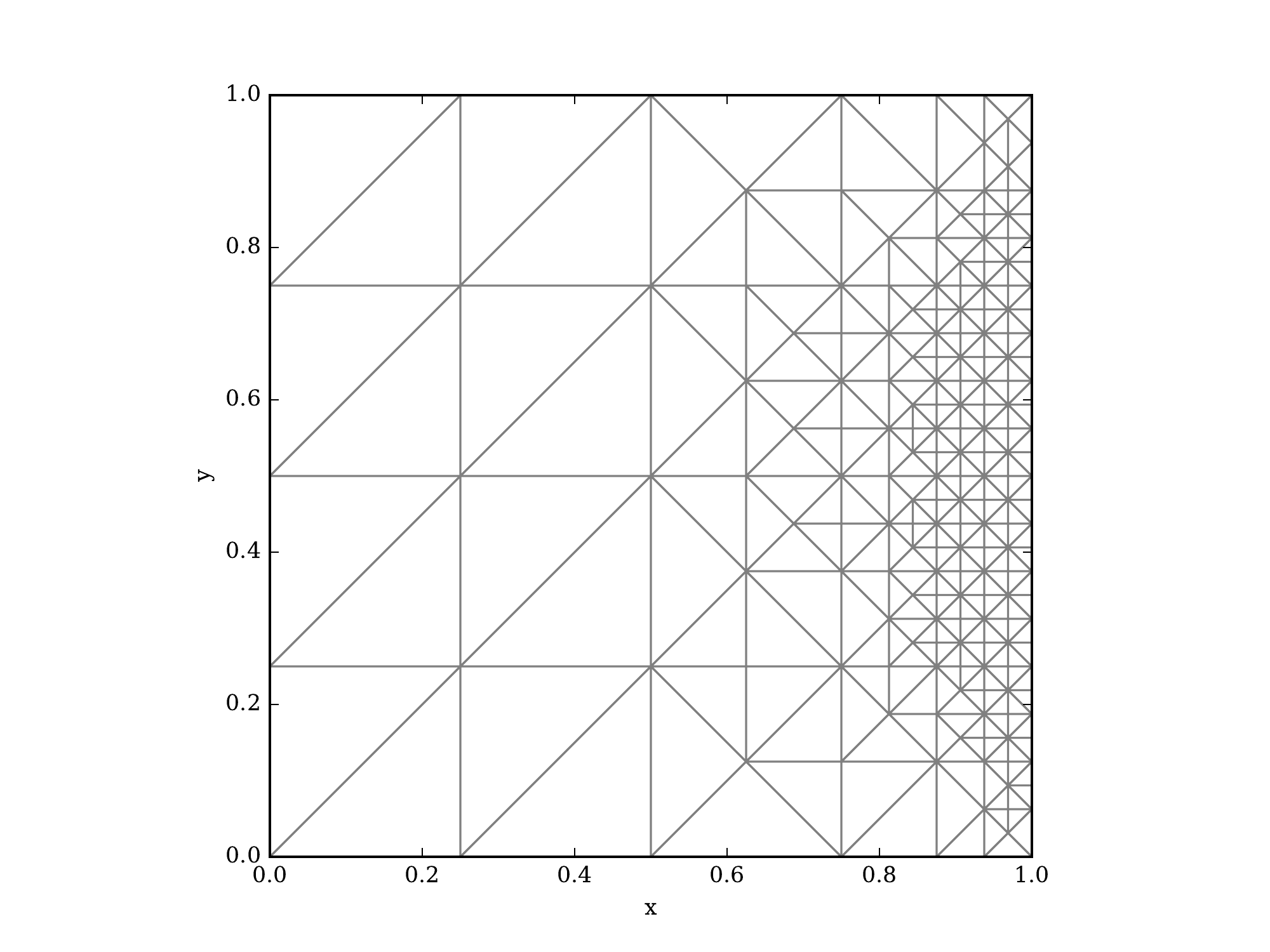}}\qquad
	\subfloat[ref. $\#$ 2, $\theta$ = 0.4]{\includegraphics[width=5cm, trim={3cm 1cm 3cm 1cm}, clip]{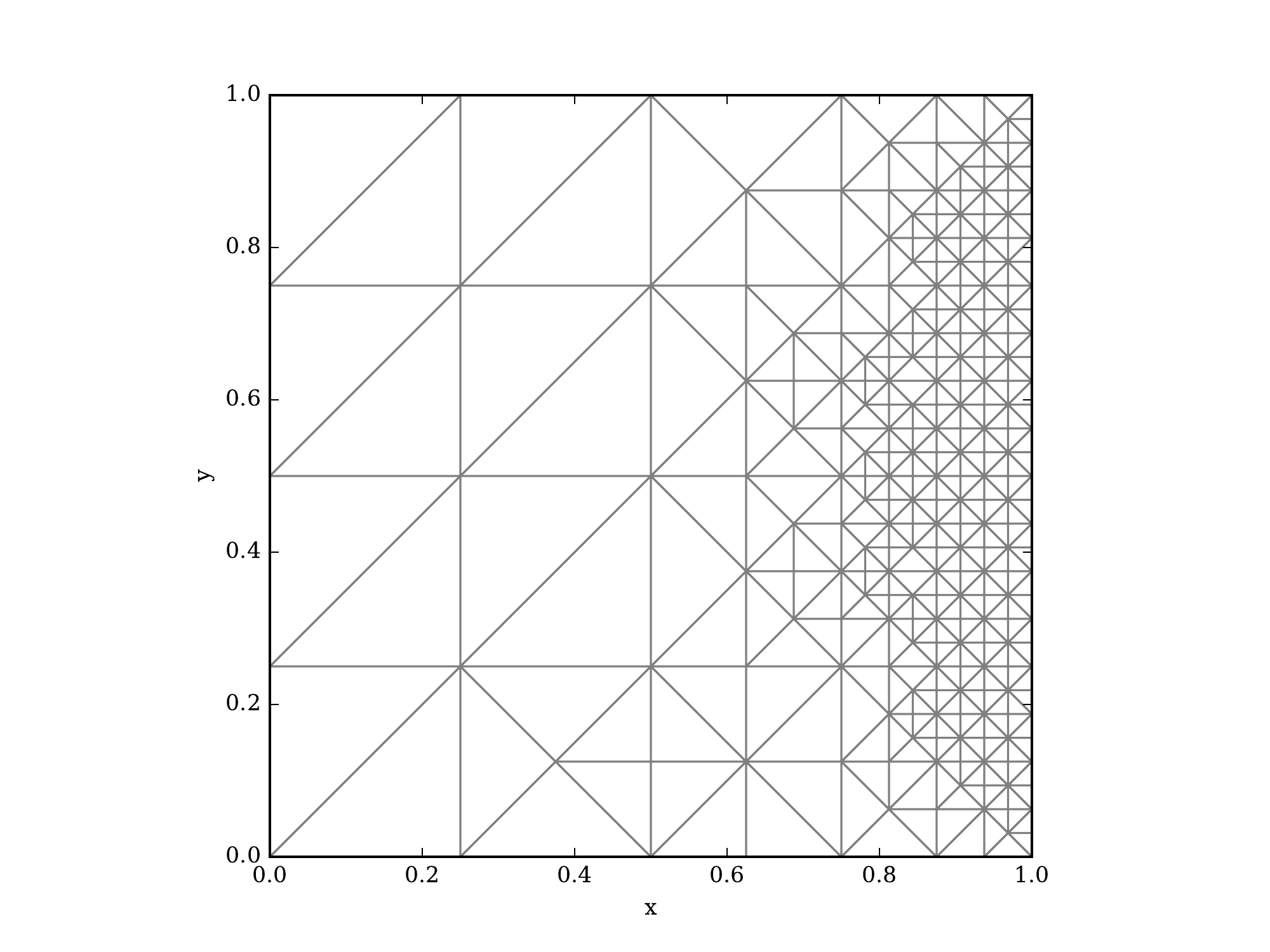}}\\
	\subfloat[ref. $\#$ 4, $\theta$ = 0.3]{\includegraphics[width=5cm, trim={3cm 1cm 3cm 1cm}, clip]{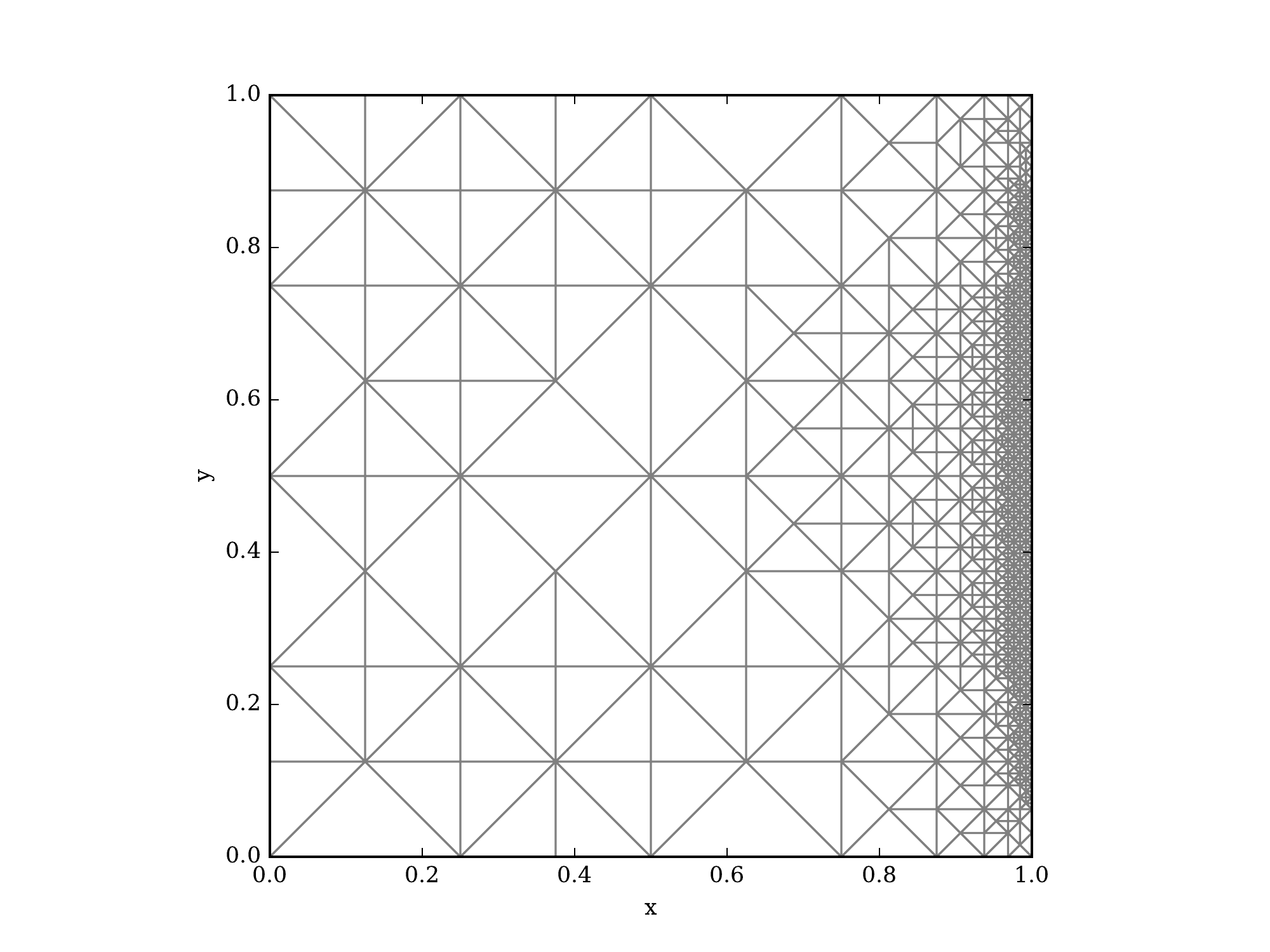}}\qquad
	\subfloat[ref. $\#$ 4, $\theta$ = 0.4]{\includegraphics[width=5cm, trim={3cm 1cm 3cm 1cm}, clip]{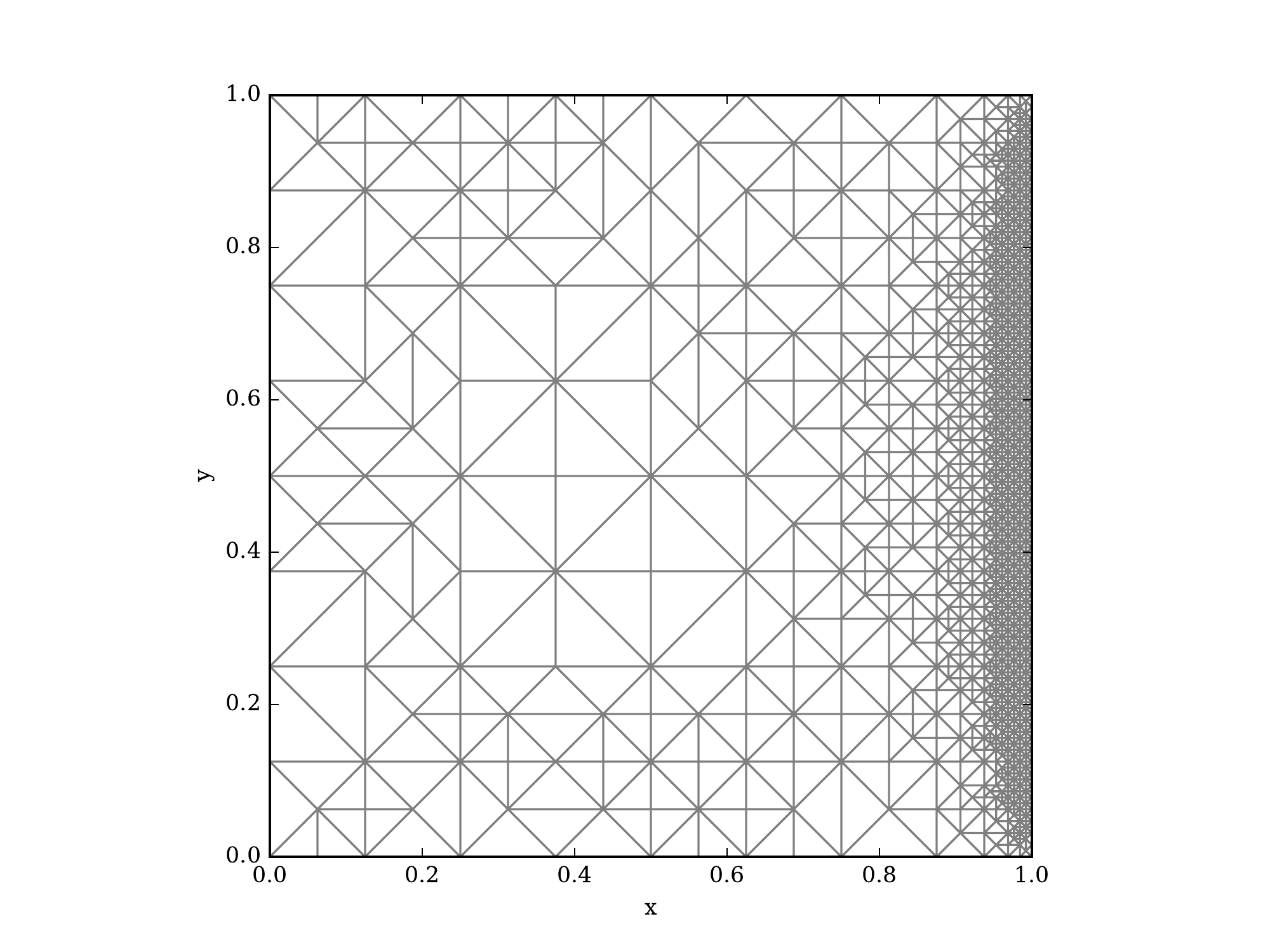}}\\
	\subfloat[ref. $\#$ 6, $\theta$ = 0.3]{\includegraphics[width=5cm, trim={3cm 1cm 3cm 1cm}, clip]{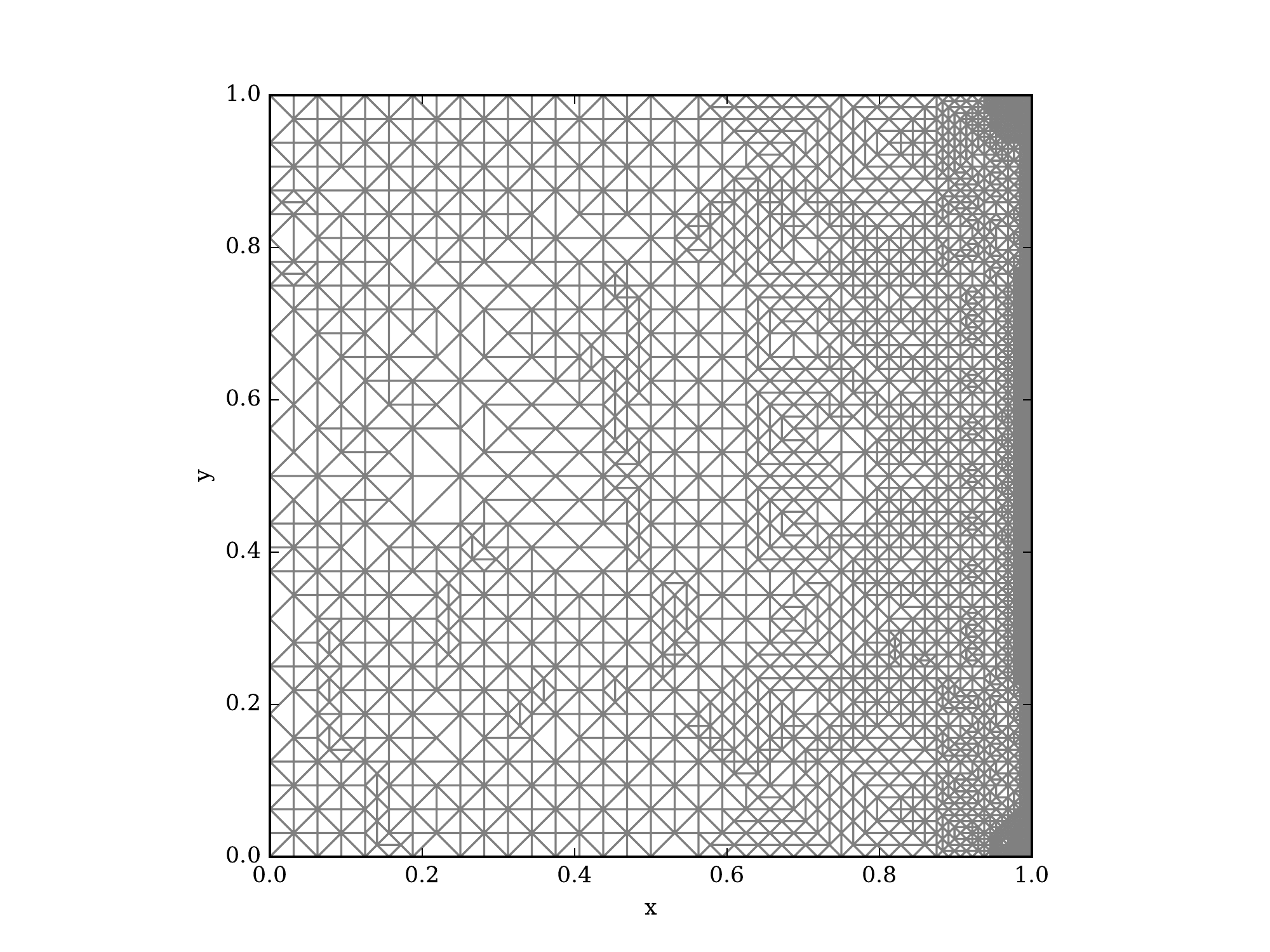}}\qquad
	\subfloat[ref. $\#$ 6, $\theta$ = 0.4]{\includegraphics[width=5cm, trim={3cm 1cm 3cm 1cm}, clip]{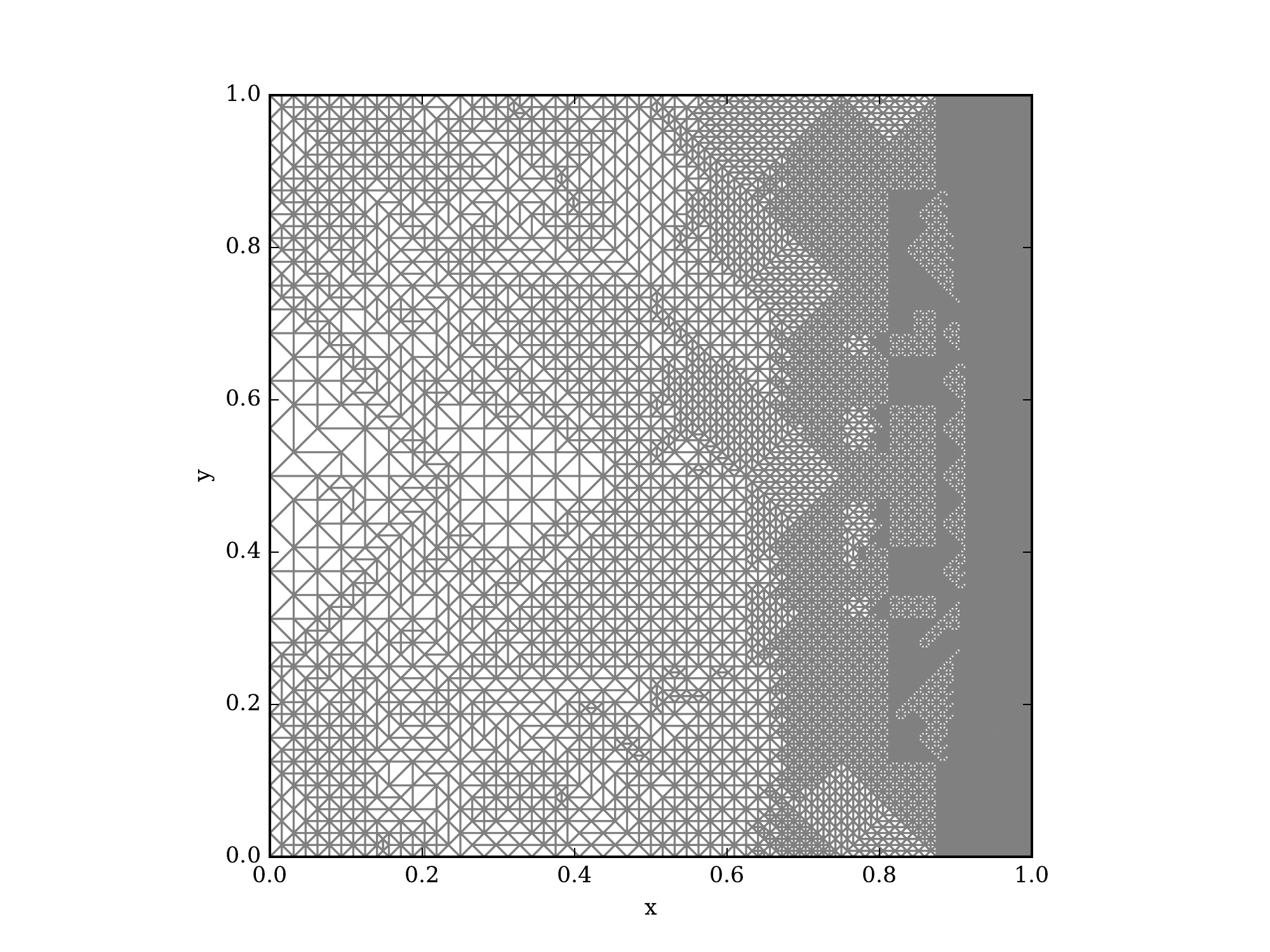}}\\
	\caption{Example \ref{ex:example-50-sicom-paper}. Evolution of the mesh on the refinement steps 2, 4, 6 with different 
	$\theta$ = 0.3 (left column) and $\theta$ = 0.4 (right column).}
	\label{fig:example-50-sicom-paper-mesh-bulk-30-40}
\end{figure}
\end{example}

According to the theory of asymptotic expansion \cite{Holmes2013}, the solution of the original 
problem depends continuously on $\varepsilon$. Using this fact, the authors of \cite{ChenSunXu2005, 
SunChenXu2010} 
use the so-called multilevel-homotopic-adaptive finite element method (MHA FEM) with respect to the 
diffusion constant $\varepsilon$, such that grids are iteratively adapted to a better approximate 
solution as the diffusion parameter gets decreased. 
The method gets initialised with $\tilde{\varepsilon} = \varepsilon_0$ (e.g., $\varepsilon_0 = 1$), 
and the mesh of the ad hoc diffusion problem 
$$- \tilde{\varepsilon} \, \Delta \, u + {\bf a} \cdot \nabla \, u + \lambda\, u = 0$$ 
gets refined adaptively. On the next iteration, the value of $\tilde{\varepsilon}$ is 
reduced, and the mesh adaptation is performed for the problem with the updated 
$\tilde{\varepsilon}$. Iterations continue until the desired value of $\varepsilon$ is reached, i.e., 
until condition $\tilde{\varepsilon} > \varepsilon$ is satisfied. 
For the reader's convenience, the steps described above are summarised 
in Algorithm \ref{alg:mha-fem}} (see also \cite{SunChenXu2010}). 
We use the MHA FEM 
in  
combination with the 
streamline 
diffusion FEM, where mesh adaptation is driven by the local distribution of the error indicator following 
from the majorant. Such an approach can also serve as a preprocessing for a  
problem so that 
one can generate problem-adapted meshes once {realistic} 
problems are concerned.

\begin{algorithm}[!t]
\caption{\quad Multilevel-homotopic-adaptive finite element method (MHA FEM)}
\label{alg:mha-fem}
\begin{algorithmic} 
\STATE {\bf Input:} $\mathcal{K}_h$ \COMMENT{initial mesh}
\STATE $\quad \,\qquad$ $\varepsilon_0$ \COMMENT{initial value of the diffusion parameter}
\STATE $\quad \,\qquad$ $\varepsilon$ \COMMENT{goal value of the diffusion parameter}
\STATE $\quad \qquad$ $N_{\rm ref}$ \COMMENT{number of refinement steps}
\vspace{4pt}
\STATE $\tilde{\varepsilon} = \varepsilon_0$\\[2pt]
\WHILE{$\tilde{\varepsilon} \geq \varepsilon$}
\vspace{4pt}
\FOR{$j = 1 : N_{\rm ref}$}
\vspace{4pt}
	\STATE Solve $- \tilde{\varepsilon} \, \Delta \, u + {\bf a} \cdot \nabla \, u + \lambda\, u = 0$, i.e., 
	reconstruct function $v_{\tilde{\varepsilon}}$ approximating $u$ \\[2pt]
	\STATE Compute the majorant 
	$\maj{} (v_{\tilde{\varepsilon}}, \flux; \beta)$ of the error $[ u - v_{\tilde{\varepsilon}} ]^2_{\maj{}}$ \\[2pt]
	\STATE Mark elements of $\mathcal{K}_h$ using the indicator following from $\maj{} (v_{\tilde{\varepsilon}}, \flux; \beta)$ \\[2pt]
	\STATE Refine mesh $\mathcal{K}_h$ using refinement procedure $\mathcal{R}$, i.e., 
	$\mathcal{K}_h = \mathcal{R} (\mathcal{K}_h)$ \\[2pt]
\ENDFOR
\STATE Update $\tilde{\varepsilon} = \tfrac{\tilde{\varepsilon}}{10}$\\[2pt]
\IF{$\tilde{\varepsilon} < \varepsilon$}
	\STATE $\tilde{\varepsilon} = \varepsilon$
\ENDIF
\ENDWHILE
\vspace{5pt}
\STATE {\bf Output:} $v_{\varepsilon}$ \COMMENT{approximate solution} \\
$\qquad \qquad$ $\mathcal{K}_{h_{\rm ref}}$ \COMMENT{refined mesh}
\end{algorithmic}
\end{algorithm}

\begin{example}
\label{ex:ex-1-chen-paper}
\rm
Assume now that $d = 2$, ${\bf a} = \big(1, 0 \big)^{\rm T}$, $\varepsilon = 1e\minus4$, 
$\lambda = 0$, and the RHS $f$ as well as and the boundary data $g = 0$ are chosen 
such that the solution of the problem is defined by  
$$u = \big(x^2 - e^{\rfrac{x - 1}{\varepsilon}}\big) \, y \, (1 - y)$$ 
(see Figure \ref{fig:ex-1-chen-paper-approximate-solution-ref} illustrating the approximate solution 
on the uniformly refined mesh). 
It exhibits an exponentially regular boundary (of the width $O(\varepsilon)$) near $\{x = 1\}$ 
perpendicular to the convection direction $\big(1, 0 \big)^{\rm T}$ (see also
\cite{SunChenXu2010}).

\begin{figure}[!ht]
	\centering
	\subfloat[$v \in P^1$, ref. 2: 1089 ND, 2048 EL]{
	\includegraphics[width=8cm, trim={2cm 1.5cm 2cm 2cm}]{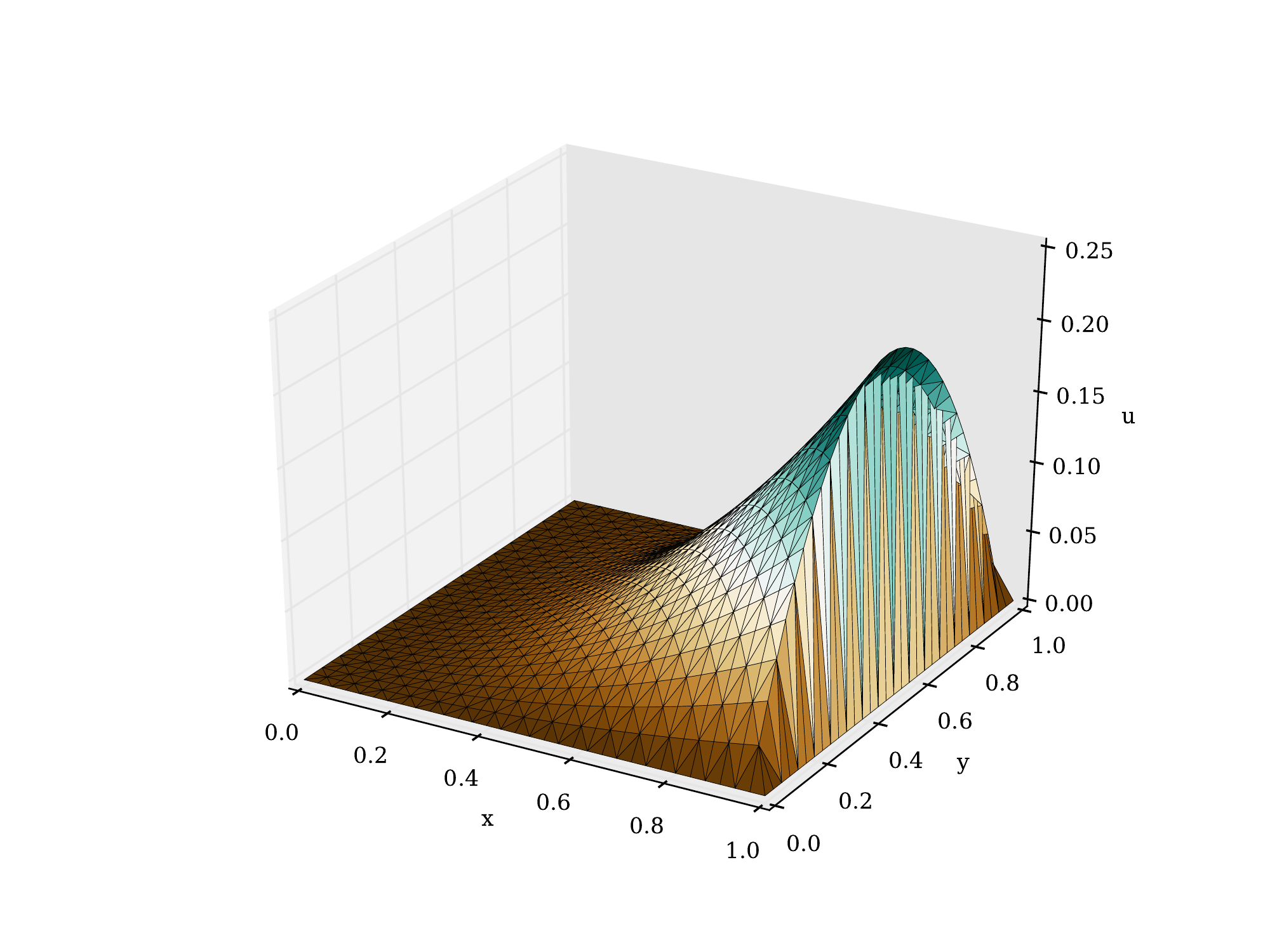}
	\label{fig:ex-1-chen-paper-approximate-solution-ref-2-b}}
	\subfloat[$v \in P^1$, ref. 3: 13353 ND, 26295 EL]{
	\includegraphics[width=8cm]{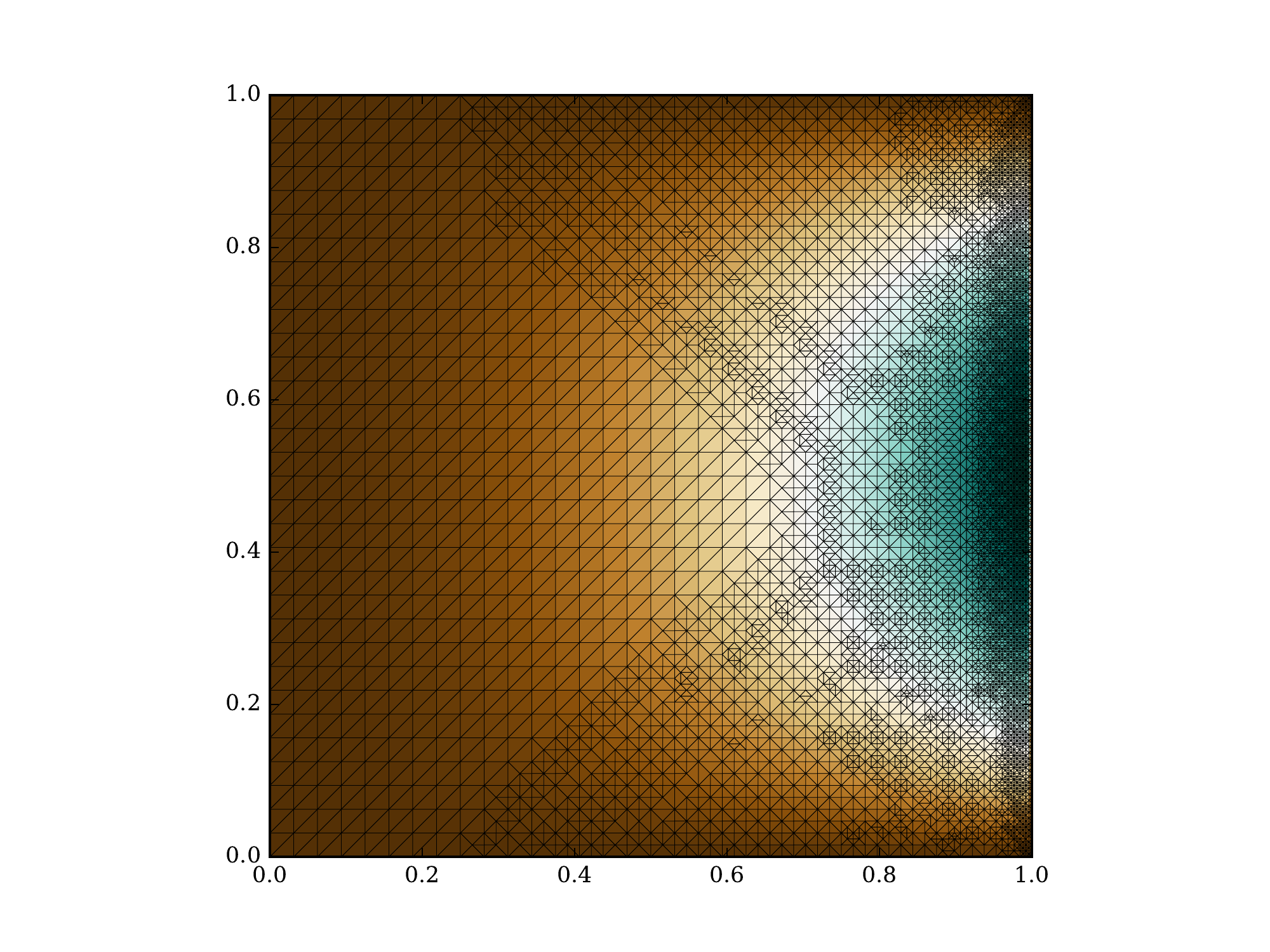}
	\label{fig:ex-1-chen-paper-approximate-solution-adapt-ref-2-a-stab}}
	\caption{Example \ref{ex:ex-1-chen-paper}. 
	(a) Approximate solution on the mesh with 1089 nodes and 2048 elements obtained by uniform refinement.
	(b) Mesh with 6833 nodes and 13404 elements with approximate solution obtained by adaptive refinement.}	
	\label{fig:ex-1-chen-paper-approximate-solution-ref}	
\end{figure}

We test two different scenarios of updating the parameter $\varepsilon$. The first one contains 4 
adaptive refinement steps for each of the levels $\varepsilon$, and the diffusion parameter 
decreases by $10$ on each level. Table \ref{tab:ex-1-chen-paper-sicom-paper-efficiency-indices-first-method}
contains the numerical results obtained by such an approach. Another way of `switching' from level to
level is refining until $h_{\rm min}$ becomes small enough (mesh becomes fine enough) to represent 
the current boundary layer $O (\varepsilon)$. The numerical testing of this approach is summarised in 
Table \ref{tab:ex-1-chen-paper-sicom-paper-efficiency-indices-second-method}. 
We see that with the first approach, we are able to reach the error ${[ e ]}_{\overline{\rm M}} = 1.49e\minus02$
with $36625$ d.o.f., whereas with the second approach we need $2951$ d.o.f. (about 12 times less) 
to reach a similar error, namely ${[ e ]}_{\overline{\rm M}} = 5.16e\minus02$, on the homotopic level 
$\varepsilon = 1e\minus3$. Therefore, choosing the number of adaptive refinements based on the 
approximation or mesh properties {seems to be the more} advantageous approach. 
Such a strategy can be used as an efficient mesh-generation technique to provide an initial setting for 
strongly convection-dominated problems (see Figure \ref{fig:example-53-convergence}). 
Here, each presented mesh corresponds to the final mesh on each of the $\varepsilon$ levels. 

\begin{table}[!t]
\centering
\footnotesize
\begin{tabular}{c|c|ccc|ccc}
ref. $\#$ & d.o.f. $\#$    & 
${[ e ]}_{\overline{\rm M}}$   & $\overline{\rm M}$    & $I_{\rm eff}(\overline{\rm M})$ & 
${[ e ]}_{\underline{\rm M}}$ & $\underline{\rm M}$  & $I_{\rm eff}(\underline{\rm M})$ \\
\midrule
\multicolumn{8}{c}{ $\varepsilon = 1$ } \\
\midrule
      1 &      120 &  3.19e-02 & 3.47e-02 &       1.09 &  3.20e-02 & 3.17e-02 &       0.99 \\
      2 &      183 &  2.48e-02 & 2.76e-02 &       1.11 &  2.49e-02 & 2.46e-02 &       0.99 \\
      3 &      269 &  1.87e-02 & 2.03e-02 &       1.08 &  1.88e-02 & 1.86e-02 &       0.99 \\
      4 &      403 &  1.62e-02 & 1.80e-02 &       1.11 &  1.63e-02 & 1.62e-02 &       0.99 \\
\midrule
\multicolumn{8}{c}{ $\varepsilon= 1e\minus1$ } \\
\midrule
      1 &      553 &  1.75e-02 & 2.33e-02 &       1.34 &  2.50e-02 & 1.95e-02 &       0.78 \\
      2 &      804 &  1.24e-02 & 1.65e-02 &       1.33 &  1.88e-02 & 1.41e-02 &       0.75 \\
      3 &     1214 &  1.01e-02 & 1.34e-02 &       1.32 &  1.52e-02 & 1.14e-02 &       0.75 \\
      4 &     1841 &  8.26e-03 & 1.08e-02 &       1.31 &  1.24e-02 & 9.26e-03 &       0.75 \\
\midrule
\multicolumn{8}{c}{ $\varepsilon= 1e\minus2$ } \\
\midrule      
      1 &     2556 &  2.67e-02 & 5.89e-02 &       2.21 &  4.01e-02 & 3.79e-02 &       0.94 \\
      2 &     3652 &  1.75e-02 & 4.14e-02 &       2.36 &  2.80e-02 & 2.61e-02 &       0.93 \\
      3 &     5494 &  1.33e-02 & 3.17e-02 &       2.39 &  2.14e-02 & 1.98e-02 &       0.93 \\
      4 &     8132 &  1.06e-02 & 2.47e-02 &       2.32 &  1.67e-02 & 1.55e-02 &       0.93 \\
\midrule
\multicolumn{8}{c}{ $\varepsilon= 1e\minus3$ } \\
\midrule            
      1 &    11311 &  3.60e-02 & 8.07e-02 &       2.24 &  5.69e-02 & 5.65e-02 &       0.99 \\
      2 &    16294 &  2.40e-02 & 6.24e-02 &       2.59 &  4.07e-02 & 4.03e-02 &       0.99 \\
      3 &    24662 &  1.98e-02 & 4.81e-02 &       2.43 &  3.15e-02 & 3.11e-02 &       0.99 \\
      4 &    36625 &  1.49e-02 & 3.96e-02 &       2.65 &  2.46e-02 & 2.42e-02 &       0.98 \\
\end{tabular}
\caption{Example \ref{ex:ex-1-chen-paper}. 
Majorant, minorant, and corresponding efficiency indices. 
Bulk 
marking with parameters $\theta = 0.1$, 4 refinements for each $\varepsilon$.}
\label{tab:ex-1-chen-paper-sicom-paper-efficiency-indices-first-method}
\end{table}

\begin{table}[!t]
\centering
\footnotesize
\begin{tabular}{c|c|ccc|ccc}
ref. $\#$ & d.o.f. $\#$    & 
${[ e ]}_{\overline{\rm M}}$   & $\overline{\rm M}$    & $I_{\rm eff}(\overline{\rm M})$ & 
${[ e ]}_{\underline{\rm M}}$ & $\underline{\rm M}$  & $I_{\rm eff}(\underline{\rm M})$ \\
\midrule
\multicolumn{8}{c}{ $\varepsilon = 1$ } \\
\midrule
      0 &       81 &  3.67e-02 & 3.71e-02 &       1.01 &  3.68e-02 & 3.65e-02 &       0.99 \\
\midrule
\multicolumn{8}{c}{ $\varepsilon= 1e\minus1$ } \\
\midrule
   0 &      120 &  5.00e-02 & 5.95e-02 &       1.19 &  7.32e-02 & 5.68e-02 &       0.78 \\
\midrule
\multicolumn{8}{c}{ $\varepsilon= 1e\minus2$ } \\
\midrule      
      0 &      170 &  1.01e-01 & 2.12e-01 &       2.11 &  2.34e-01 & 2.07e-01 &       0.88 \\
      1 &      250 &  8.22e-02 & 1.42e-01 &       1.72 &  1.52e-01 & 1.40e-01 &       0.92 \\
      2 &      362 &  5.97e-02 & 9.43e-02 &       1.58 &  1.00e-01 & 9.26e-02 &       0.92 \\
      3 &      549 &  4.13e-02 & 6.57e-02 &       1.59 &  6.99e-02 & 6.40e-02 &       0.92 \\
\midrule
\multicolumn{8}{c}{ $\varepsilon= 1e\minus3$ } \\
\midrule      
     0 &      816 &  9.96e-02 & 2.31e-01 &       2.32 &  2.29e-01 & 2.25e-01 &       0.98 \\
      1 &     1196 &  8.22e-02 & 1.67e-01 &       2.03 &  1.65e-01 & 1.61e-01 &       0.98 \\
      2 &     1914 &  6.52e-02 & 1.25e-01 &       1.91 &  1.23e-01 & 1.19e-01 &       0.97 \\
      3 &     2951 &  5.16e-02 & 9.71e-02 &       1.88 &  9.50e-02 & 9.16e-02 &       0.96 \\
\midrule
\multicolumn{8}{c}{ $\varepsilon= 1e\minus4$ } \\
\midrule      
     0 &     4450 &  1.05e-01 & 2.97e-01 &       2.84 &  2.90e-01 & 2.79e-01 &       0.96 \\
      1 &     6825 &  9.33e-02 & 2.31e-01 &       2.47 &  2.28e-01 & 2.13e-01 &       0.94 \\
      2 &    10540 &  8.19e-02 & 1.77e-01 &       2.16 &  1.76e-01 & 1.64e-01 &       0.93 \\
      3 &    15817 &  6.70e-02 & 1.36e-01 &       2.04 &  1.34e-01 & 1.25e-01 &       0.94 \\
      4 &    25318 &  5.41e-02 & 1.04e-01 &       1.92 &  1.02e-01 & 9.76e-02 &       0.95 \\
      5 &    37928 &  4.30e-02 & 8.18e-02 &       1.90 &  8.18e-02 & 7.75e-02 &       0.95 \\
\midrule
\multicolumn{8}{c}{ $\varepsilon= 1e\minus5$ } \\
\midrule      
      0 &    57437 & 1.01e-01 & 2.69e-01 &        2.66 & 2.64e-01 & 2.58e-01 &        0.97 \\
      1 &    86256 & 9.11e-02 & 2.06e-01 &        2.26 & 2.05e-01 & 1.95e-01 &        0.95 \\
      2 &   128442 & 7.47e-02 & 1.54e-01 &        2.06 & 1.56e-01 & 1.45e-01 &        0.93 \\ 
      3 &   197588 & 6.10e-02 & 1.17e-01 &        1.93 & 1.17e-01 & 1.12e-01 &        0.95 \\
\end{tabular}
\caption{Example \ref{ex:ex-1-chen-paper}. 
Majorant, minorant, and corresponding efficiency indices for different homotopic levels for 
bulk parameter $\theta = 0.1$.}
\label{tab:ex-1-chen-paper-sicom-paper-efficiency-indices-second-method}
\end{table}

\begin{figure}[!t]
\centering
\captionsetup[subfigure]{oneside, labelformat=empty}
\subfloat[]{
\includegraphics[width=6cm, trim={0cm 0cm 0cm 0cm}, clip]{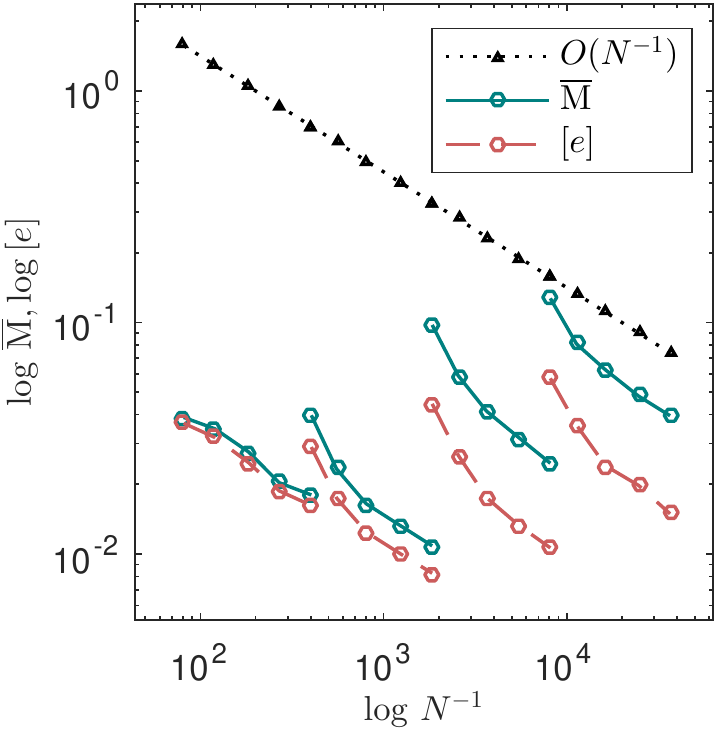}}
\qquad
\subfloat[]{
\includegraphics[width=6cm, trim={0cm 0cm 0cm 0cm}, clip]{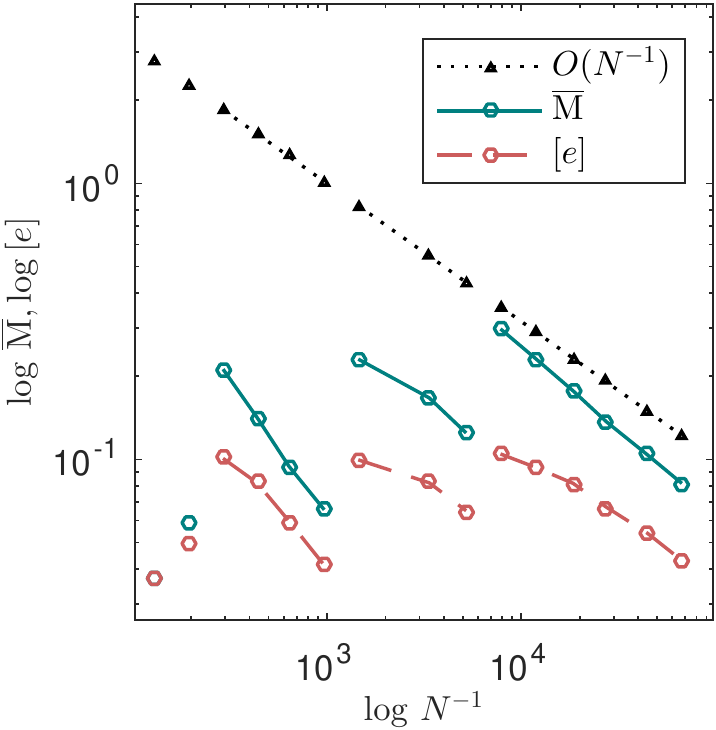}}
\caption{Example \ref{ex:ex-1-chen-paper}. 
Convergence of the error, majorant, and minorant corresponding to 
(a) Table \ref{tab:ex-1-chen-paper-sicom-paper-efficiency-indices-first-method} and 
(b) Table \ref{tab:ex-1-chen-paper-sicom-paper-efficiency-indices-second-method}.}
\label{fig:example-53-convergence}
\end{figure}

\begin{figure}[!t]
	\centering
	\captionsetup[subfigure]{oneside, labelformat=empty}
	\subfloat[$\varepsilon = 1$: \quad 81 d.o.f.]{
	\includegraphics[width=5cm, trim={3cm 1cm 3cm 1cm}, clip]{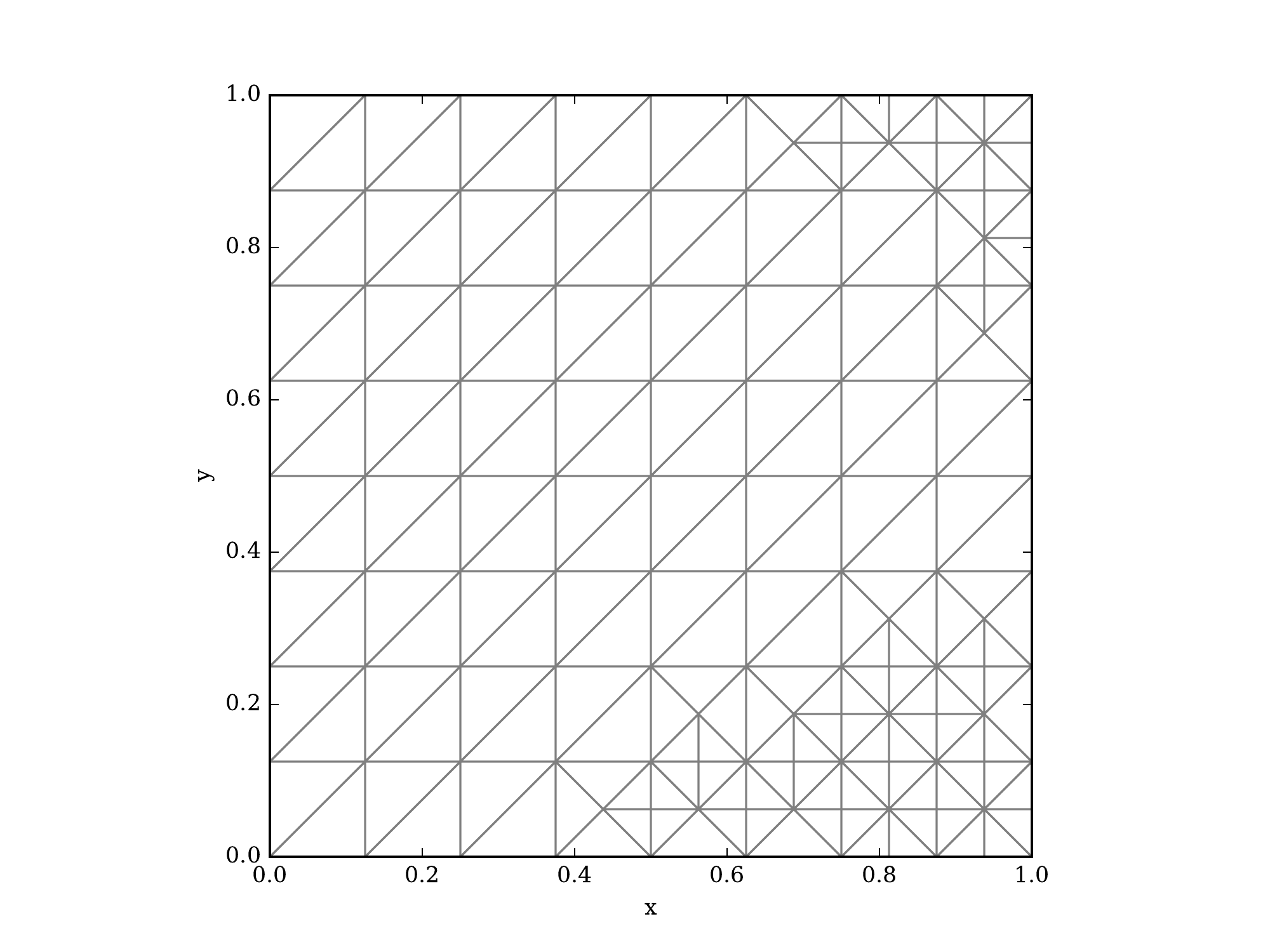}} \qquad
	\subfloat[$\varepsilon = 1e-1$: \quad 120 d.o.f.]{
	\includegraphics[width=5cm, trim={3cm 1cm 3cm 1cm}, clip]{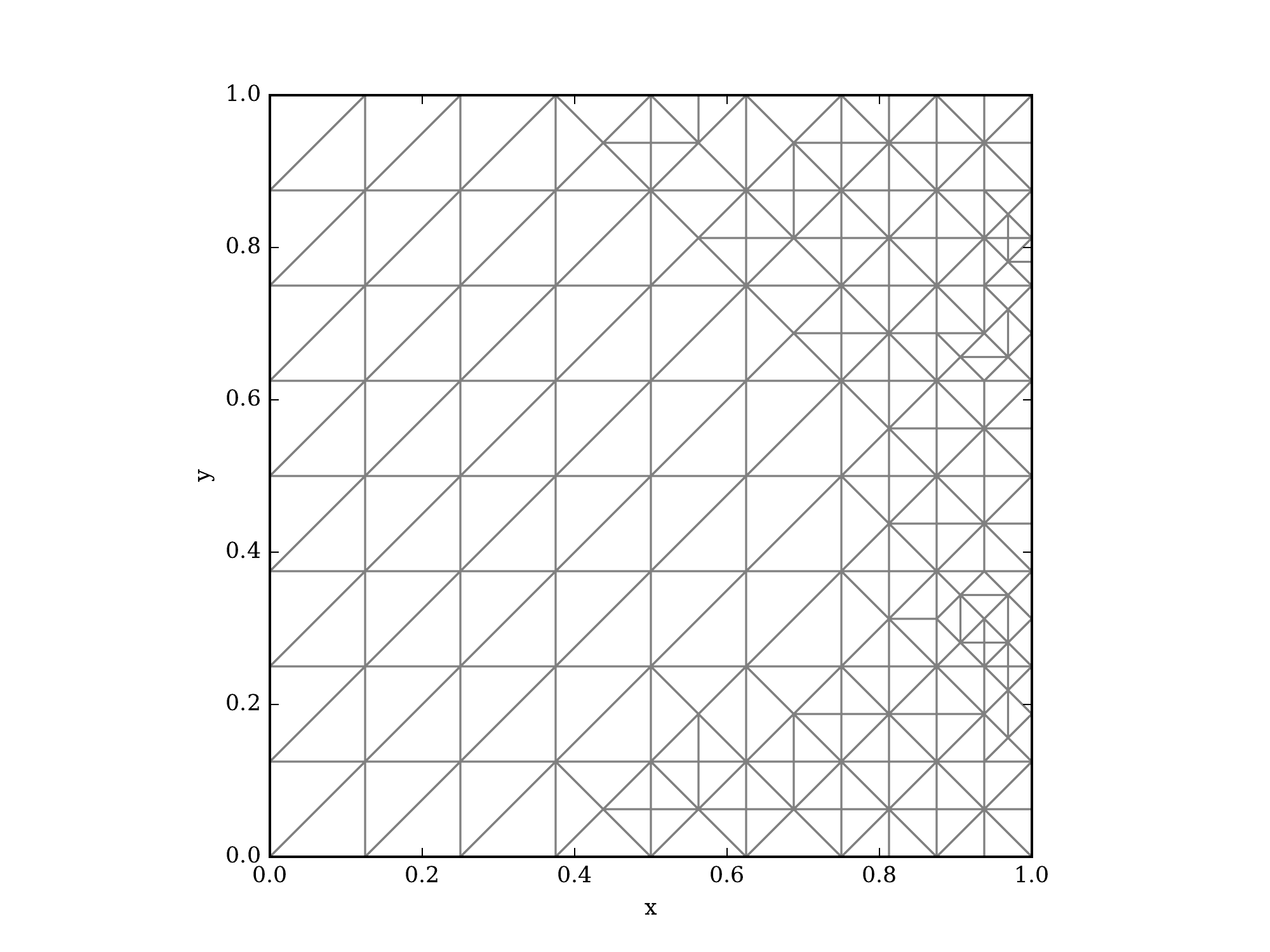}} \\
	\subfloat[r$\varepsilon = 1e-3$: \quad 2951 d.o.f.]{
	\includegraphics[width=5cm, trim={3cm 1cm 3cm 1cm}, clip]{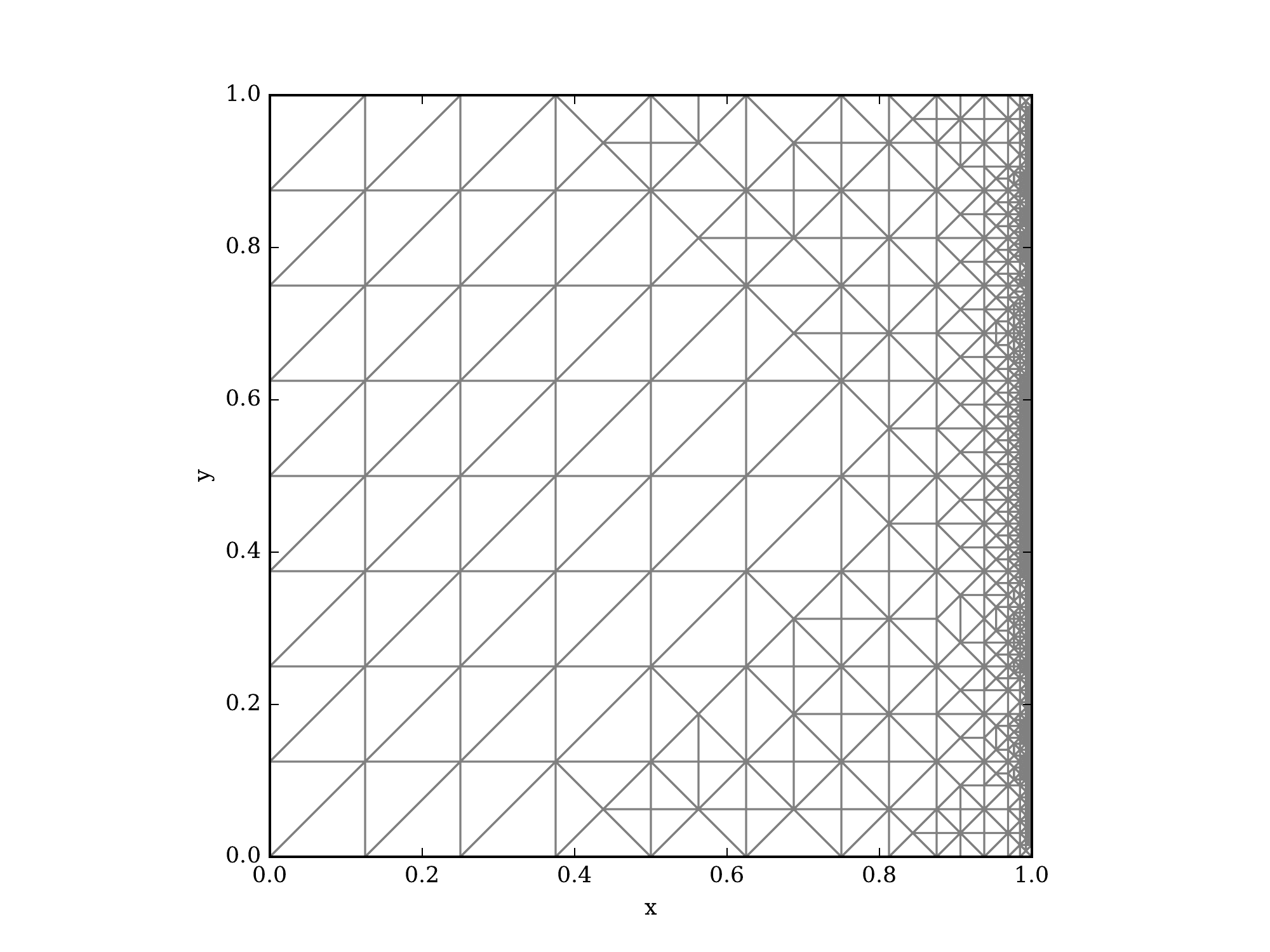}} \qquad
	\subfloat[$\varepsilon = 1e-4$: \quad 37928 d.o.f.]{
	\includegraphics[width=5cm, trim={3cm 1cm 3cm 1cm}, clip]{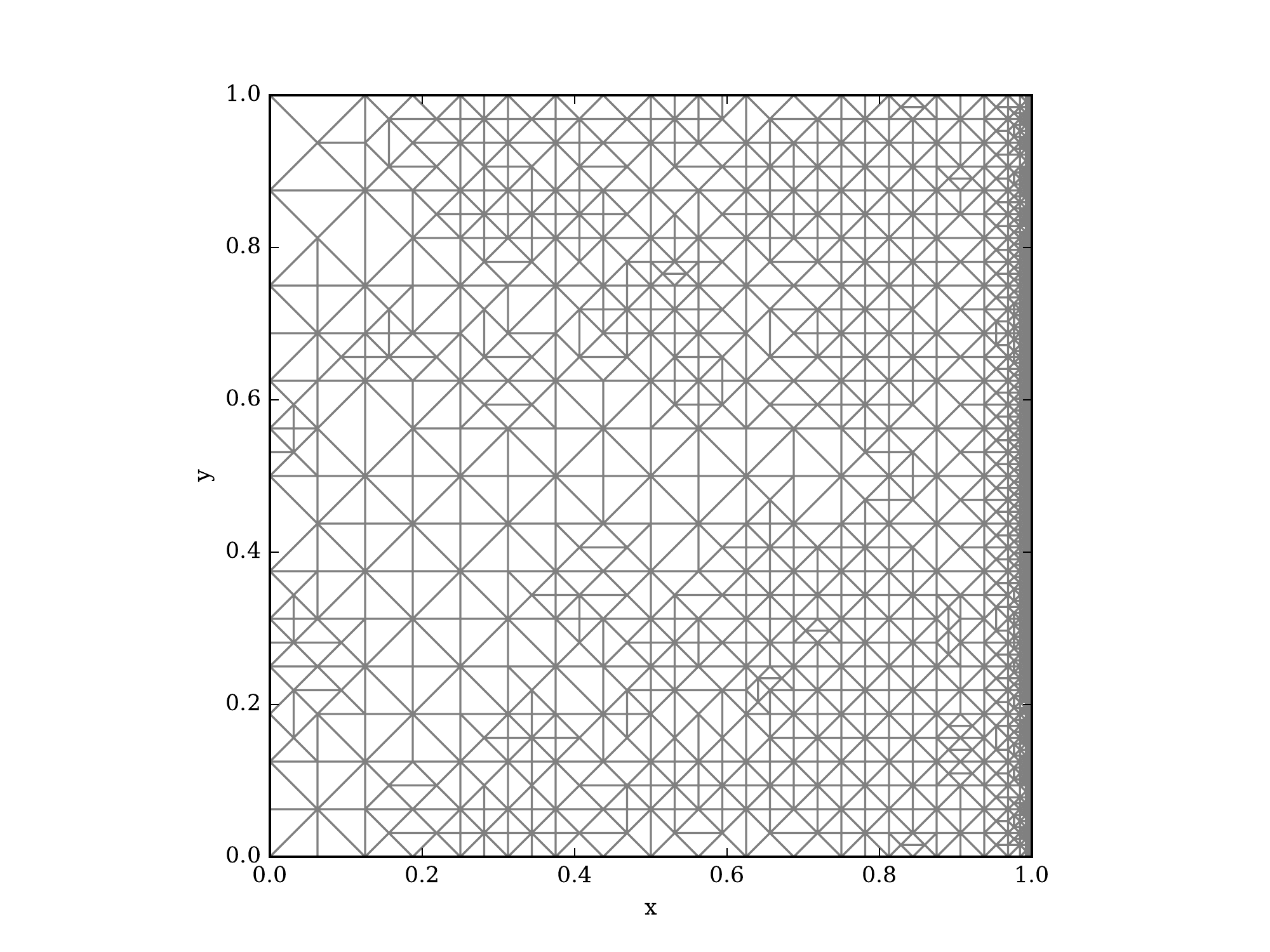}} \caption{Example \ref{ex:ex-1-chen-paper}. 
Mesh obtained on the different homotopy levels using error indicators following from
$\overline{\rm M}_{\mu_{\rm opt}}$ for $v \in \Pone$, $\flux \in \Ptwo$ with bulk marking $\theta = 0.1$.}
\label{fig:example-53-comparison-for-meshes-maj-error}
\end{figure}

\begin{table}[!ht]
\centering
\footnotesize
\begin{tabular}{c|cc|cc}
 $N_i$ & 
 $[ e ]_i$ & $r ([ e ]_i) = 
 \tfrac{\ln ([ e ]_{i} / [ e ]_{i + 1})}{\ln \Big( N_{i}^{-{1}/{d}} / N_{i + 1}^{-{1}/{d}} \Big)}$ & 
 $\| e \|_{L_{2}}$      & $r(\| e \|_{L_{2}})$ \\
 \midrule
\multicolumn{5}{c}{ $\varepsilon= 1e\minus2$ } \\
\midrule
      170 & 1.01e-01 &   1.06 & 4.46e-04 &   5.16 \\
      250 & 8.22e-02 &   1.72 & 1.65e-04 &   5.04 \\
      362 & 5.97e-02 &   1.77 & 6.48e-05 &   3.42 \\
\midrule
\multicolumn{5}{c}{ $\varepsilon= 1e\minus3$ } \\
\midrule
      816 & 9.96e-02 &   1.00 & 4.27e-05 &   3.82 \\
     1196 & 8.22e-02 &   0.99 & 2.06e-05 &   2.76 \\
     1914 & 6.52e-02 &   1.08 & 1.08e-05 &   2.44 \\
\midrule
\multicolumn{5}{c}{ $\varepsilon= 1e\minus4$ } \\
\midrule
     4450 & 1.05e-01 &   0.53 & 7.30e-06 &   2.46 \\
     6825 & 9.33e-02 &   0.60 & 4.31e-06 &   2.61 \\
    10540 & 8.19e-02 &   0.99 & 2.44e-06 &   2.97 \\
    15817 & 6.70e-02 &   0.91 & 1.34e-06 &   2.45 \\
    25318 & 5.41e-02 &   1.13 & 7.53e-07 &   2.19 \\
 \end{tabular}
\caption{Example \ref{ex:ex-1-chen-paper}. 
Convergence rates for the error $[ e ]$ and the $L_{2}$-error measure.}
\label{tab:eq:ex-3-chen-paper-convergence-rates-uniform-stab}
\end{table}

We also present convergence rates for the $L^{2}$-norm of the errors that has been addressed in 
\cite{SunChenXu2010}. The obtained numerical examples show that the $L^{2}$-error decreases 
as $O(h^2) \sim O(N^{-1})$ (see Table \ref{tab:eq:ex-3-chen-paper-convergence-rates-uniform-stab}).
\begin{figure}[!t]
	\centering
	\captionsetup[subfigure]{oneside, labelformat=empty}
	\subfloat[]{
	\includegraphics[width=5cm, trim={3cm 1cm 3cm 1cm}, clip]{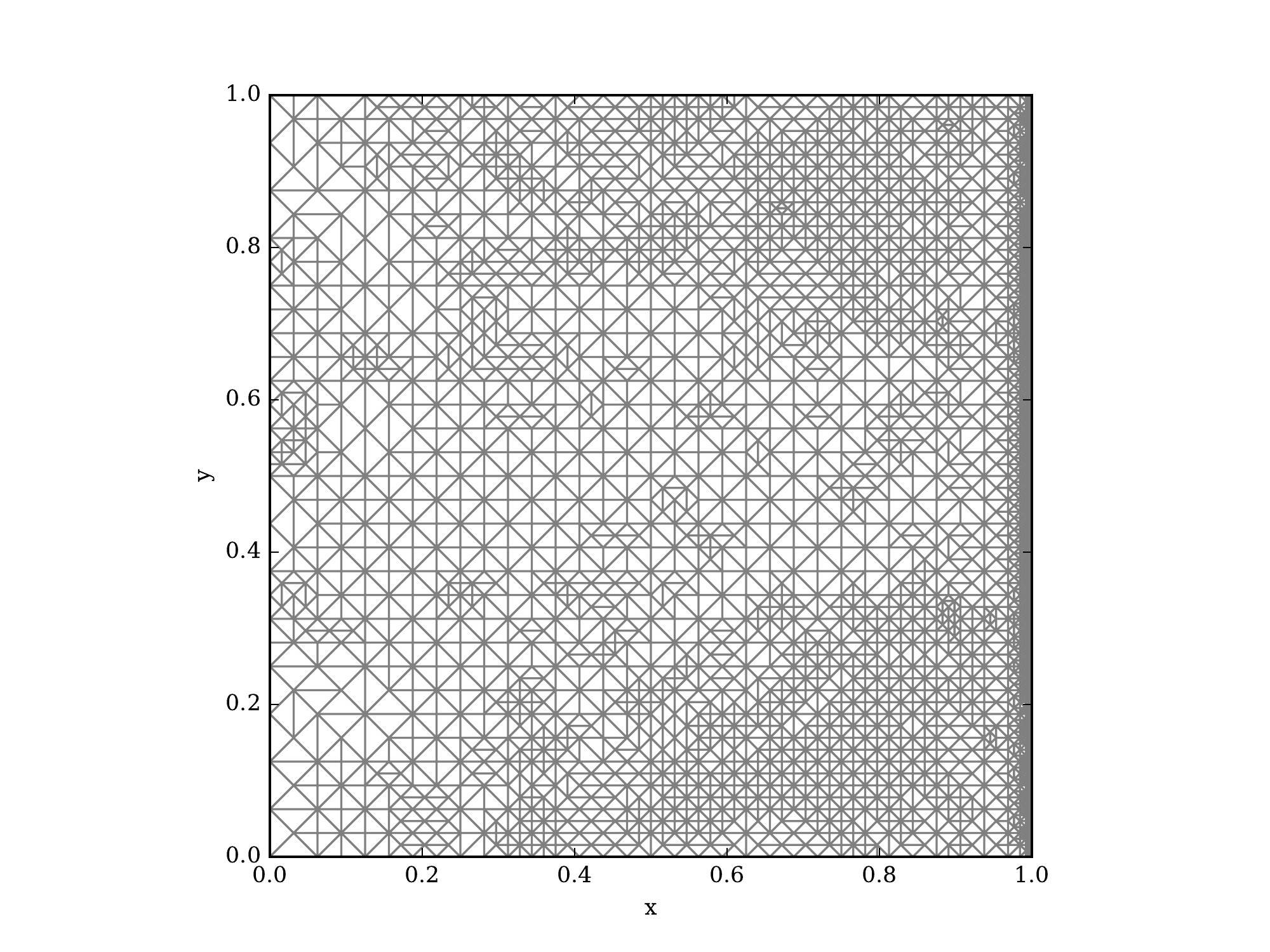}
	\label{fig:example-53-the-final-mesh}}
	\qquad
	\subfloat[]{
	\includegraphics[width=5cm, trim={3cm 1cm 3cm 1cm}, clip]{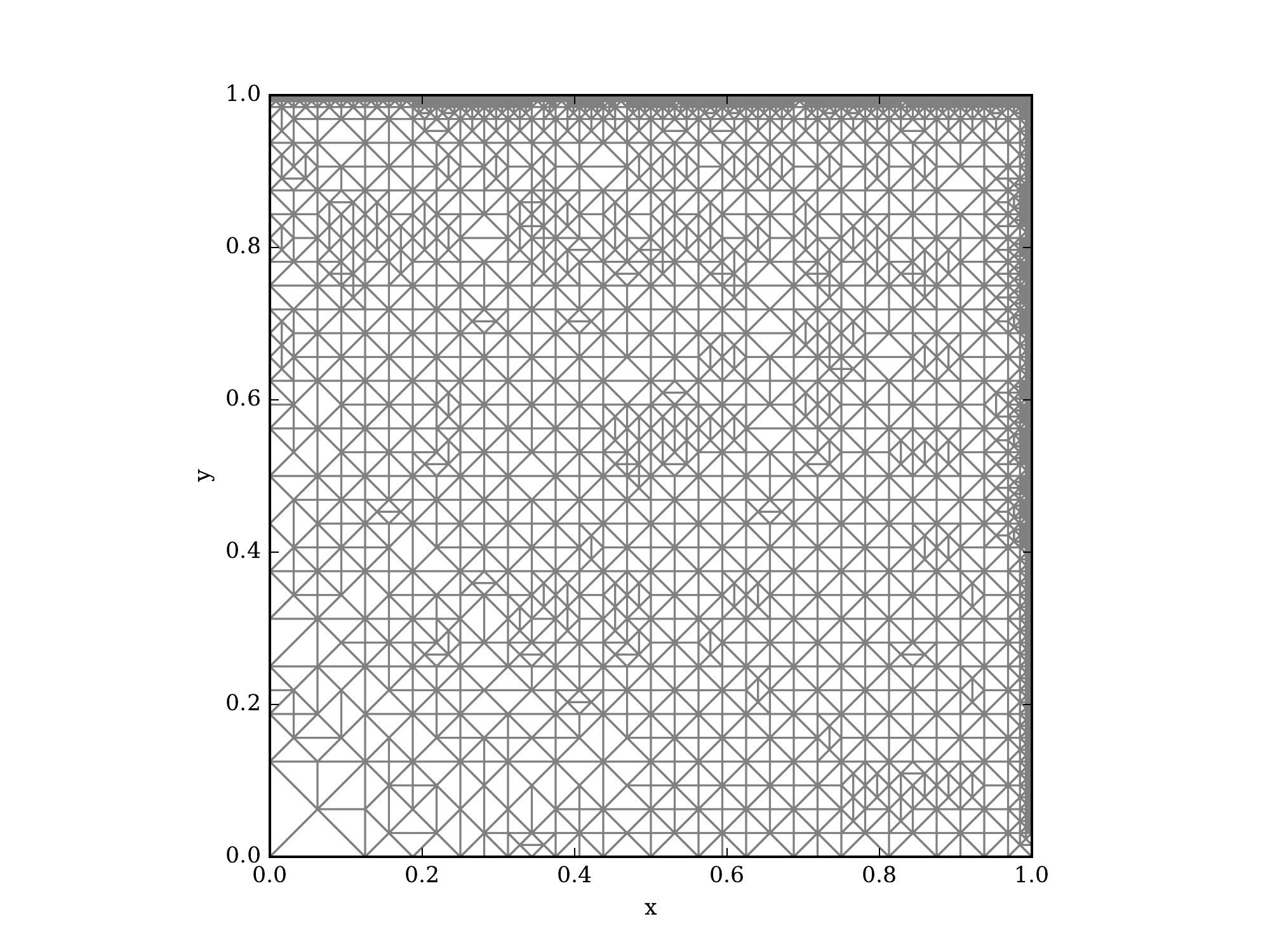}
	\label{fig:example-55-the-final-mesh}}
	\caption{The final, adaptive  
	mesh for (a) Example \ref{ex:ex-1-chen-paper} and (b) Example \ref{ex:ex-2-chen-paper}.}
\end{figure}

\end{example}

\begin{example}
\label{ex:ex-2-chen-paper}
\rm
Let us also consider problem
\eqref{eq:generalized-statement-static}
(on which 
the MHA FEM algorithm
was tested in \cite[Example 4]{SunChenXu2010}). We now choose 
$d = 2$, 
${\bf a} = \big(2, 3 \big)^{\rm T}$, 
$\varepsilon = 1e\minus4
$, 
$\lambda = 1$, 
and the exact solution is defined by the function 
$u = \big(x\,y^2 - y^2 \, e^{\rfrac{2(x - 1)}{\varepsilon}}
- x \, e^{\rfrac{3(y - 1)}{\varepsilon}} + e^{\rfrac{2(x - 1) + 3(y - 1)}{\varepsilon}}
\big)$ (it possesses regular boundary layers at $\{x = 1\}$ and 
$\{y = 1\}$).

First, we illustrate the behaviour of the majorant and minorant in Table 
\ref{tab:ex-2-chen-paper-sicom-paper-efficiency-indices-second-method}. Table 
\ref{tab:ex-2-chen-paper-convergence-rates-stab} provides the information on the 
$L_{2}$-error order of convergence. It is easy to observe that the obtained e.o.c. corresponds 
to the expected theoretical $O(h^{p+1})$, which is $O(h^2)$ for linear Lagrangian elements. 
The level of the error $\| e \|_{L_{2}} := 1.97e\minus05$ for 27817 d.o.f. is comparable to the 
$L^{2}$-error obtained in \cite{SunChenXu2010}, i.e., $\| e \|_{L_{2}} := 8.38e\minus05$ for 79104 d.o.f. 
($\varepsilon= 1e\minus4$). However, in our experiment this accuracy is obtained for $2.8$ times 
less d.o.f. then in \cite{SunChenXu2010}. 
We also present the mesh obtained by the adaptive 
refinement (driven by the indicator, which yields from the majorant), see Figure 
\ref{fig:example-55-the-final-mesh}. From the comparison of Figure 
\ref{fig:example-55-the-final-mesh} to Figure 9 in \cite[Example 5]{SunChenXu2010}, it is
easy to see that in our case the refinement is rather concentrated on the region with a boundary layer {rather than the region, where solution is smooth.}

\begin{table}[!t]
\centering
\footnotesize
\begin{tabular}{c|c|ccc|ccc}
ref. $\#$ & d.o.f. $\#$    & 
${[ e ]}_{\overline{\rm M}}$   & $\overline{\rm M}$    & $I_{\rm eff}(\overline{\rm M})$ & 
${[ e ]}_{\underline{\rm M}}$ & $\underline{\rm M}$  & $I_{\rm eff}(\underline{\rm M})$ \\
\midrule
\multicolumn{8}{c}{ $\varepsilon = 1$ } \\
\midrule
       0 &       81 &  2.04e-02 & 2.28e-02 &       1.12 &  2.31e-02 & 2.09e-02 &       0.91 \\
\midrule
\multicolumn{8}{c}{ $\varepsilon= 1e\minus1$ } \\
\midrule
     0 &      119 &  4.45e-01 & 6.65e-01 &       1.49 &  8.98e-01 & 6.09e-01 &       0.68 \\
\midrule
\multicolumn{8}{c}{ $\varepsilon= 1e\minus2$ } \\
\midrule      
       0 &      169 &  6.68e-01 & 2.07e+00 &       3.10 &  2.81e+00 & 1.75e+00 &       0.62 \\
      1 &      241 &  6.53e-01 & 1.54e+00 &       2.36 &  2.08e+00 & 1.36e+00 &       0.65 \\
      2 &      360 &  5.74e-01 & 1.16e+00 &       2.03 &  1.54e+00 & 1.05e+00 &       0.68 \\
      3 &      559 &  4.96e-01 & 8.91e-01 &       1.80 &  1.14e+00 & 8.20e-01 &       0.72 \\
\midrule
\multicolumn{8}{c}{ $\varepsilon= 1e\minus3$ } \\
\midrule      
     0 &      862 &  6.40e-01 & 2.71e+00 &       4.24 &  3.56e+00 & 2.43e+00 &       0.68 \\
      1 &     1320 &  6.66e-01 & 2.07e+00 &       3.12 &  2.73e+00 & 1.90e+00 &       0.70 \\
      2 &     2003 &  6.48e-01 & 1.62e+00 &       2.51 &  2.12e+00 & 1.52e+00 &       0.72 \\
      3 &     3134 &  5.79e-01 & 1.23e+00 &       2.13 &  1.59e+00 & 1.17e+00 &       0.74 \\
\midrule
\multicolumn{8}{c}{ $\varepsilon= 1e\minus4$ } \\
\midrule      
     0 &     4883 &  5.56e-01 & 3.59e+00 &       6.45 &  4.61e+00 & 3.34e+00 &       0.72 \\
      1 &     7445 &  6.36e-01 & 2.82e+00 &       4.44 &  3.68e+00 & 2.67e+00 &       0.73 \\
      2 &    11693 &  6.66e-01 & 2.19e+00 &       3.28 &  2.87e+00 & 2.09e+00 &       0.73 \\
      3 &    17940 &  6.53e-01 & 1.71e+00 &       2.61 &  2.25e+00 & 1.65e+00 &       0.74 \\
      4 &    27817 &  6.00e-01 & 1.32e+00 &       2.19 &  1.71e+00 & 1.29e+00 &       0.75 \\
      5 &    43735 &  5.21e-01 & 1.02e+00 &       1.96 &  1.29e+00 & 1.00e+00 &       0.77 \\
\midrule
\multicolumn{8}{c}{ $\varepsilon= 1e\minus5$ } \\
\midrule      
\end{tabular}
\caption{Example \ref{ex:ex-2-chen-paper}. 
Majorant, minorant, and corresponding efficiency indices for different homotopic levels with  
bulk parameter $\theta = 0.1$.}
\label{tab:ex-2-chen-paper-sicom-paper-efficiency-indices-second-method}
\end{table}

\begin{table}[!ht]
\centering
\footnotesize
\begin{tabular}{c|cc}
 $N_i$ & $\| e \|_{L_{2}}$      & 
 $r(\| e \|_{L_{2}}) := \tfrac{\ln (\| e \|_{i} / \| e \|_{i + 1})}{\ln \Big( N_{i}^{-\rfrac{1}{d}} / N_{i + 1}^{-\rfrac{1}{d}} \Big)}$ \\
\midrule
\multicolumn{3}{c}{ $\varepsilon= 1e\minus2$ } \\
\midrule
      169 & 5.72e-03 &   3.63 \\
      241 & 3.01e-03 &   3.26 \\
      360 & 1.56e-03 &   3.01 \\
\midrule
\multicolumn{3}{c}{ $\varepsilon= 1e\minus3$ } \\
\midrule
      862 & 9.44e-04 &   2.63\\
     1320 & 5.39e-04 &   2.60 \\
     2003 & 3.14e-04 &   2.77 \\
\midrule
\multicolumn{3}{c}{ $\varepsilon= 1e\minus4$ } \\
\midrule
   4883 & 1.61e-04 &   2.20 \\
     7445 & 1.01e-04 &   2.31 \\
    11693 & 6.00e-05 &   2.45 \\
    17940 & 3.55e-05 &   2.68 \\
    27817 & 1.97e-05 &   2.67 \\
 \end{tabular}
\caption{Example \ref{ex:ex-2-chen-paper}. 
Convergence rates of the $L_{2}$-norm.}
\label{tab:ex-2-chen-paper-convergence-rates-stab}
\end{table}

\end{example}

\subsection{Time-dependent problems}

In the next subsection, we consider the time-dependent reaction-convection-diffusion problem 
that follows from \eqref{eq:fokker-planck-equation}--\eqref{eq:parabolic-robin-bc}. We assume, 
that the balance equation \eqref{eq:fokker-planck-equation} is written in the form 
{
\begin{equation}
u_t - \varepsilon \, \Delta \, u + {\bf a} \cdot \nabla \, u + \lambda\, u = 0,
\label{eq:time-dependent-problem}
\end{equation}
}
such that the reaction and convection are assumed to be independent of each other.  
Since the incremental time-stepping method possesses several disadvantages, for instance, 
time-consuming sequentiality which leads to complications in the parallelisation,  
we use a variational space-time approach to solve 
\eqref{eq:fokker-planck-equation}--\eqref{eq:parabolic-robin-bc}. 

For our numerical examples, the generalised formulation \eqref{eq:mix:state-generalized-statement} is 
discussed in the Petrov--Galerkin setting with different test and ansatz spaces (following the spirit 
of \cite{Steinbach2015}): find 
$$u \in V_{u_0} := \Big\{ \, v \in L_{2}(0, T; H^1_0 (\Omega)) \cap H^1(0, T; H^{-1}(\Omega)) \;| \;
v(x, 0) = u_0(x) \; \mbox{for} \;x \in \Omega \, \Big\},$$ 
such that 
\begin{equation}
	\varepsilon \, (\nabla {u}, \nabla {\eta})_{Q_T} 
	+ (\lambda \, u, \eta)_{Q_T} 
	+ ( {\bf a} \cdot \nabla \, u, \eta)_{Q_T} 
	+ (u_t, \eta)_{Q_T} =: a(u, \eta) = \ell (\eta) 
	:= (f, \eta)_{Q_T}, 
	\label{eq:diffusion-generalized-statement}
\end{equation}
for all $\eta \in W := L_{2}(0, T; H^1_0 (\Omega))$.
The stability condition is provided by the extension of \cite[Theorem 2.1]{Steinbach2015}, 
i.e., having that $W \subset  V_{u_0}$, the inequality 
$$c \, \| u \|_ {V_{u_0}} \leq {\rm sup}_{v \in W, v \neq 0} \tfrac{a(u, v)}{\| v \|_{W}}$$
holds with a positive constant $c$. After `homogenisation' of the problem by splitting 
$u = \bar{u}(x, t) + \bar{u}_0, (x, t) \in Q_T$, where $\bar{u}_0$ is an extension of 
$u_0 \in H^{1}_0(\Omega)$ to the rest of the space-time cylinder. Hence, we aim 
to find $\bar{u} \in  V_{0}$ such that
\begin{equation}
a(\bar{u}, v) = \textlangle f, v_h \textrangle_{Q_T} - a(\bar{u}_0, v), \quad \forall v \in W.
\label{eq:petrov-galerkin-formulation}
\end{equation}
Using the fact that $V_{0} \subset W$, we can state existence and uniqueness of the solution of the
weak formulation \eqref{eq:petrov-galerkin-formulation} based on 
\cite[Theorem 3.7]{Steinbach2008}.

Let the discretisation spaces be $V_{h0} \subset V_0 $ and $W_{h} \subset W$ such that 
the inclusion $V_{h0} \subset W_{h}$ is satisfied. Then, the discrete version of 
\eqref{eq:petrov-galerkin-formulation} can be presented as follows: find $\bar{u}_h \in V_{h0}$
satisfying
\begin{equation}
a(\bar{u}_h, v_h) = \textlangle f, v_h \textrangle_{Q_T} - a(\bar{u}_0, v_h), \quad \forall v_h \in W.
\label{eq:petrov-galerkin-formulation-discrete}
\end{equation}
The unique solvability of \eqref{eq:petrov-galerkin-formulation-discrete} follows from
the discrete stability condition 
$$c \, \| u_h \|_ {V_{h0}} 
\leq {\rm sup}_{v_h \in W_h, v_h \neq 0} 
       \tfrac{a(u_h, v_h)}{\| v_h \|_{W}}, 
       \quad \forall u_h \in V_{h0}, \quad c > 0,$$
which is proved in \cite[Theorem 3.1]{Steinbach2015}. 

Finally, we introduce the FE spaces for $u_h \in V_{h0}$ and $w_h \in W_{h}$, i.e., 
$$V_{h0} = \Pone (Q_h) \cap V_0, \quad W_{h} = \Pone (Q_h) \cap W,$$ 
where $\Pone (Q_h)$ is the FE space with a basis provided by piecewise linear 
and continuous functions (e.g., Lagrange polynomials). This brings us to the 
final approximation result for $\bar{u} \in V_0$ and $\bar{u}_h \in V_{h0}$.  
{\color{blue}
According to \cite[Theorem 3.3]{Steinbach2015}, by assuming that 
$u \in H^2(Q)$, we obtain the following energy 
error estimate
\begin{equation}
\| \bar{u} - \bar{u}_h \|_W \leq \bar{c} \, h \, |u|_{H^2(Q)}
\label{eq:energy-estimate-h2}
\end{equation}
with a $\bar{c} >0$ independent on $h$. 
More general estimate follows as a corollary \cite[Corollary 3.4]{Steinbach2015}. If
$\bar{u} \in H^s_0(Q_T)$, $s \in [1, 2]$, then the estimate 
\begin{equation}
\| \bar{u} - \bar{u}_h \|_W \leq \bar{C} \, h^{s-1} \,| \bar{u}|_{H^s_0(Q_T)}, \quad 
s \in [1, p+1],
\label{eq:energy-estimate-hs}
\end{equation}
follows. Here, $\bar{C}$ is a generic positive constant, independent from the mesh parameter $h$.
This theoretical result is confirmed by results of numerical tests (see Table 
\ref{tab:example-39-40-41-comparison-for-mu-opt} and 
\ref{tab:example-39-40-41-mu-opt-error-id}).}

Since the majorant is defined as an integral over the whole cylinder $Q_T$, it comes rather natural 
to use it for the error control. 
It also does not dependent on the discretisation used, therefore 
can be applied to the solutions reconstructed with fully unstructured grids.  
Next, 
we demonstrate numerical results obtained by applying functional type error 
estimates in case of a space-time discretisation. 

\begin{example}
\label{ex:example-39-40-41}
\rm
Consider a problem defined on the space-time cylinder $Q_T := \Omega \times (0, T)$ such that 
$\Omega = (0, 1)$, $T = 2$,  $\partial \Omega = \Gamma_N$. For the diffusion-reaction-convection 
operator, see \eqref{eq:time-dependent-problem},
we select $\varepsilon = 1$, ${\bf a = 1}$, and three different 
cases for $\lambda (x, t)$ (see Figure 
\ref{fig:lambda}) given as follows 
\begin{alignat}{2}
(a) \quad & \, \lambda (x, t) = \tfrac{1}{0.05 \, \sqrt{2 \pi}} e^{- 200 \, \left(x - 0.2 \right)^2}, \\
(b) \quad & \, \lambda (x, t) =  0.001 \, (x + 0.001) \, (t \sin t + 1), \\
(c) \quad & \, \lambda (x, t) =  1000 \, (x + 0.001) \, (t \sin t + 1).
\end{alignat}
The initial condition as well as the RHS 
are set to 
$$u_0 = \sin 2\, \pi \, +\cos \pi x \quad 
\mbox{and} 
\quad 
f = u_0 \, (\cos (2 t) - 2 t \sin(2 t)) 
+ (\pi^2 \, \cos (\pi x) - \pi \, \sin \pi x + \lambda \,  u_0) \, \left(t \, \cos 2 \, t  + 1\right).$$
The corresponding solution is illustrated in Figure \ref{fig:example-39-solution}. 

\begin{figure}[!t]
	\centering
	\subfloat[]{
	\includegraphics[width=6cm]{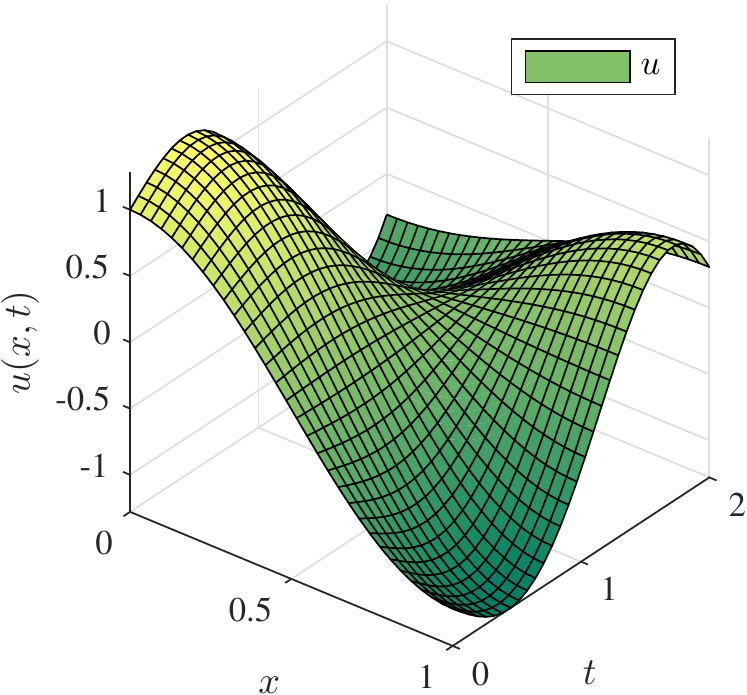}}
	\caption{Example \ref{ex:example-39-40-41}. The exact solution $u = \left( \sin 2\, \pi \, +\cos (\pi x) \right)\, \left(t \, \cos 2 \, t  + 1\right)$.}
	\label{fig:example-39-solution}
\end{figure}

\begin{figure}[!t]
	\centering
	\subfloat[]{
	\includegraphics[width=5cm]{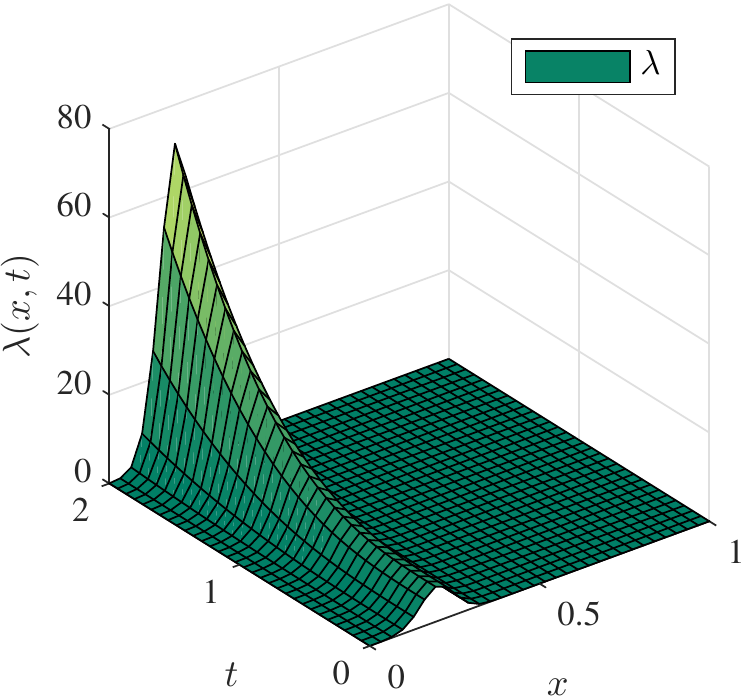}
	\label{fig:example-39}}
	\quad
	\subfloat[]{
	\includegraphics[width=5cm]{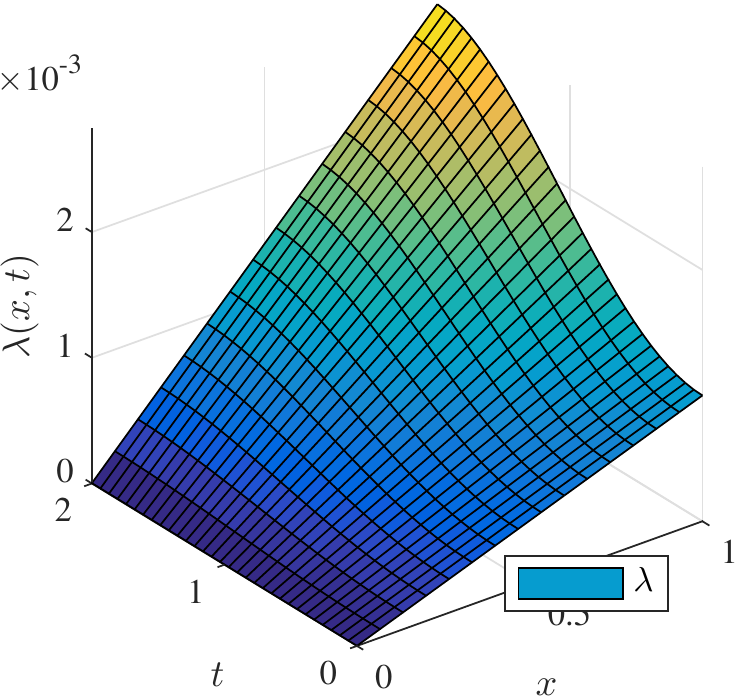}
	\label{fig:example-40-b}}
	\quad
	\subfloat[]{
	\includegraphics[width=5cm]{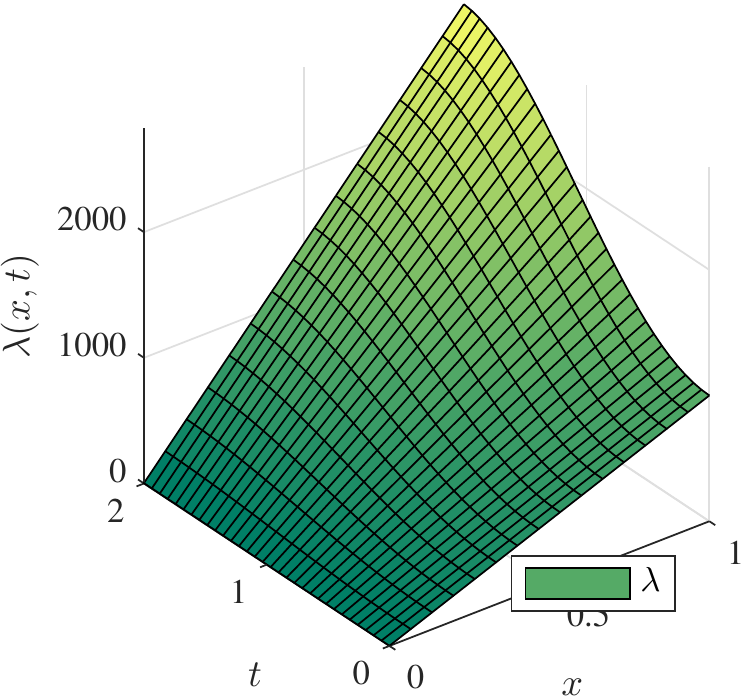}
	\label{fig:example-41-c}}
	\caption{Example \ref{ex:example-39-40-41}. Reaction functions $\lambda(x, t)$
	for the case (a), (b), and (c).}
	\label{fig:lambda}
\end{figure}

\begin{table}[!ht]
\centering
\footnotesize
\newcolumntype{g}{>{\columncolor{gainsboro}}c} 	
\newcolumntype{k}{>{\columncolor{lightgray}}c} 	
\newcolumntype{s}{>{\columncolor{silver}}c} 
\newcolumntype{a}{>{\columncolor{ashgrey}}c}
\newcolumntype{b}{>{\columncolor{battleshipgrey}}c}
\begin{tabular}{c|ccg|ccc}
\midrule
ref.   & 
d.o.f. $\#$& $[e]$ & $\Ieff (\overline{\rm M}_{\mu_{\rm opt}})$ & 
d.o.f. $\#$& $[e]$ & $\Ieff (\overline{\rm M}_{\mu = 1})$ \\
\midrule
  & \multicolumn{6}{c}{(a) 
$\lambda (x, t) = \tfrac{1}{0.05 \, \sqrt{2 \pi}} \, 
e^{- 200 \, \left(x - 0.2 \right)^2}$} \\
\midrule
      1 &      271 &  2.54e-01 & 1.45 & 258 &  3.01e-01 & 1.04 $\cdot$ 1e20 \\
      3 &      887 &  1.35e-01 & 1.38  & 787 &  2.22e-01 & 3.14 $\cdot$ 1e24 \\
      5 &     2864 &  7.40e-02 & 1.37 & 2542 &  1.30e-01 & 2.11 $\cdot$ 1e26 \\
      7 &     8692 &  4.10e-02 & 1.39 & 8159 &  7.32e-02 & 6.34 $\cdot$ 1e26 \\
      9 &    25825 &  2.34e-02 & 1.42 & 25425 &  3.90e-02 & 1.03 $\cdot$ 1e26 \\
\midrule
  & \multicolumn{6}{c}{(b)
  $\lambda (x, t) =  0.001 \, (x + 0.001) \, (t \sin t + 1)$} \\
\midrule
      1 &      272 &  2.53e-01 & 1.49 &  257 &  3.00e-01 & 39.65\\
      3 &      869 &  1.39e-01 & 1.46 &  772 &  2.36e-01 & 42.15\\
      5 &     2768 &  7.54e-02 & 1.45 &  2411 &  1.64e-01 &  56.30 \\
      7 &     8450 &  4.15e-02 & 1.48 &  7467 &  6.17e-02 & 64.74\\
      9 &    25240 &  2.38e-02 & 1.57 & 22598 &  3.70e-02 &  88.56  \\
\midrule
 & \multicolumn{6}{c}{(c) 
 $\lambda (x, t) =  1000 \, (x + 0.001) \, (t \sin t + 1)$} \\
\midrule
      1 &   272 &  2.53e-01 & 1.49 &      963 &  2.06e-02 &  82.87 \\
      3 &    869 &  1.39e-01 & 1.46 &  2841 &  1.28e-02 &  136.27 \\
      5 &    2768 &  7.54e-02 & 1.45 & 9111 &  8.59e-03 & 213.33 \\
      7 &    8450 &  4.15e-02 & 1.48 & 27594 &  5.39e-03 & 329.32 \\
      9 &    25240 &  2.38e-02 & 1.57 &82646 &  1.18e-03 & 1829.54 \\
\end{tabular}
\caption{Example \ref{ex:example-39-40-41}. 
Comparison of $\overline{\rm M}_{\mu_{\rm opt}}$ and $\overline{\rm M}_{\mu = 1}$
$v \in \Pone$, $\flux \in \Ptwo$ with bulk marking $\theta = 0.6$.}
\label{tab:example-39-40-41-comparison-for-mu-opt}
\end{table}

\begin{figure}[!t]
	\centering
	\captionsetup[subfigure]{oneside, labelformat=empty}
	\subfloat[ref. based on $\overline{\rm M}_{\mu_{\rm opt}}$ for (a)]{
	\includegraphics[width=5.5cm, trim={6cm 0.5cm 6cm 1cm}, clip]{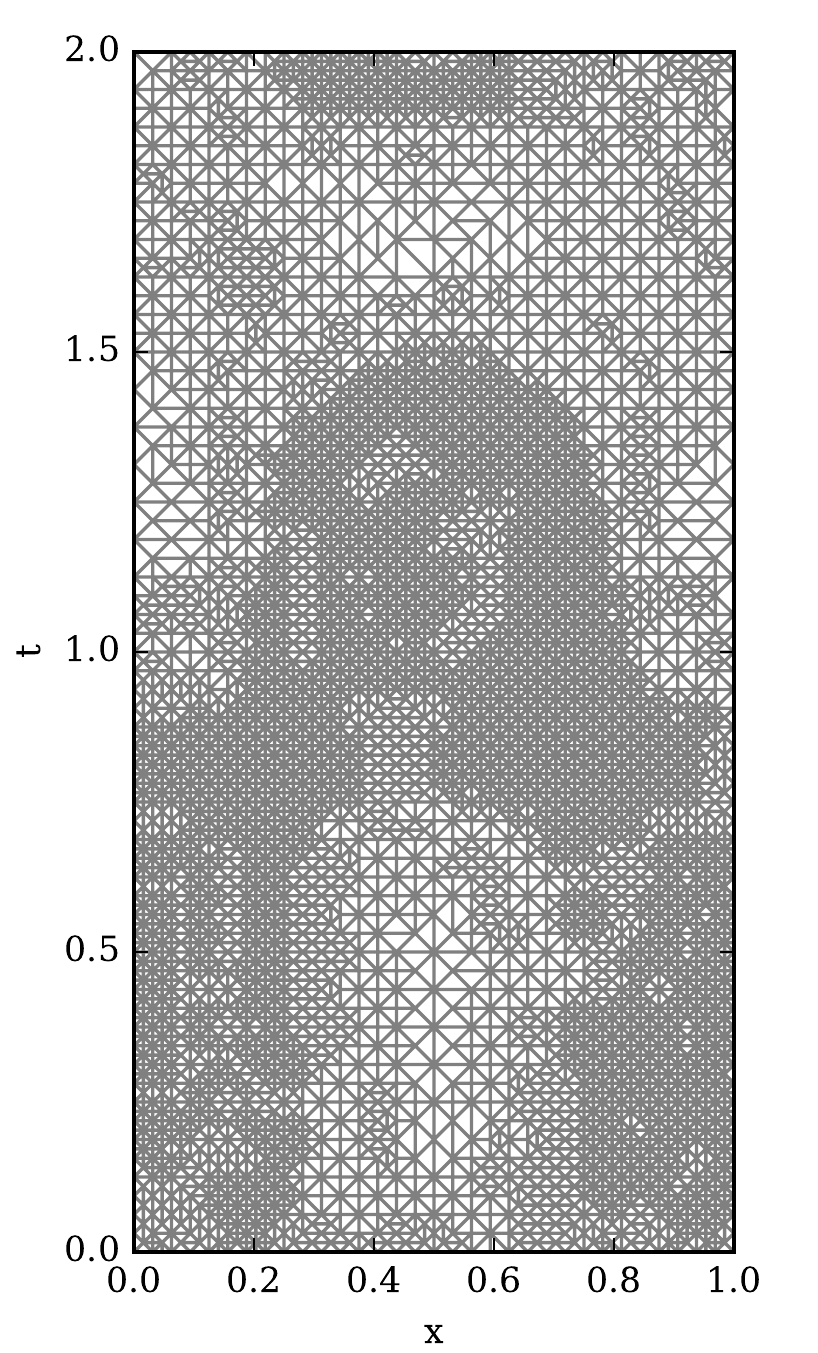}} \quad
	\subfloat[ref. based on $\overline{\rm M}_{\mu_{\rm opt}}$ for (b)]{
	\includegraphics[width=5.5cm, trim={6cm 0.5cm 6cm 1cm}, clip]{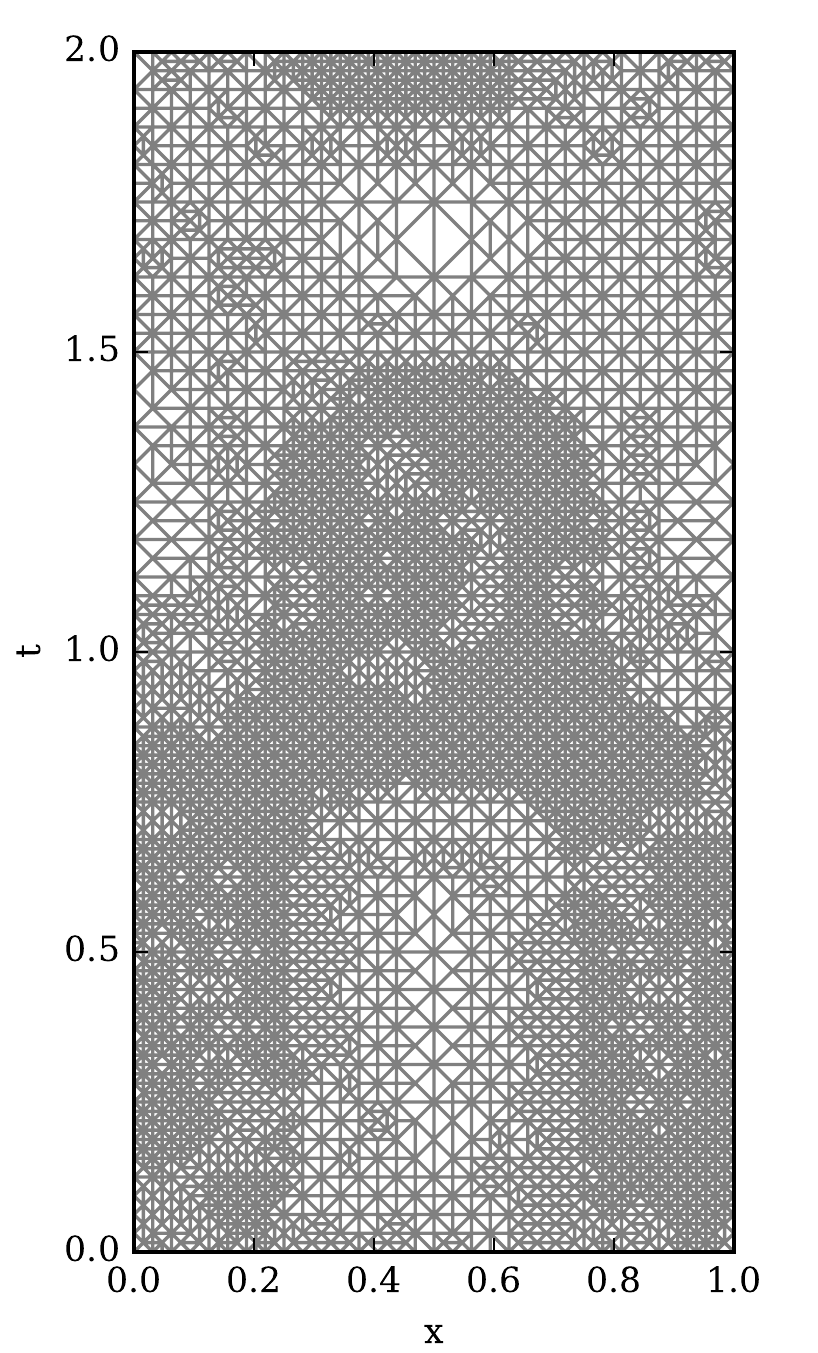}} \quad
	\subfloat[ref. based on $\overline{\rm M}_{\mu_{\rm opt}}$ for (c)]{
	\includegraphics[width=5.5cm, trim={6cm 0.5cm 6cm 1cm}, clip]{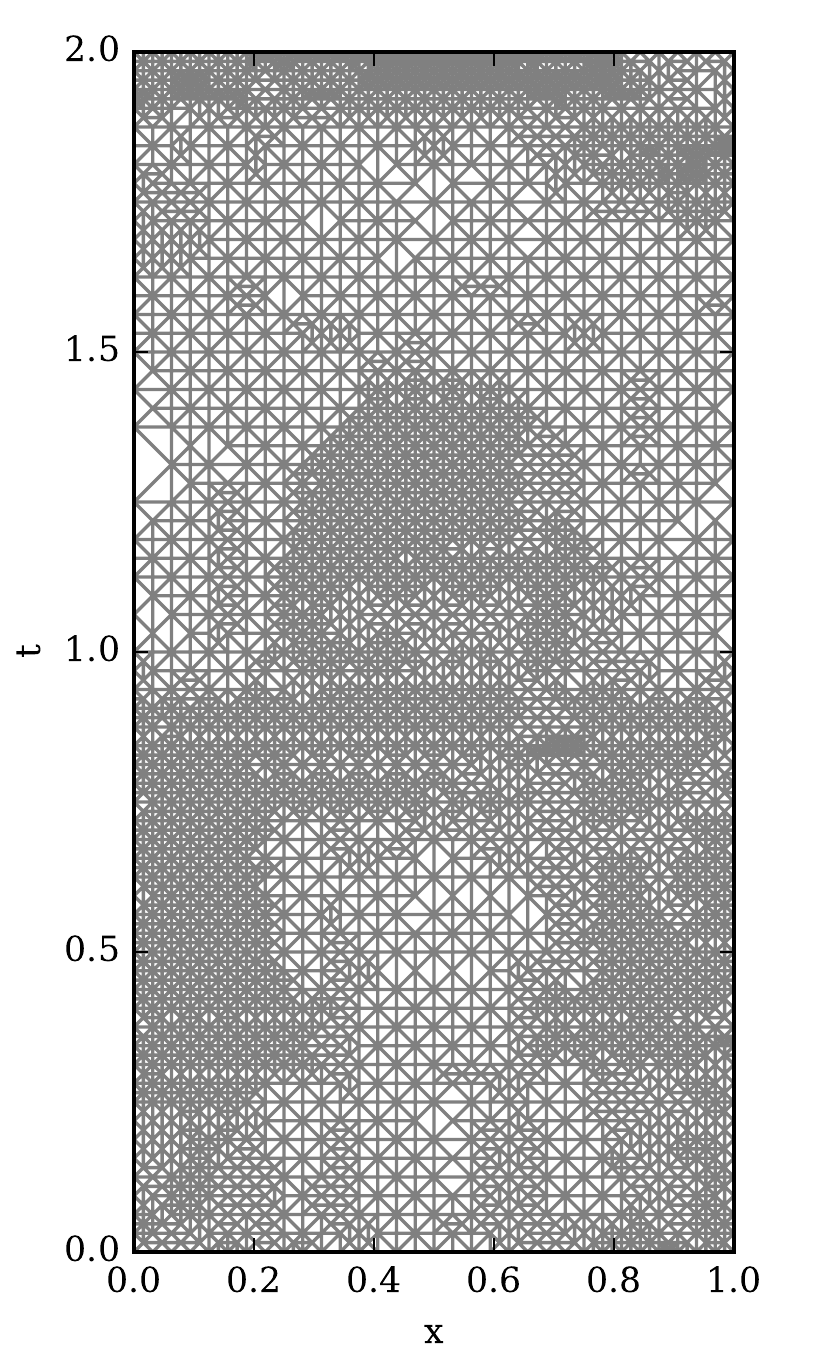}} \\
	\subfloat[ref. based on $\overline{\rm M}_{\mu = 1}$ for (a)]{
	\includegraphics[width=5.5cm, trim={6cm 0.5cm 6cm 1cm}, clip]{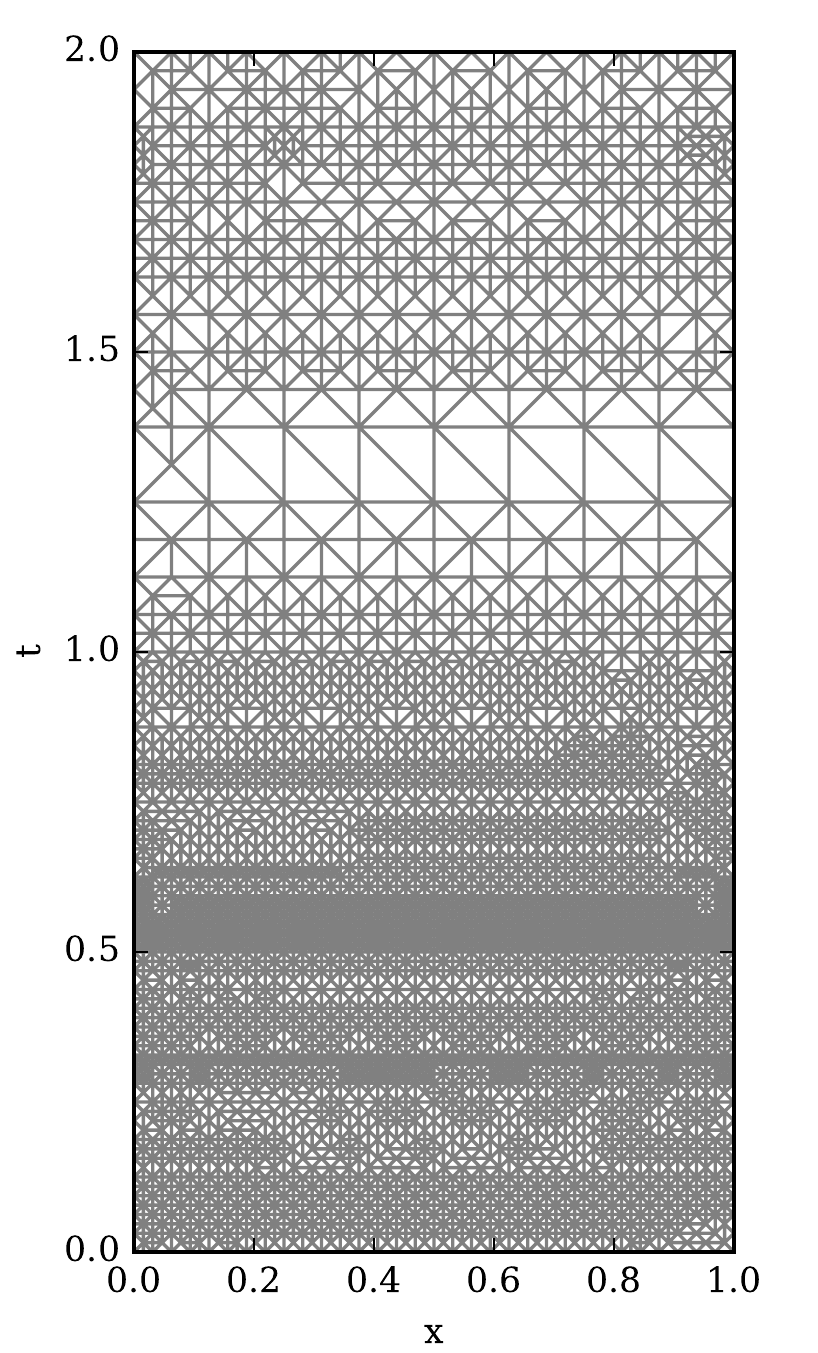}} \quad
	\subfloat[ref. based on $\overline{\rm M}_{\mu = 1}$ for (b)]{
	\includegraphics[width=5.5cm, trim={6cm 0.5cm 6cm 1cm}, clip]{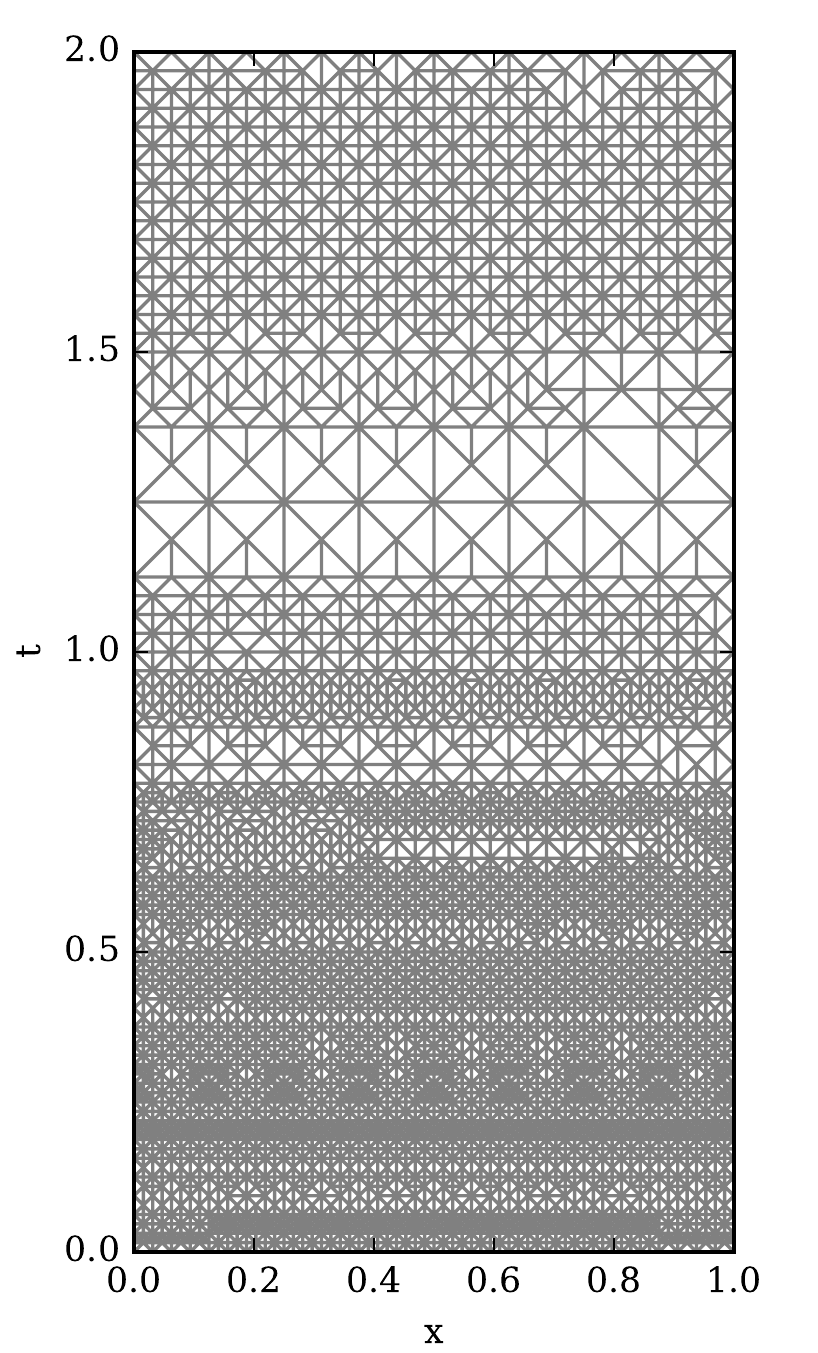}} \quad
	\subfloat[ref. based on $\overline{\rm M}_{\mu = 1}$ for (c)]{
	\includegraphics[width=5.5cm, trim={6cm 0.5cm 6cm 1cm}, clip]{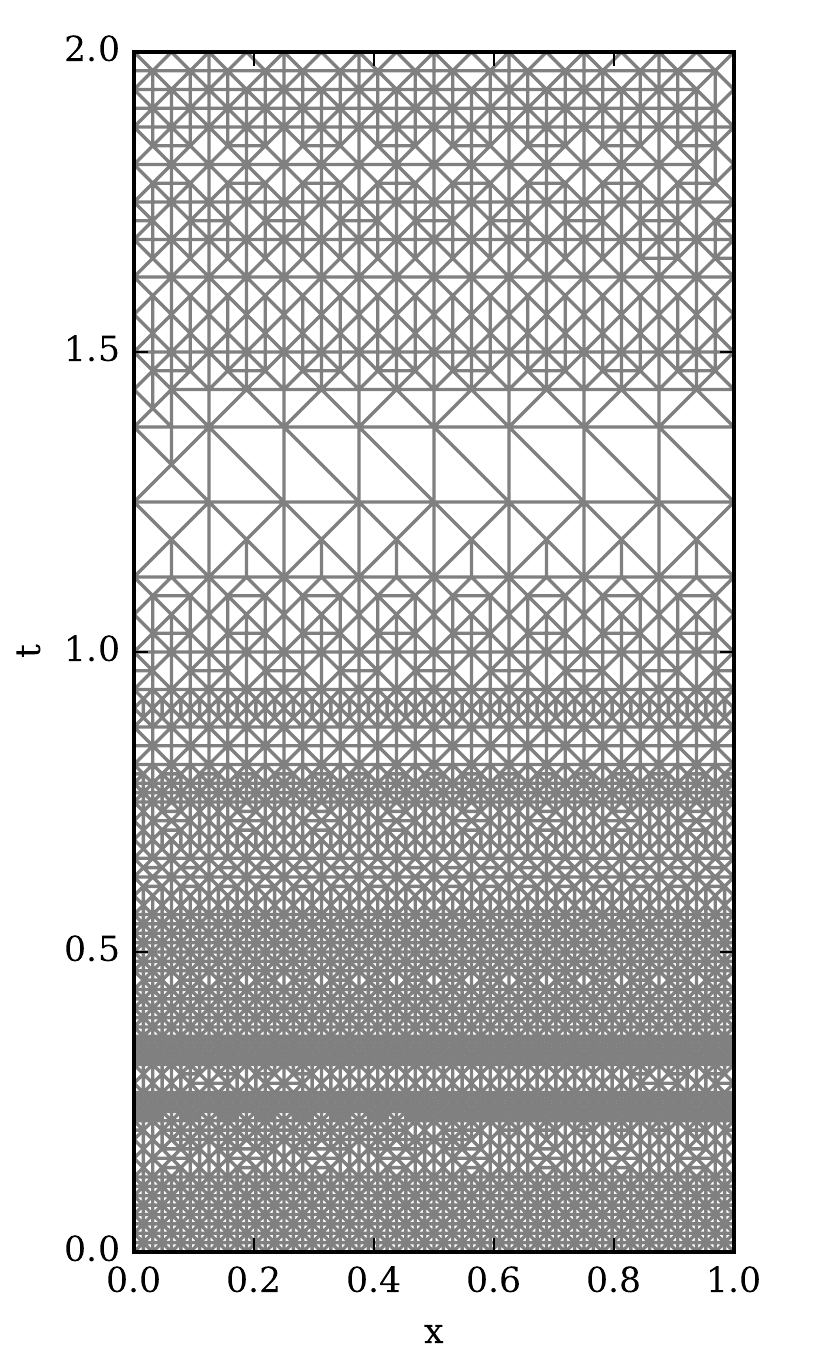}} 
\caption{Example \ref{ex:example-39-40-41}. The mesh obtained on the sixth refinement step using 
$\overline{\rm M}_{\mu_{\rm opt}}$ (first row) and $\overline{\rm M}_{\mu = 1}$ (second row). 
Here, $v \in \Pone$, $\flux \in \Ptwo$ with bulk marking $\theta = 0.6$.}
\label{fig:example-39-40-41-comparison-for-meshes-mu-opt}
\end{figure}

In order to provide robustness of the majorant w.r.t. changing values of $\lambda$, analogously to the 
elliptic case, we consider the additional function $\mu \in L^{\infty}_{[0, 1]} (Q_T)$, 
which is used 
to balance the contribution of $\lambda$ into the ${\bf r}_{\rm eq}$ term, i.e., 
\begin{alignat*}{2}
\int_{Q_T} {\bf r}_{\rm eq} \, e \dxt
= \int_{Q_T} \mu \, {\bf r}_{\rm eq} \, e \dxt 
  + \int_{Q_T} (1 - \mu) \,{\bf r}_{\rm eq}\, e  \dxt
\leq \| \tfrac{\mu}{\delta} \, {\bf r}_{\rm eq} \|_{\Omega} \| \delta e\|_{Q_T} 
+ \tfrac{\CFriedrichs}{\sqrt{\varepsilon}} \, \| (1 - \mu) \,{\bf r}_{\rm eq} \|_{Q_T} \, 
\| \nabla e \|_{\varepsilon, Q_T}.
\end{alignat*}
After setting parameters $\alpha_1 = 1 + \beta$,  
$\alpha_2 = 1 + \tfrac{1}{\beta}$, with $\beta ={\rm const}$, 
and $\alpha_3 = 2$ in \eqref{eq:estimate-space-time},
assuming that the initial conditions are satisfied exactly, i.e.,  $e(0, x) = 0$,
the final form of the estimate can be written as follows
\begin{alignat}{2}
	\overline{\rm M}^2_{\mu} (v, \flux) 
	& := \gamma \Big\| \, \tfrac{\mu}{\sqrt{\lambda}} \, {\bf r}_{\rm eq} (v, \flux) \, \Big\|^2_{Q_T} +
	                 \Big(1 + \tfrac{1}{\beta}\Big) \tfrac{\CF^2}{\varepsilon} 
	                 \| \, (1 - \mu) \, {\bf r}_{\rm eq} (v, \flux) \, \|^2_{Q_T} +
			(1 + \beta) \| \, {\bf r}_{\rm d} (v, \flux) \,\|^2_{\varepsilon^{ \minus 1}, Q_T} 
			 + 2\, \tfrac{\Ctr^2}{\varepsilon} \, \| \flux \cdot  {\bf n} \|^2_{S_T} \nonumber\\
	& = \Int_{Q_T} \Big ( \gamma  \tfrac{\mu^2}{\lambda} \, {\bf r}^2_{\rm eq} (v, \flux) +
	                 \big(1 + \tfrac{1}{\beta}\big) \CF^2 (1 - \mu)^2 \, {\bf r}^2_{\rm eq} (v, \flux) +
		        (1 + \beta) \varepsilon^{-1} {\bf r}^2_{\rm d} (v, \flux) \Big) \dx \dt 
		         + 2\, \tfrac{\Ctr^2}{\varepsilon} \, \| \flux \cdot  {\bf n} \|^2_{S_T}.
	\label{eq:majorant-mu-simplified}
\end{alignat}
The minimisation of \eqref{eq:majorant-mu-simplified} w.r.t. $\mu$ is equivalent to the following variational 
problem: find $\hat{\mu} \in L^{\infty}_{[0, 1]} (Q_T)$ such that
\begin{equation}
\mu_{\rm opt} :={\rm arg}\!\!\!\!\!\!\!\!\!\!\min\limits_{\mu \in L^{\infty}_{[0, 1]} (Q_T)} \!\!\! \Upsilon(\mu), \quad 
\Upsilon(\mu) := \IntQT \Big (\gamma  \tfrac{\mu^2}{\lambda} \, {\bf r}^2_{\rm eq} (v, \flux) \, \: +
	                 \big(1 + \tfrac{1}{\beta}\big) \CF^2 (1 - \mu)^2 \, {\bf r}^2_{\rm eq} (v, \flux) \Big ) 
	                 \dx \dt.
\end{equation}
Functional $\Upsilon(\mu)$ admits its minimum in
$\mu_{\rm opt} (x, t) 
= \frac{\CF^2(1 + \beta)\lambda}{\beta \gamma + \CF^2(1 + \beta)\lambda}.$
Table \ref{tab:example-39-40-41-comparison-for-mu-opt} provides the efficiency indices of  
${\overline{\rm M}}_{\mu_{\rm opt}}$ and ${\overline{\rm M}}_{\mu = 1}$ for different cases (a)--(b). 
The results correspond to the approximation $v \in \Pone$ and the auxiliary flux $\flux \in \Ptwo$. 
The majorant
${\overline{\rm M}}_{\mu = 1}$ grows dramatically in the case (a) and diverges 
according to the efficiency indices for the cases (b) and (c). At the same time, 
${\overline{\rm M}}_{\mu_{\rm opt}}$ performs efficiently for all considered $\lambda$. 
It stays efficient and robust
even if the reaction function changes its values drastically in different parts of the domain. The 
advantage of using an error estimate with optimal weighting function $\mu$ is illustrated also by 
the meshes in Figure \ref{fig:example-39-40-41-comparison-for-meshes-mu-opt}. 
Using 
${\overline{\rm M}}_{\mu = 1}$ in all cases (a)--(c) results in a refinement procedure generating 
meshes that do not make changes in the solution visible (see the lower row of 
Figure \ref{fig:example-39-40-41-comparison-for-meshes-mu-opt}). The meshes presented in the upper 
row, however, reflect the approximate solution accurately.

If we assume that $\lambda = 0$ and ${\bf a = 0}$ (${\bf F = 0}$ in the problem statement 
\eqref{eq:fokker-planck-equation}--\eqref{eq:parabolic-robin-bc}), 
the test case reduces to the diffusion equation of parabolic type with right-hand side
$f = u_0 \, (\cos 2 t - 2 t \sin 2 t ) + \pi^2 \, \cos \pi x  \, \left(t \, \cos 2 \, t  + 1\right).$
Assuming that $f \in L_{2}(Q)$ and $u_0 \in \HD{1}{0}(\Sigma_0)$ such a problem is uniquely 
solvable in $V^{\Delta_x}_{0}$, where
$$V^{\Delta_x}_{0} := 
\big\{ v \in H^{1, 1}(Q_T) \, | \, 
\Delta_x v \in L^{2}(Q_T) \; \wedge \;  v |_{\Omega} = 0 \, \text{ for a.a. } t \in \, (0, T) \big\},$$
and $u$ depends continuously on $t$ in the $\HD{1}{0}$ norm. If $u, v \in V^{\Delta_x}_0$ 
are provided, the time-dependent diffusion equation \eqref{eq:fokker-planck-equation} imposes the error 
identity (see \cite{LMR:AnjamPauly:2016}):
\begin{alignat}{2}
\| \dvrg( \varepsilon \nabla_x e) \|^2_{Q_T} 
+  \| \partial_t e \|^2_{Q_T} 
+  \| \nabla_x e \|^2_{\varepsilon, \Sigma_T} 
=: |\!|\!|  e |\!|\!|^2_{\mathcal{L}} 
\equiv \EI^2 (v)
:= \| \nabla_x (u_0 - v) \|^2_{\varepsilon, \Sigma_0} 
+ \|  \dvrg( \varepsilon \nabla_x v) + f - \partial_t v\|^2_{Q_T}.
\label{eq:strong-norm-error-identity}
\end{alignat}

The results corresponding to the performance of the majorant $\overline{\rm M}_{\mu_{\rm opt}}$ 
as well as $\EI$ in comparison to the errors they bound are illustrated in Table 
\ref{tab:example-39-40-41-mu-opt-error-id}. We see that for the setting (I) we get very close to the value one 
for the efficiency indices of $\overline{\rm M}_{\mu_{\rm opt}}$, and $\EI$ is always identical to the values
of $|\!|\!|  e |\!|\!|^2_{\mathcal{L}}$. However, an approximation of $v$ by linear affine elements does not 
provide the condition $v \in V^{\Delta_x}_0$ required for application of $\EI$. Therefore, stagnation of 
$|\!|\!|  e |\!|\!|^2_{\mathcal{L}}$, as well as $\EI$, is expected. For the setting (I\!I), we see that the 
strong norm of the error and the error identity begin to converge as $O(h)$. We also 
illustrate the results for the setting (I\!I\!I) in order to show that the performance of 
$\overline{\rm M}_{\mu_{\rm opt}}$ improves for smoother approximations of $\flux$ (see the $5$-th 
column). As we have observed in the numerical experiments,
considering even smoother approximations of the flux 
does not pay off in this case, i.e., the minimisation of the majorant becomes even more time consuming, 
whereas the efficiency indices improve only slightly. The error identity stays quantitatively identical to the error norm 
$|\!|\!|  e |\!|\!|^2_{\mathcal{L}}$ and, since it is only dependent on the approximation $v$, it does not 
require any computational time overhead. At the same time, the restriction on the regularity of
$u, v \in V^{\Delta_x}_0$ still remains, which makes the application of the error identity rather restricted 
to a certain class of problems.

One of the ways to study the efficiency of the majorant when it comes to generating the mesh refinement procedure 
is to compare the numerical results for the refinement steps of 
adaptive meshes, where the refinement 
is based on the element-wise distribution of  
the majorant, with those meshes refined 
w.r.t. 
the the element-wise errors.
Figure \ref{fig:example-39-40-41-comparison-for-meshes-maj-error} provides such a comparison of 
meshes obtained on the refinement steps 4, 6, and 8. 
We see that the meshes obtained in the upper row are topologically close
to the meshes in the second row.

\begin{table}[!ht]
\centering
\footnotesize
\newcolumntype{g}{>{\columncolor{gainsboro}}c} 	
\newcolumntype{k}{>{\columncolor{lightgray}}c} 	
\newcolumntype{s}{>{\columncolor{silver}}c} 
\newcolumntype{a}{>{\columncolor{ashgrey}}c}
\newcolumntype{b}{>{\columncolor{battleshipgrey}}c}
\begin{tabular}{c|cccg|ccg}
\midrule
ref.   & 
d.o.f. $\#$& $[e]$                 & $\overline{\rm M}_{\mu_{\rm opt}}$ & $\Ieff (\overline{\rm M}_{\mu_{\rm opt}})$ & 
$|\!|\!| e |\!|\!|$ & $\EI$                                                   & $\Ieff (\EI)$ \\
\midrule
\multicolumn{8}{c}{(I) $v \in \Pone$ and $\flux \in \Ptwo$ } \\
\midrule
      1 &      272 &  2.54e-01 & 3.84e-01 &       1.51 &  8.42e+00 & 8.42e+00 &       1.00 \\
      3 &      887 &  1.36e-01 & 1.99e-01 &       1.46 &  8.42e+00 & 8.42e+00 &       1.00 \\
      5 &     2867 &  7.42e-02 & 1.08e-01 &       1.45 &  8.42e+00 & 8.42e+00 &       1.00 \\
      7 &     8738 &  4.08e-02 & 5.99e-02 &       1.47 &  8.42e+00 & 8.42e+00 &       1.00 \\
      9 &    26132 &  2.33e-02 & 3.53e-02 &       1.52 &  8.42e+00 & 8.42e+00 &       1.00 \\
\midrule
\multicolumn{8}{c}{(I\!I) $v \in \Ptwo$ and $\flux \in \Ptwo$} \\
\midrule
      1 &     1033 &  1.44e-02 & 7.60e-02 &       5.26 &  8.08e-01 & 8.08e-01 &       1.00 \\
      3 &     3672 &  4.09e-03 & 1.83e-02 &       4.49 &  4.28e-01 & 4.28e-01 &       1.00 \\
      5 &    12382 &  1.14e-03 & 5.35e-03 &       4.70 &  2.32e-01 & 2.32e-01 &       1.00 \\
      7 &    40392 &  3.28e-04 & 1.53e-03 &       4.67 &  1.27e-01 & 1.27e-01 &       1.00 \\
      9 &   132711 &  9.92e-05 & 4.74e-04 &       4.78 &  6.85e-02 & 6.85e-02 &       1.00 \\
\midrule
\multicolumn{8}{c}{(I\!I\!I) $v \in \Ptwo$ and $\flux \in \Pthree$} \\
\midrule
      1 &     1029 &  1.53e-02 & 4.73e-02 &       3.10 &  8.46e-01 & 8.45e-01 &       1.00 \\
      3 &     3789 &  3.98e-03 & 1.09e-02 &       2.73 &  4.25e-01 & 4.24e-01 &       1.00 \\
      5 &    12703 &  1.10e-03 & 3.33e-03 &       3.01 &  2.26e-01 & 2.26e-01 &       1.00 \\
      7 &    41749 &  3.17e-04 & 1.01e-03 &       3.20 &  1.23e-01 & 1.23e-01 &       1.00 \\
      9 &   141873 &  9.40e-05 & 2.92e-04 &       3.10 &  6.63e-02 & 6.63e-02 &       1.00 \\
\end{tabular}
\caption{Example \ref{ex:example-39-40-41}. Efficiency of $\overline{\rm M}_{\mu_{\rm opt}}$
and $\EI$ with bulk marking $\theta = 0.6$.}
\label{tab:example-39-40-41-mu-opt-error-id}
\end{table}

\begin{figure}[!t]
	\centering
	\captionsetup[subfigure]{oneside, labelformat=empty}
	\subfloat[ref. 4: based on $\overline{\rm M}_{\mu_{\rm opt}}$]{
	\includegraphics[width=5.5cm, trim={6cm 0.5cm 6cm 1cm}, clip]{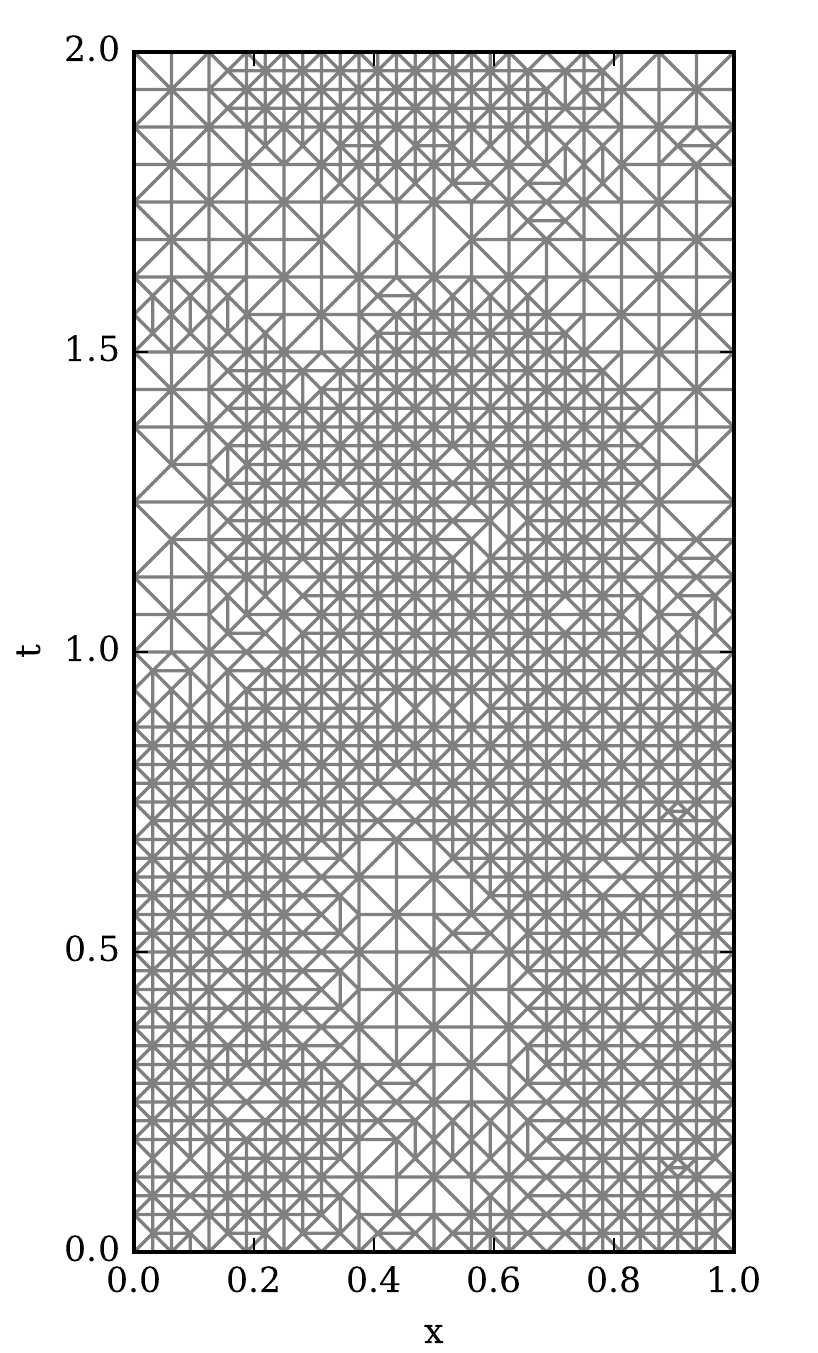}} \quad
	\subfloat[ref. 6: based on $\overline{\rm M}_{\mu_{\rm opt}}$]{
	\includegraphics[width=5.5cm, trim={6cm 0.5cm 6cm 1cm}, clip]{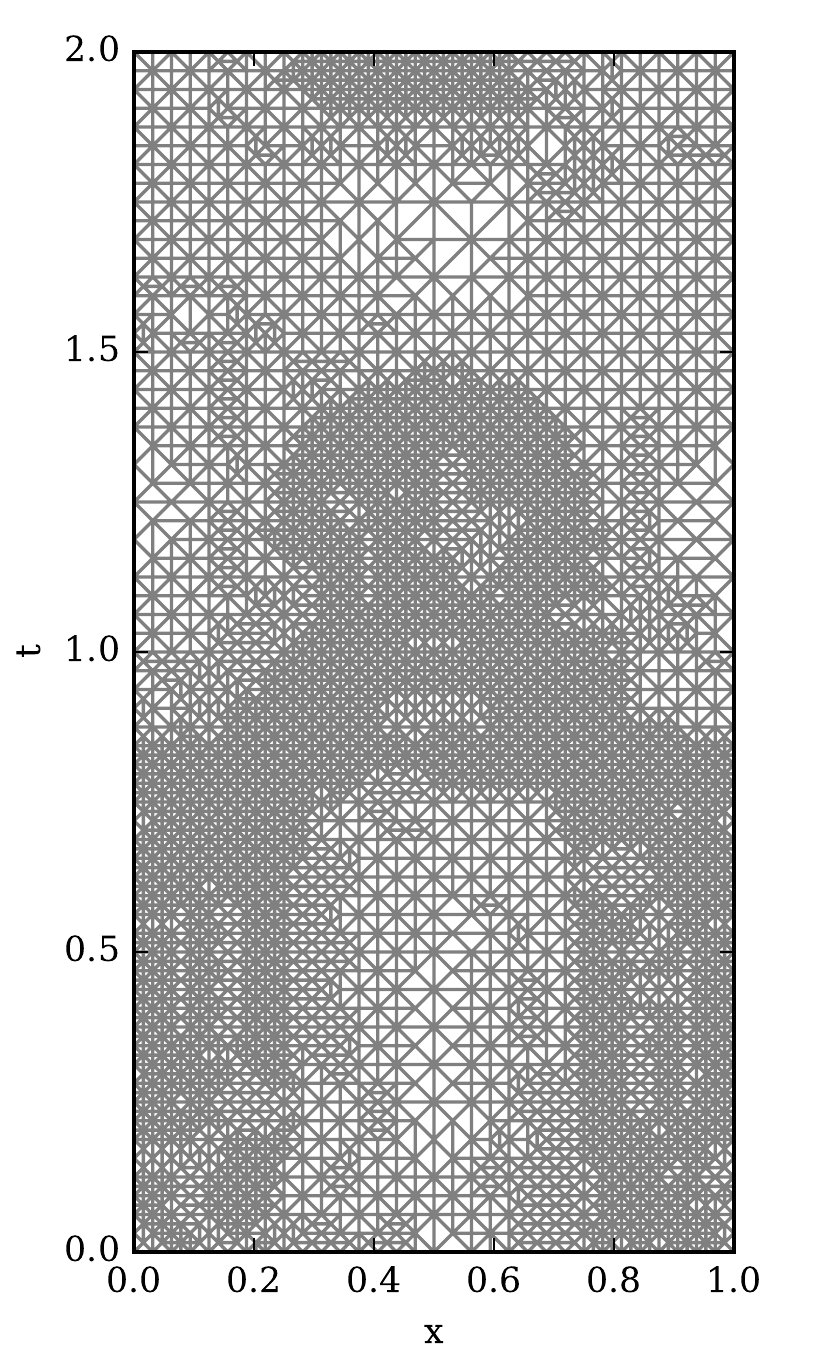}} \quad
	\subfloat[ref. 8: based on $\overline{\rm M}_{\mu_{\rm opt}}$]{
	\includegraphics[width=5.5cm, trim={6cm 0.5cm 6cm 1cm}, clip]{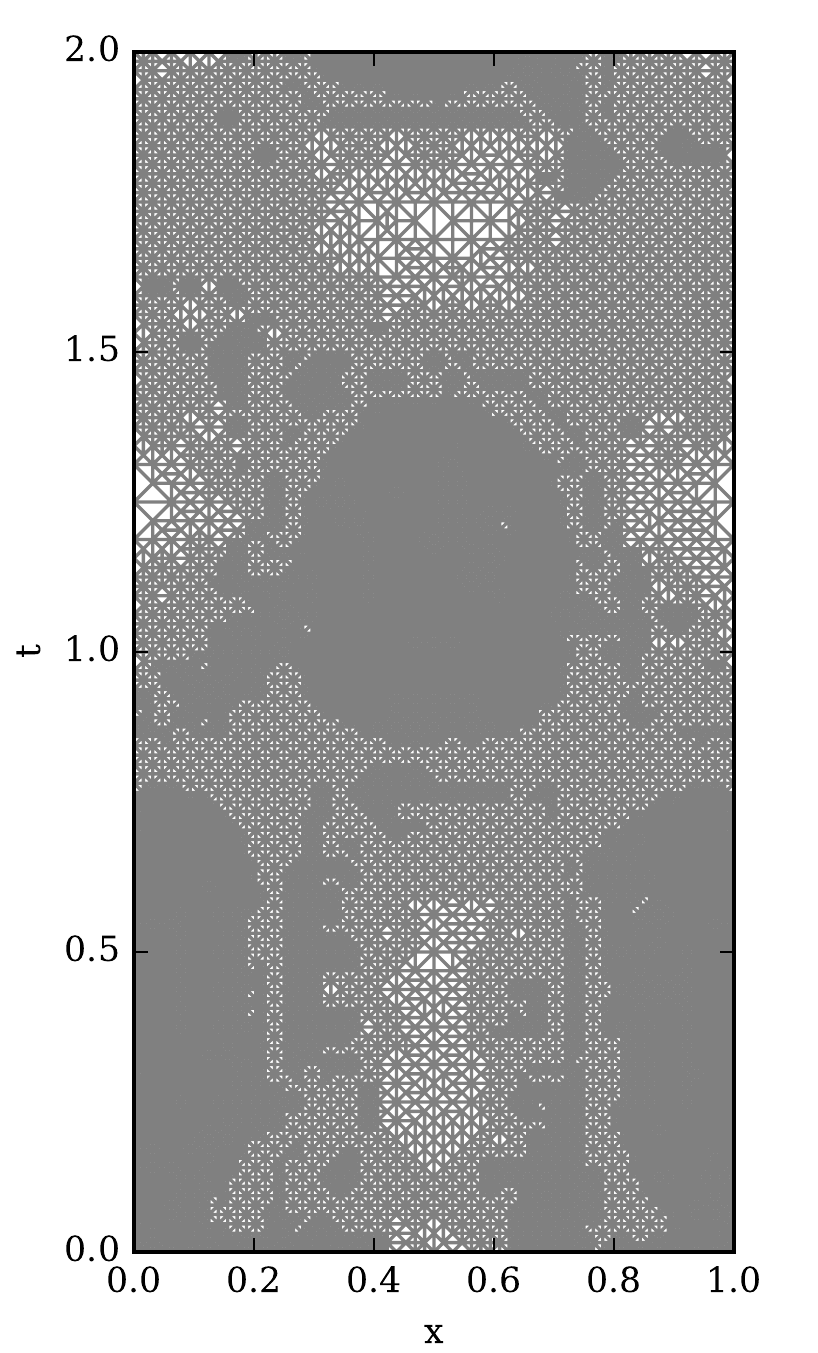}} \\
	\subfloat[ref. 4: based on ${[e]}$]{
	\includegraphics[width=5.5cm, trim={6cm 0.5cm 6cm 1cm}, clip]{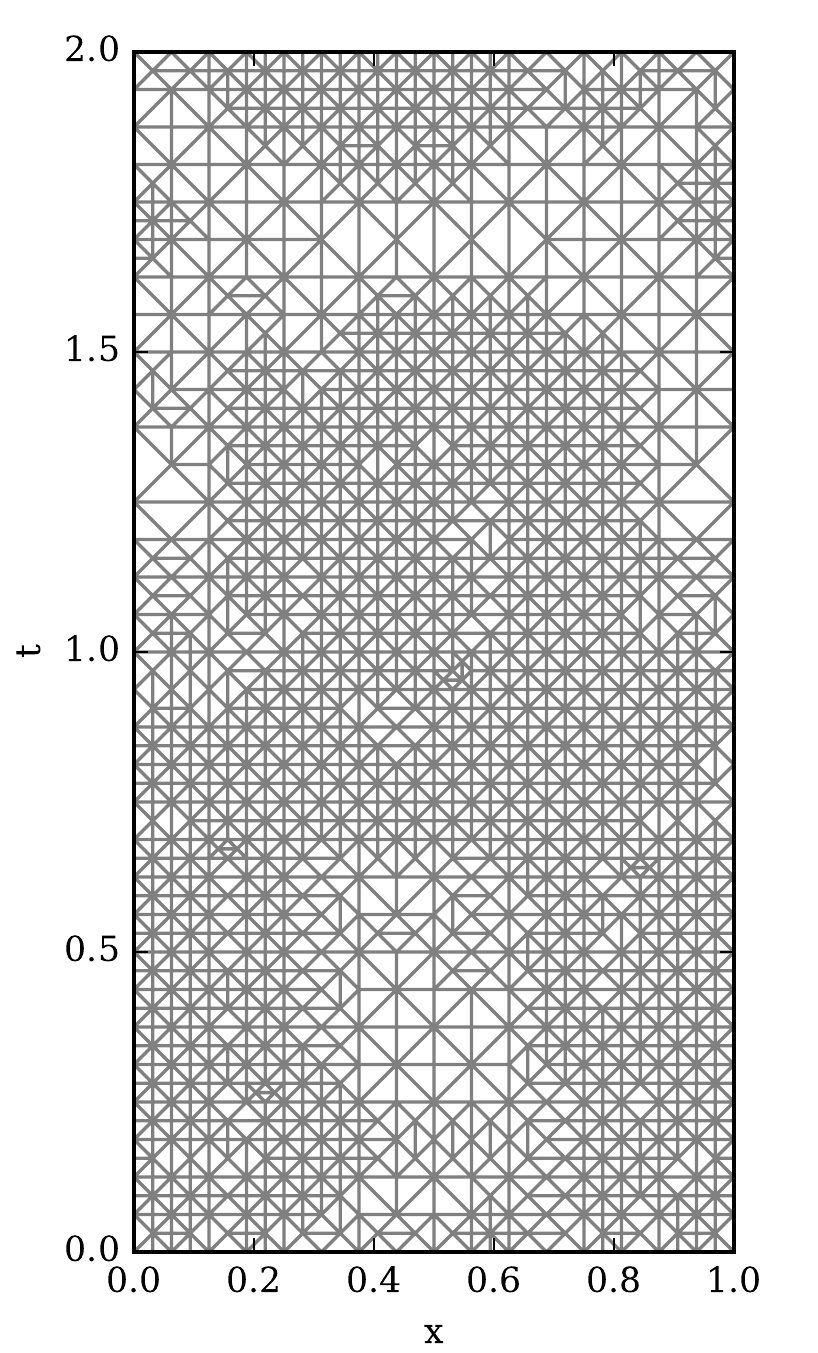}} \quad
	\subfloat[ref. 6: based on ${[e]}$]{
	\includegraphics[width=5.5cm, trim={6cm 0.5cm 6cm 1cm}, clip]{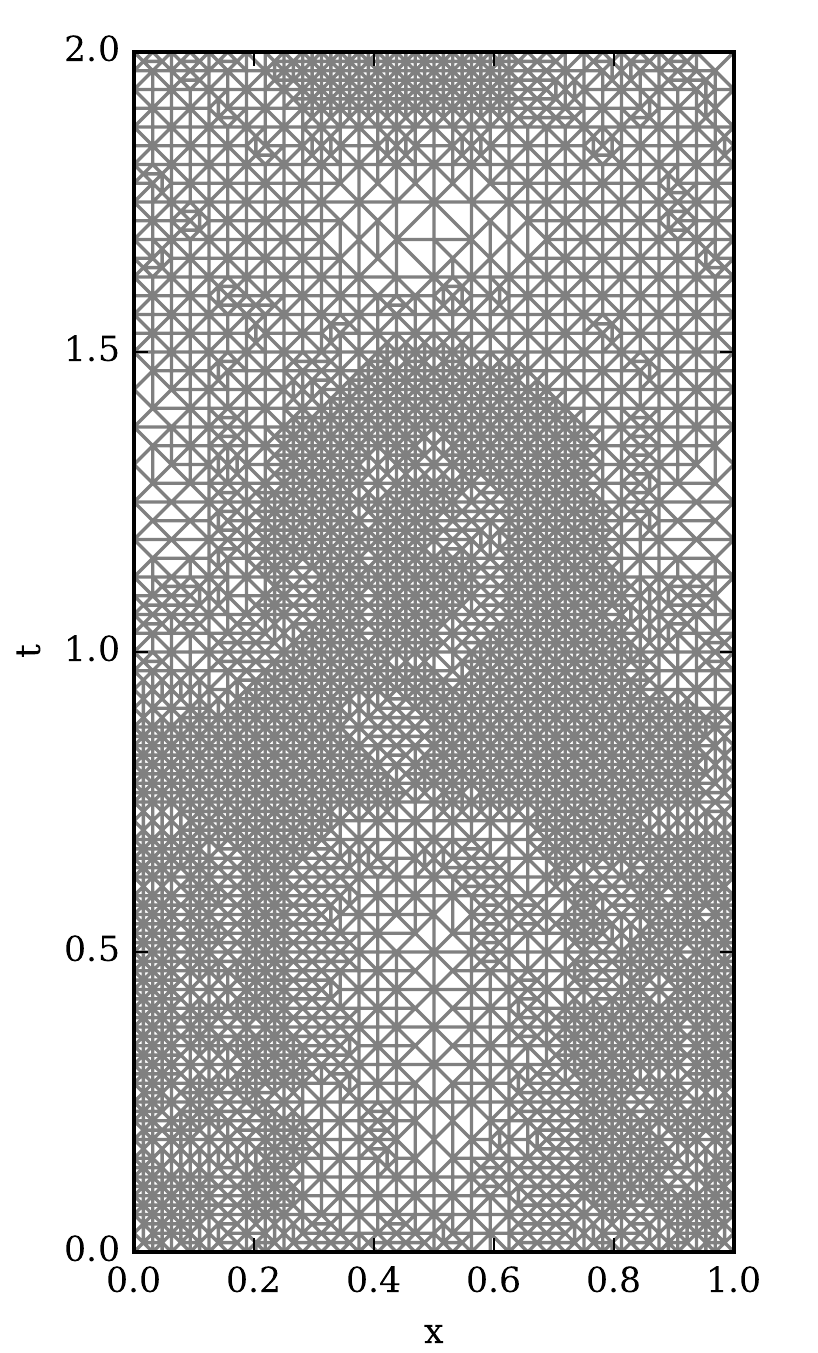}} \quad
	\subfloat[ref. 8: based on ${[e]}$]{
	\includegraphics[width=5.5cm, trim={6cm 0.5cm 6cm 1cm}, clip]{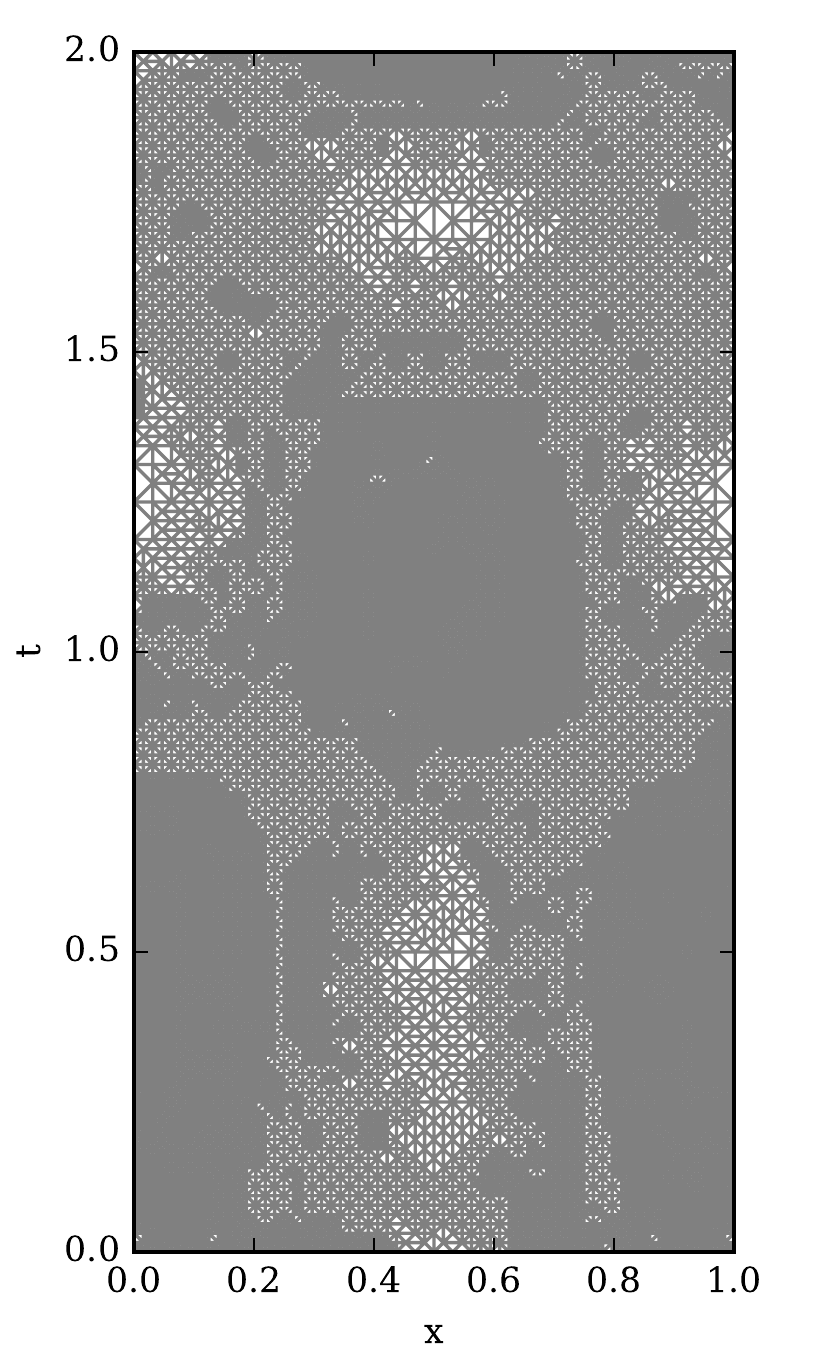}} 
\caption{Example \ref{ex:example-39-40-41}. The mesh obtained on the refinement steps (ref.) 4, 6, and 8 using 
$\overline{\rm M}_{\mu_{\rm opt}}$ (first row) and ${[e]}$ (second row). 
Here, $v \in \Pone$, $\flux \in \Ptwo$ with bulk marking $\theta = 0.6$.}
\label{fig:example-39-40-41-comparison-for-meshes-maj-error}
\end{figure}

\end{example}

\section{Conclusions}

We have studied functional type a posteriori error estimation for stationary and time-dependent linear
convection-diffusion problems. This work was motivated by a decision-making Fokker-Planck model 
problem appearing in computational neuroscience, which has been discussed in 
\cite{CarrilloCordierMancini2011}. We derived guaranteed two-sided estimates for the elliptic as well as 
parabolic problem and extended the a posteriori error analysis earlier derived in 
\cite{RepinDeGruyter2008, RepinTomar2010} by studying the numerical aspects as well as properties 
of the two-sided error bounds applied to this class of problems.

We presented a set of various numerical examples for different parameters and boundary conditions
showing the sharpness of the upper and lower bounds in practice. In particular, we addressed 
convection-dominated problems including adaptive mesh refinement based on a bulk marking criterion 
from \cite{Doerfler1996}. Furthermore, we presented numerical experiments applying the 
multilevel-homotopic-adaptive finite element method in the stationary case and space-time discretisation 
for the time-dependent reaction-convection-diffusion problem. The numerical results in both cases, 
static as well as time-dependent, successfully demonstrated the capability of the two-sided error bounds 
to provide a basis for the implementation of reliable and efficient computational methods.

\section*{Acknowledgments}

The authors gratefully acknowledge the financial support by 
the Austrian Science Fund (FWF) through the NFN S117-03 
project, and by the Academy of Finland, grant 295897.

\bibliographystyle{plain}
\bibliography{lib,my_lib}

\end{document}